\documentclass[
        parskip = half, %
        twoside,
        footinclude = false, %
        ]{scrartcl}

\usepackage{ifthen} %
\usepackage[normalem]{ulem}
\newboolean{boldeqnitalic}
\setboolean{boldeqnitalic}{true}
\usepackage{defdb1} %
\usepackage{rublkm_paper} %
\usepackage[export]{adjustbox}
\usepackage{booktabs}
\usepackage{subfig}
\usepackage[font=footnotesize,labelfont=bf]{caption}
\usepackage{xurl}                %
\usepackage[title]{appendix}     %
\usepackage{etoolbox}
\AtBeginEnvironment{appendices}{%
  \counterwithin{figure}{section}\counterwithin{table}{section}%
  }
\AtBeginDocument{\let\cite\citep}
\makeatletter
\newcommand*{\authormark}{}
\newcommand*{\markauthor}[1]{%
  \renewcommand{\authormark}{#1}%
  \ignorespaces
}
\newcommand*{\titlemark}{}
\newcommand*{\marktitle}[1]{%
  \renewcommand{\titlemark}{#1}%
  \ignorespaces
}
\makeatother

\markauthor{Sharan N. Ramesh, Lea Sa\ss mannshausen, Axel Klawonn, Daniel Balzani}
\marktitle{Pharmaco-mechanical simulation of atherosclerotic arteries}

\clearpairofpagestyles
\rohead{\pagemark}
\rehead{\pagemark}

\lohead{\titlemark}
\lehead{\authormark}

\KOMAoptions{
   headsepline = true
}
\setkomafont{pagehead}{\normalfont\normalcolor\small}

\newcommand{\colsr}[1]{#1}

\begin{document}

\clearpage
\setcounter{page}{1}
\begin{center}
    {\LARGE Pharmaco-mechanical computational simulation of atherosclerotic arteries: effects of calcium channel blockers on local stress distributions}

    \vspace{5mm}

    S. Nurani Ramesh$^{1}$, L. Sa\ss mannshausen$^{2,3}$, A. Klawonn$^{2,3}$, D. Balzani$^{1\star}$

    \vspace{3mm}

    {\small $^1$Chair of Continuum Mechanics, Ruhr-Universit\"at Bochum,\\
    Universit\"atsstra{\ss}e~150, 44801~Bochum, Germany}\\[3mm]
    {\small $^2$Department of Mathematics and Computer Science, University of Cologne,\\
    Weyertal 86-90, 50931 Cologne, Germany}\\[3mm]
    {\small $^3$Center for Data and Simulation Science, University of Cologne,\\
    Albertus-Magnus-Platz, 50923 Cologne, Germany}\\[3mm]

    \vspace{3mm}

    {\small ${}^{\star}$E-mail address of corresponding author: daniel.balzani@rub.de}

    \vspace{10mm}

    \begin{minipage}{15.0cm}
    \textbf{Abstract}\hspace{3mm}
    This paper presents a computational model integrating transmural drug transport, smooth muscle cell contraction, nonlinear anisotropic material behavior, and residual stresses to investigate the effects of calcium channel blockers~(CCBs) on local stress distributions in atherosclerotic arteries.
    The balance of linear momentum is coupled with a three-dimensional reaction-diffusion equation for the drug transport.
    Smooth muscle contraction is modeled phenomenologically based on a chemo-mechanical model that incorporates the effects of stretch on intracellular calcium dynamics and calcium sensitization.
    A novel drug-tissue interaction model is proposed to capture the effects of CCBs on intracellular calcium and thus on smooth muscle contraction.
    Material parameters are identified by fitting to experimental data from the literature and finite element simulations are conducted on five exemplary atherosclerotic arterial geometries reconstructed from histological cross-sections.
    Simulation results show that although CCBs reduce the stresses in the media through reduced smooth muscle contraction, they shift load to the fibrous cap, whose stress increases in all five geometries (median of up to 28.9~kPa).
    This increase varies strongly between plaque morphologies and becomes smaller with lower smooth muscle tone before treatment, lower drug concentrations in the arterial wall and increasing axial pre-stretch.
    A blood pressure reduction of 10~mmHg already lowers it to a median of at most 3.8~kPa in four of the five geometries.
    Our findings suggest that CCBs may increase the rupture risk of vulnerable plaques with highly loaded fibrous caps under insufficient blood pressure reduction, which remains to be confirmed experimentally and clinically.
    \end{minipage}
\end{center}

    \medskip{}

\section{Introduction}

    Cardiovascular diseases (CVDs) remain the foremost cause of mortality globally, claiming an estimated 19.8 million lives in 2022 and representing approximately 32\% of all deaths worldwide.
    Over three-quarters of these fatalities occurred in low- and middle-income countries, disproportionately affecting the less fortunate~\cite{Who:2024:ghe}.
    Among the modifiable risk factors for CVD, hypertension is the most significant, affecting 1.28 billion adults aged 30--79 years and serving as a primary driver of premature death~\cite{Who:2021:hyp,MilBunKelReeKeaReyCheHe:2016:gdo}.
    In response to this burden, there has been a steady increase in the adoption of and treatment with antihypertensive drugs.
    This is evidenced by a rise in the global median consumption of antihypertensive drugs from 184.78 standard units per 1,000 inhabitants per day in 2010 to 325.6 in 2021~\cite{Oec:2023:hag}.
    This, alongside initiatives like the WHO Global Hearts program that placed 7.5 million individuals on protocol-based hypertension management by 2023~\cite{Who:2023:htp}, has helped significantly reduce mortality rates.
    In view of this evolving therapeutic landscape, there is
    a strong interest in
    sophisticated simulations to evaluate the consequences of antihypertensive interventions on intramural stress distributions, particularly in the context of atherosclerosis.
    This can facilitate the development of more precise and effective strategies to counteract the progression of CVDs.

    Recent advancements in computational models, along with hardware and software frameworks, have enabled the efficient simulation of large problems involving intricate bio-chemo-mechanical phenomena.
    Here, we endeavor to capitalize on these advancements to simulate the effects of antihypertensive drugs, specifically amlodipine, on the local stress distribution in atherosclerotic arteries in a finite element setting.
    Several studies have attempted to leverage finite element analysis to predict stress distribution in atherosclerotic arteries with the aim of assessing plaque rupture risk.
    \citet{MilLemDetBurMarReiBur:2021:cpc} analyzed plaque rupture risk using patient-specific atherosclerotic geometries derived from optical coherence tomography~(OCT) and found that stress concentrations were much higher in plaques that ruptured later and concluded that this could predict rupture.
    They considered the artery and plaque components to be isotropic and linear elastic.
    \citet{DorOtsNadVilKarZanDijZijSoeDae:2020:bsp} used linear elastic, transversely isotropic material models to obtain stress distributions and attempted to identify a multifactorial plaque rupture risk metric that combines biomechanical and morphological features.
    \citet{LiWanWuLiuHuaHeLiuRen:2020:stf} conducted one way coupled FSI simulations on carotid artery plaques constructed using contrast-enhanced ultrasound images. They used an isotropic third order Yeoh material model for the plaque components and found that maximum principal stresses increased with increasing stenosis and decreasing fibrous cap thickness.
    \citet{CorDeGroDinAikWeiCar:2022:teo} simulated plaque stresses from microCT-derived human coronary atheroma geometries using the nonlinear anisotropic Holzapfel-Gasser-Ogden material model and found that a reduction in plaque thickness and the presence of micro-calcifications significantly increased the maximum principal stresses and lead to a vulnerable plaque.
    Recently, \citet{ZhaLvrMaeWanGaoXuGuoZhuHuaZha:2024:pra} used OCT rupture data and fluid-structure interaction (FSI) simulations, considering the arterial wall and plaque components to be hyperelastic and anisotropic to demonstrate a correlation between high plaque stresses and rupture.
    A significant drawback in many of these studies is the use of simplified, isotropic material models.
    Uniaxial tensile test studies on atherosclerotic artery tissues have shown that they exhibit not only anisotropic characteristics but also significant nonlinear behavior, even at low strains~\cite{HolSomReg:2004:amp,YanBacZheNaqWooTenBilTan:2009:ivi,TenZhaHuaFenYuaLuSutBroJinGil:2014:mpc}.
    \citet{TenYuaFenZhaBroWanLuGil:2015:tic} investigated the suitability of various isotropic hyperelastic constitutive laws in describing the mechanical response of atherosclerotic tissues and found that apart from the neo-Hookean model, all other models were able to accurately fit experimental data.
    Furthermore, \citet{YanBacZheNaqWooTenBilTan:2009:ivi} compared the stresses obtained using isotropic and anisotropic Mooney-Rivlin models with material parameters obtained by fitting to experimental uniaxial tensile test data and found that stress values from the anisotropic models were substantially higher than the isotropic model.
    This was also found in numerical studies of overstretched arteries as occurring during angioplasty, cf.~\citet{SchPanBal:2015:ioi}.
    Taking into account more sophisticated anisotropic tissue models, \citet{BalSchOrt:2017:mft} proposed an uncertainty quantification framework for the calculation of sharpest bounds possible on the rupture probability.
    Whereas significant improvements have been achieved regarding the incorporation of the blood flow, see, e.g., \citet{BalDepFauForHeiKlaQuaRheSch:2015:nmo,BalHeiKlaRheSch:2023:coa,WanBalVedUhlVar:2021:otp}, and the references therein, %
    in the context of analyzing stress concentrations in atherosclerotic arteries, still further limitations in the wall models exist.
    These include the
    active contraction due to smooth muscle cell activity as well as the presence of residual stresses resulting from adaptive processes in the wall.

    In the simulations presented here, the material response of the arterial wall forms the basis for drug-tissue interaction.
    This not only entails capturing aspects like the nonlinear hyperelastic response of the extracellular matrix and the collagen fibers, but also the active contraction of smooth muscle cells (SMCs).
    Contraction in SMCs is a result of relative sliding between actin and myosin filaments, enabled by cross-bridges.
    Cross-bridges are the temporary links between the head regions of myosin filaments with actin binding sites.
    Several approaches have been developed to model the active contraction in the arterial wall.
    \citet{CocEdwAggNitPar:2018:ama} presented a biochemical model based on calcium dynamics advanced by \citet{BoiGeoParMitAziNit:2015:sbi}.
    Notably, \citet{HerOdwYanMatTieFonLewSanCla:2024:acm} developed a sex-specific mathematical model of electrical and calcium ion signaling in mesenteric arterial myocytes and extended the framework to include a 1D tissue representation, predicting heightened female sensitivity to L-type calcium channel blockers due to differential voltage-gated potassium channel expression.
    Some others have explored multiscale models~(see \citet{DalFanHol:2024:mcm} for a comprehensive review).
    In contrast, \citet{IroHum:2020:csm} presented a graph based continuum model for cell signaling using normalized Hill functions for activation, which a few others extended to fit experimental observations~\cite{FlaJoRavAleBar:2023:acb,FlaBar:2023:asb}.
    Phenomenological approaches based on the \citet{HaiMur:1988:cpa} cross-bridge model are widely used to describe the mechanics of smooth muscle contraction due to their relative simplicity and ability to capture key features of the contraction process~\cite{StaKlaHol:2008:smc,BolSchNowSie:2012:atd,MurArnHol:2012:eam,MurHol:2014:itr,GruGruSniGalRam:2013:mma,UhlBal:2023:cmo}.
    Here, we build upon the framework established by \citet{UhlBal:2023:cmo}, where the authors modeled the stretch-dependent intracellular calcium dynamics and myosin light chain kinase~(MLCK) activation along with the calcium sensitization that occurs as a result of the inhibition of the enzyme myosin light chain phosphatase~(MLCP).

    Aside from the active response, residual stresses are inherent in arterial walls, persisting even in the unloaded state and originating from dynamic biological processes such as growth, remodeling, and adaptation.
    These processes are prominent during angiogenesis but can also be triggered by factors like hypertension, injury, or alterations to the chemo-mechanical microenvironment due to changes in blood flow, blood pressure, or extracellular matrix composition.
    Residual stresses serve to homogenize the transmural stress distribution under physiological loading, thereby optimizing load-bearing capacity and preventing localized failure~\cite{HumSch:2021:vmh}.
    Neglecting them in computational simulations leads to inaccurate predictions of arterial mechanics, such as overestimated intimal stresses or underestimated adventitial contributions.
    Mainly three different approaches to model residual stresses can be found in the literature: the opening angle method~\cite{ChuFun:1986:rsi}, the kinematic growth model~\cite{RodHogMcC:1994:sdf}, and the constrained mixture model~\cite{HumRaj:2002:acm,Hum:2021:cmm}.
    The opening angle, referring to how a radial cut in an arterial ring causes it to spring open into a sector, has historically been widely used to experimentally quantify residual stresses.
    However, a single parameter is generally insufficient to characterize the three-dimensional, layer-specific nature of residual stresses~\cite{HolSomAueRegOgd:2007:l3r}. {Even more detailed numerical approaches based on virtually closing an opened artery, cf., e.g., \citet{BalSchGro:2007:nso}, require knowledge on the geometry of the stress-free state which is usually not available.}
    Kinematic growth models overcome these limitations by employing a multiplicative decomposition of the deformation gradient to separate elastic and inelastic (growth) contributions, enabling simulation of volumetric changes over time~\cite{MenKuh:2012:fig}.
    They were further extended to include anisotropic effects, allowing for more detailed modeling~\cite{ZahBal:2018:acg}.
    These kinematic models however do not account for residual stresses that occur due to remodeling, i.e., change in the stress state that occurs without any change in volume.
    Constrained mixture models~(CMMs), first proposed by \citet{HumRaj:2002:acm}, treat the arterial wall as an evolving mixture of constituents (e.g., elastin, collagen etc.), each with their own individual stress-free reference configurations, capturing turnover, and remodeling dynamics.
    \colsr{However, CMMs require detailed knowledge of constituent-specific parameters, such as constituent specific mass fractions, turnover rates, and deposition stretches, that are difficult to measure directly and must instead be assumed or identified inversely from limited experimental data~\cite{ValHum:2009:pss,NavUtrRivLatCelGar:2023:aif}.
    Moreover, CMMs are computationally expensive as they must track the full deposition history of each constituent.}
    Recently, \citet{CyrAydHum:2016:ahc} presented a less computationally intensive, temporally homogenized version of the constrained mixture model where the individual constituents have their own elastic-growth split similar to kinematic growth models.

    To allow the analysis of the arterial wall mechanics under the influence of antihypertensive drugs, another important aspect to be modeled is the transport and distribution of these drugs in the wall.
    Transmural drug transport has been widely modeled in the context of drug-eluting stents (DES), where drugs are released from a polymeric coating into the arterial wall to inhibit smooth muscle cell proliferation and prevent restenosis.
    Early models, such as those by \citet{PonDem:2009:mom}, focused on one-dimensional advection-diffusion-reaction equations to describe drug dynamics across two porous layers (the stent coating and the arterial wall), incorporating convection due to transmural pressure gradients and metabolic degradation in the tissue.
    These approaches highlighted the role of diffusivity and residence time in optimizing release profiles.
    More advanced multiphysics frameworks have integrated FSI to couple blood flow, arterial wall deformation, and drug transport.
    For instance, \citet{CalBraBazHug:2008:mmf} proposed a model treating blood as an incompressible Navier-Stokes fluid and the arterial wall as a linear poroelastic medium, with drug advection-diffusion governed by NURBS-based isogeometric analysis for patient-specific geometries.
    To account for drug-tissue interactions, \citet{McgMckWadMcc:2013:maw} developed a coupled system with reversible binding kinetics between free and bound drug states in the extracellular matrix and target cells, solving the resulting Volterra integral equation to predict cellular and interstitial concentrations.
    Extending this, \citet{BozChoBar:2014:mtt} incorporated multilayer arterial structures (endothelium, subendothelial space, internal elastic lamina, media) with anisotropic diffusion, convection, and nonlinear binding, showing that paclitaxel is transport-dominated while sirolimus is binding-limited, underscoring the need for drug-specific release kinetics. %
    \citet{HwaWuEde:2003:iot} studied the effect of drug lipophilicity on transmural drug transport, finding that hydrophilic drugs are susceptible to advection-based transport, whereas lipophilic drugs are more diffusion-dominated.
    Although these studies focus on antiproliferative agents, their transport models are adaptable to antihypertensive drugs.
    Drug transport is modeled here using a three-dimensional reaction-diffusion equation with an isotropic diffusion tensor.
    In the present work, we aim to integrate all the aforementioned aspects into a unified computational model to simulate antihypertensive drug effects on local stress distributions in atherosclerotic arteries.

\section{Methods}
    \subsection{Mathematical model}
    We model the mechanical problem based on the balance of linear momentum neglecting inertia effects and body forces, where the passive response in the arterial wall is described by
    a hyperelastic, anisotropic material model, with two independent fiber families.
    Smooth muscle activation is considered to be mainly driven
    by stretch stimulus, and it is assumed that stimuli due to neurotransmitters and hormones remain constant.
    Active contraction
    is modeled phenomenologically based on the cross-bridge phosphorylation model by \citet{HaiMur:1988:cpa}, using the stretch-dependent model developed by \citet{UhlBal:2023:cmo}.
    Transmural drug transport is modeled by means of a reaction–diffusion equation.
    This then results in the following system of coupled nonlinear partial differential equations
    \begin{align}
    \div{\Bsigma(\bu,c)} = \bzero, \nonumber \\
    \dfrac{\partial c}{\partial t}
    \colsr{+ c \div{\bv}} 
    \colsr{-}
    \div{\left( \bD_{\Bt}\,\grad{c} \right)}
    +
    r(c) = 0,
    \label{eq:pdes}
    \end{align}
    where $\bu$ is the displacement vector, $\Bsigma$ is the Cauchy stress tensor,
    $c$ is the drug concentration \colsr{per unit current volume}, \colsr{$\frac{\partial (\cdot)}{\partial t}$ is the material time derivative, $\bv$ is the velocity,} $\bD_{\Bt}$ is the diffusion tensor, and $r(c)$ is the reaction term.
    The Cauchy stress tensor is computed by
    \eb
    \Bsigma = \bF\bS\bF^T/(\text{det}\bF)
    \quad\text{with}\quad
    \bS = 2\frac{\partial\Psi(\bu,c)}{\partial\bC},
    \ee
    where $\bF$ denotes the deformation gradient, $\bS$ the 2nd Piola-Kirchhoff stress tensor, $\Psi$ the strain energy density, and $\bC=\bF^T\bF$ the right Cauchy-Green deformation tensor.
    The divergence and gradient operators with respect to the deformed configuration are denoted by div and grad.
    The problem-specific definition of the strain energy density will be provided in the following subsections.
    The framework developed by the authors~\cite{BalHeiKlaKneRamRheSasUhl:2024:acf} is used to numerically solve the coupled problem monolithically.

    \subsubsection{Continuum-mechanical model of the arterial wall}\label{sec:cont_mech_model}
    A point in the \\ 
    three-dimensional continuum is denoted by $\bX$ in the reference configuration $\B$ and by $\bx$ in the current configuration $\Bt$.
    The time-dependent motion of the point is defined by the nonlinear deformation map $\bx = \varphi\left(\bX,t\right)$ and the displacement vector is given as $\bu = \bx-\bX$.
    The deformation gradient is then defined as $\bF\defeq \pp{\bx}{\bX}$.
    Residual stresses in the arterial wall~\cite{Ber:1961:tse} are accounted for using the kinematic growth model first proposed by \citet{RodHogMcC:1994:sdf}, and extended by \citet{ZahBal:2018:acg} to account for the anisotropic growth in arteries.
    Here, the deformation gradient is decomposed into a growth part and an elastic part
    \begin{align}
    \bF =
    \bF_{\text{e}}\bF_{\text{g}},
    \quad \text{with} \quad
    \bF_{\text{g}} =
    \bF_{\text{g}}^{(3)}\bF_{\text{g}}^{(2)}\bF_{\text{g}}^{(1)}, \nonumber \\
    \bF_{\text{g}}^{(i)} =
    \bI +
    (\vartheta^{(i)} - 1)\,\ba_{\text{g}}^{(i)} \otimes \ba_{\text{g}}^{(i)},
    \end{align}
    where $\vartheta^{(i)}$ is the growth factor in the $\ba_{\text{g}}^{(i)}$ growth direction.
    Details on the growth formulation to automatically compute reasonable residual stresses will be given in a following subsection.
    The elastic part of the deformation gradient $\bF_{\text{e}}$
    is responsible for the development of stresses.
    The elastic right Cauchy-Green tensor $\bC_{\text{e}}$ is defined as $\bC_{\text{e}} \defeq \bF_{\text{e}}{^T}\bF_{\text{e}}$.
    A coordinate-invariant formulation is used to model the effects of two non-orthogonal fiber families described by the unit vectors that lie in the circumferential-axial plane of the artery given in the reference configuration by {$\bA^{(f)}$}.
    Their directions in the intermediate configuration are given by $\ba^{(f)} = \bF_{\text{g}} \bA^{(f)}/\|\bF_{\text{g}} \bA^{(f)}\|$, where $f=1,2$ denotes the two fiber families.
    The generalized structure tensor $\bM^{(f)} \defeq \kappa_{\text{disp}} \bI + (1-3\kappa_{\text{disp}})\,\ba^{(f)} \otimes \ba^{(f)}$, where $\kappa_{\text{disp}}$ is the dispersion parameter that modulates the intensity of fiber dispersion~\cite{GasOgdHol:2006:hmo}.
    The principal invariants of $\bC_{\text{e}}$ and the mixed invariants with $\bM^{(f)}$ are defined as
    \begin{align}
    I_{1} \defeq \tr\bC_{\text{e}},
    \quad I_{3} \defeq \det\bC_{\text{e}}, \nonumber \\
    J_{4}^{(f)} \defeq \tr[\bC_{\text{e}}\bM^{(f)}],
    \quad J_{5}^{(f)} \defeq \tr[\bC_{\text{e}}^{2}\bM^{(f)}].
    \end{align}
    The strain energy density stored within the arterial wall tissue is defined to be additively decoupled into a passive and an active part, i.e.,
    \begin{align}
    \Psi \defeq \Psi_{\text{pass}} + \Psi_{\text{act}}, \nonumber \\
    \Psi_{\text{pass}} =
    \Psi_{\text{grnd}} +
    \sum_{f=1}^{2}
    \Psi_{\text{coll}}^{(f)}, \nonumber \\
    \Psi_{\text{act}} = \sum_{f=1}^{2}
    \Psi_{\text{smc}}^{(f)},
    \end{align}
    wherein the passive response, is governed by an
    isotropic elastin-based ground substance and two transversely isotropic collagen fiber families.
    Thus, it is described using a
    polyconvex strain energy density function \citep[see][]{BalNefSchHol:2006:apf} for the isotropic ground matrix and the collagen fibers
    \begin{align}
    \Psi_{\text{grnd}} &=
    \alpha_{1} \left( \dfrac{I_{1}}{I_{3}^{1/3}} - 3 \right) +
    \alpha_{2} \left( I_3^{\alpha_{3}} + \dfrac{1}{I_3^{\alpha_{3}}} - 2 \right), \nonumber \\
    \Psi_{\text{coll}}^{(f)} &=
    \alpha_{4} \langle K_{3}^{(f)} - 2 \rangle^{\alpha_{5}},
    \end{align}
    wherein $\alpha_{1}, \alpha_{2}$, and $\alpha_{4}$ are stress-like parameters and $\alpha_{3}$ and $\alpha_{5}$ are dimensionless parameters, and $K_{3}^{(f)} = I_{1}J_{4}^{(f)} - J_{5}^{(f)}$; $\langle \cdot \rangle$ denotes the Macauley bracket.

    The active contraction of the SMCs is modeled by the formulation presented in \citet{UhlBal:2023:cmo}, which includes an extension of the cross-bridge phosphorylation model by~\citet{HaiMur:1988:cpa} to describe a stretch-dependent chemical response and an extended form of \citet{MurArnHol:2012:eam}, cf. the three-element Hill muscle model~\cite{Hil:1938:tho,Fun:1970:mro}, for the chemo-mechanical interaction.
    The~\citet{HaiMur:1988:cpa} model postulates that cross-bridges, the myosin head-actin links responsible for the power-strokes driving the active contraction, exist in four distinct states based on whether they are attached and phosphorylated.
    The fraction of cross-bridges in the four states: $n_{\text{A}}$~(free dephosphorylated), $n_{\text{B}}$~(free phosphorylated), $n_{\text{C}}$~(attached phosphorylated), and $n_{\text{D}}$~(attached dephosphorylated) that sum up to unity are linked by the following system of ordinary differential equations
    \eb
    {\setlength{\arraycolsep}{2pt}
    \begin{aligned}
    &\begin{bmatrix}\dot{n}_\text{A}\\ \dot{n}_\text{B}\\ \dot{n}_\text{C}\\ \dot{n}_\text{D}\end{bmatrix}
    =\mathbf{K}\begin{bmatrix}n_\text{A}\\ n_\text{B}\\ n_\text{C}\\ n_\text{D}\end{bmatrix}, \\
    \mathbf{K}
    &=\begin{bmatrix}-k_{\text{M}}&k_{\text{P}}&0&k_{\text{d}}\\ k_{\text{M}}&-k_{\text{P}}-k_{\text{a}}&k_{\text{d}}&0\\ 0&k_{\text{a}}&-k_{\text{d}}-k_{\text{P}}&k_{\text{M}}\\ 0&0&k_{\text{P}}&-k_{\text{M}}-k_{\text{d}}\end{bmatrix}.
    \end{aligned}}
    \ee
    where $k_{\text{M}} = k_{\text{MLCK}}$ is the rate of phosphorylation, $k_{\text{P}} = k_{\text{MLCP}}$ is the rate of dephosphorylation, $k_{\text{a}}$ is the rate of attachment of phosphorylated cross-bridges to actin filaments, and $k_{\text{d}}$ is the rate of detachment of cross-bridges from actin filaments.
    The calcium- and stretch-dependent SMC activation model developed by \citet{UhlBal:2023:cmo} considers phenomenological equations to describe the calcium dynamics as well as the change in calcium sensitivity.
    The calcium-dependent rate of phosphorylation $k_{\text{MLCK}}$ is given by %
    \begin{align}
    k_{\text{MLCK}} = 
    \eta \dfrac{[\mathrm{Ca}^{2+}]_{\text{c}}^2}{[\mathrm{Ca}^{2+}]_{\text{c}}^2 +
     (\text{Ca}_{50}^{2})}, \nonumber \\
    [\mathrm{Ca}^{2+}]_{\text{c}} =
    \gamma_{1} \langle \lambda -
    \bar{\lambda}_{\text{c}} \rangle^{2},
    \label{eq:Ca}
    \end{align}
    where $[\mathrm{Ca}^{2+}]_{\text{c}}$ is the instantaneous intracellular calcium concentration, $\eta$ is a material parameter representing the maximum rate of phosphorylation, $\text{Ca}_{50}$ is the half activation calcium concentration, i.e., the concentration at which the rate of phosphorylation is half the maximum rate, $\gamma_{1}$ is a material parameter that controls the magnitude of the calcium spikes due to stretch stimulus, and $\lambda=\left({J_{4}}\right)^{1/2}$ is the \colsr{generalized dispersed-fiber stretch due to the generalized structure tensor used to compute it}.
    $\bar{\lambda}_{\text{c}}$ is the threshold stretch below which the calcium concentration reaches zero.
    To model the reduction in intracellular calcium concentration that occurs over time due to calcium sequestration into the sarcoplastic reticulum and by the action of various calcium pumps and exchangers in the cell membrane, a four parameter sigmoid function is used for the evolution of $\bar{\lambda}_{\text{c}}$
    \begin{align}
    \dot{\bar{\lambda}}_{\text{c}} =
    \dot{\bar{\lambda}}_{\text{c,min}} +
    \dfrac{\dot{\bar{\lambda}}_{\text{c,max}} - \dot{\bar{\lambda}}_{\text{c,min}}}{1 + e^{\gamma_{2} \left(\left[\mathrm{Ca}^{2+}\right]_{\text{tar}} - \left[\mathrm{Ca}^{2+}\right]_{\text{c}} - \tau_{\text{c}} \right)}}, \nonumber \\
    \left[\mathrm{Ca}^{2+}\right]_{\text{tar}} =
    \gamma_{3} \dfrac{\lambda^{2}}{\lambda^{2} +\lambda^{2}_{\text{50,c}}}.
    \label{eq:CaTar}
    \end{align}
    Herein the parameters $\dot{\bar{\lambda}}_{\text{c,min}}$ and $\dot{\bar{\lambda}}_{\text{c,max}}$ are the minimum and maximum rates of change of $\bar{\lambda}_{\text{c}}$, $\gamma_{2}$ is a material parameter, $\left[\mathrm{Ca}^{2+}\right]_{\text{tar}}$ represents the steady-state calcium concentration at a constant stretch value, $\gamma_{3}$ represents the maximum target calcium concentration achievable at any material point,
    $\lambda_{\text{50,c}}$ is the stretch at which the calcium concentration is half the target concentration and $\tau_{\text{c}}$ is chosen such that $\left[\mathrm{Ca}^{2+}\right]_{\text{c}} = \left[\mathrm{Ca}^{2+}\right]_{\text{tar}}$ at steady-state.
    For brevity, we do not present the equations for rate of dephosphorylation $k_{\text{MLCP}}$ as they will remain unchanged from the original publication~\cite{UhlBal:2023:cmo}.
    Although MLCP activity is not directly affected by calcium channel blockers, it is essential to include its stretch-dependent effects in the model to accurately capture the dynamics of cross-bridge cycling and smooth muscle contraction.
    Furthermore, the above set of equations describing the action of MLCK apply for each individual fiber direction and the superscript $(f)$ is omitted for improved legibility.

    The mechanical influence of the smooth muscle cells is modeled using the active strain energy density functions ${\Psi_{\text{smc}}}^{(f)}$ which are functions of the \colsr{generalized fiber stretch}
    $\lambda^{(f)}$ and the fraction of cross-bridges in states C~($n_{\text{C}}$) and D~($n_{\text{D}}$) as
    \begin{align}
    \Psi_{\text{smc}}^{(f)} =
    \dfrac{\mu_{\text{a}}}{2} \left(n_{\text{C}} + n_{\text{D}}\right) \left(\lambda_{\text{e}}^{(f)} - 1\right)^2, \nonumber \\
    \lambda_{\text{e}}^{(f)} = \dfrac{\lambda^{(f)}}{\lambda_{\text{a}}^{(f)}}.
    \end{align}
    Here, $\mu_{\text{a}}$ is the active stiffness parameter and $\lambda_{\text{a}}^{(f)}$ is the active part of the stretch that represents fiber shortening.
    Its evolution over time is derived using Hill's equation for muscle contraction~\cite{Kat:1939:trb} as
    \begin{align}
    \dot{\lambda}_{\text{a}}^{(f)} = 
    \beta_{1} \dfrac{P_{\text{a}}^{(f)} - P_{\text{c}}^{(f)}}{P_{\text{a}}^{(f)} - \beta_{2}}, \nonumber \\
    P_{\text{a}}^{(f)} = \dfrac{\mu_{\text{a}}}{2} \left(n_{\text{C}} + n_{\text{D}}\right) \left(\lambda_{\text{e}}^{(f)} - 1\right), \nonumber \\
    P_{\text{c}}^{(f)} = \kappa n_{\text{C}},
    \end{align}
    where $\beta_{1}$ and $\beta_{2}$ are material parameters and $\kappa$ denotes the maximum achievable active force.

    \subsubsection{Transmural drug transport and drug-tissue interaction}
\begin{figure}[t]
        \centering
        \includegraphics[width=0.475\textwidth]{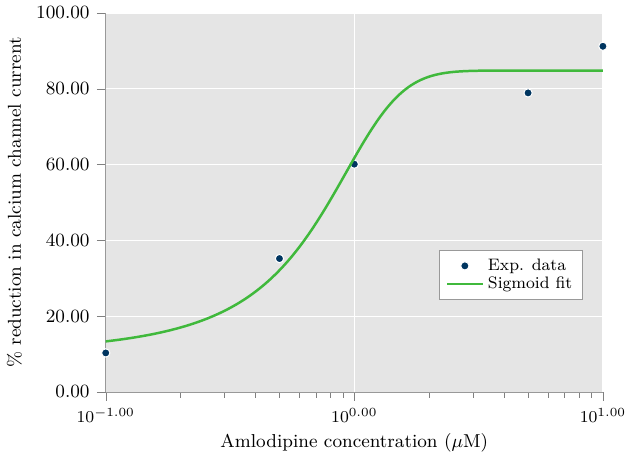}
        \caption{Amlodipine-induced reduction in intracellular calcium concentration depicted by drug concentration vs. reduction in calcium channel current. Blue dots represent experimental data from \citet{WanJiaLiLai:2008:eos}, whereas the green curve represents the sigmoid function fit to the data
        \label{fig:drug_interaction}}
    \end{figure}

    Transport of drugs through the arterial wall is a combination of solvent-driven flow~(advection) and random molecular agitation~(diffusion). 
    The susceptibility of each drug to advection- and diffusion-mediated transport varies based on its hydrophilicity~\cite{HwaWuEde:2003:iot}.
    Hydrophilic drugs are prone to advective transport, whereas lipophilic drugs are transported by diffusion.
    Here, we consider specifically amlodipine, which is a lipophilic drug and its transmural transport is hence modeled by means of a reaction-diffusion equation with a linear reaction term $r(c) = Kc$, where $K$ is the drug absorption rate.
    In general, the reaction term may be nonlinear depending on specific drug type and tissue interactions.
    Although the diffusivity is in general anisotropic with diffusion in the radial direction being an order slower~\cite{LevVukHwaEde:2004:sbt,BozChoBar:2014:mtt},
    due to a lack of data regarding detailed diffusion properties and for computational simplicity,
    the diffusion tensor was assumed to be isotropic here and is given by $\bD_{\Bt} = D_{\text{aml}} \bI$, where $D_{\text{aml}}$ is the drug diffusivity.
    This is also motivated by the fact that in our study presented later, the spatial distribution of the drug in the artery resulting from the diffusion is of much higher relevance than a high accuracy in the time-dependent diffusion processes.
    Dihydropyridine Calcium Channel Blockers~(CCBs) 
    affect smooth muscle contraction by reducing the flow of calcium ions into the cytoplasm through L-type calcium channels.
    This effect is modeled here by making the proportionality constants for the instantaneous intracellular calcium concentration $\left[\mathrm{Ca}^{2+}\right]_{\text{c}}$ and the target intracellular calcium concentration $\left[\mathrm{Ca}^{2+}\right]_{\text{tar}}$ dependent on the drug concentration $c$ as follows
    \begin{align}
    \gamma_{1}(c) = \gamma_{1,\text{min}} + \dfrac{\gamma_{1,\text{max}} - \gamma_{1,\text{min}}}{1 + \exp\left(p_{1} \left(c - c_{\text{50}}\right)\right)}, \nonumber \\
    \gamma_{3}(c) = \gamma_{3,\text{min}} + \dfrac{\gamma_{3,\text{max}} - \gamma_{3,\text{min}}}{1 + \exp\left(p_{2} \left(c - c_{\text{50}}\right)\right)},
    \label{eq:gamma}
    \end{align}
    where $\gamma_{1,\text{min}}, \gamma_{1,\text{max}}$ and $\gamma_{3,\text{min}}, \gamma_{3,\text{max}}$ are the minimum and maximum values of $\gamma_{1}$ and $\gamma_{3}$.
    $c_{\text{50}}$ is the concentration at which the parameters reach half the maximum value, and $p$ is the curvature parameter that determine the steepness of the sigmoid function.
    Both 
    $\gamma_{1}$ and $\gamma_{3}$ are made dependent on the drug concentration, since a reduction in calcium channel current would not only reduce the magnitude of transient peaks which are controlled by $\gamma_{1}$, but also result in a reduction in steady-state calcium concentration controlled by $\gamma_{3}$.
    The choice of the above sigmoid functions is motivated by experimental observations~\cite{HugWij:1993:tao,WanJiaLiLai:2008:eos}.
    To this end, Fig.~\ref{fig:drug_interaction} shows the amlodipine-induced reduction in calcium channel current with increasing drug concentration, where data obtained from \citet{WanJiaLiLai:2008:eos} shows a sigmoid trend.
    \subsubsection{Residual stresses and fiber reorientation}
    \label{sec:residual_stresses}
    Residual stresses are incorporated using the anisotropic kinematic growth model proposed in \citet{ZahBal:2018:acg}, where the evolution of the growth factors is given according to \cite{LubHog:2002:otm} by
    \begin{align}
    \dot{\vartheta}^{(i)} &=
    k_{\vartheta}^{(i)}(\vartheta^{(i)}) \, \phi^{(i)}(\BSigma_{\text{e}}), \nonumber \\
    k_\vartheta^{(i)} &= 
    \begin{cases}
        k_{\vartheta,(i)}^{+} \left[\dfrac{\vartheta^{+}_{(i)}-\vartheta^{(i)}}{\vartheta^{+}_{(i)}-1}\right]^{m^{+}_{\vartheta,(i)}} & \text{for} \,\, \phi^{(i)}>0 \\
        k_{\vartheta,(i)}^{-} \left[\dfrac{\vartheta^{(i)}-\vartheta^{-}_{(i)}}{1-\vartheta^{-}_{(i)}}\right]^{m^{-}_{\vartheta,(i)}} & \text{for} \,\, \phi^{(i)}<0 \\
        0 &\mathrm{for} \, \phi^{(i)}=0
    \end{cases} 
    \end{align}
    where $k_\vartheta^{(i)}$ is the growth velocity function, which scales the growth rate by a factor of $k_{\vartheta,(i)}^{+/-}$ and limits the growth factors to the range $[\vartheta^{-}_{(i)}, \vartheta^{+}_{(i)}]$, while the parameters $m^{+/-}_{\vartheta,(i)}$ determine the shape of function in its non-zero range. The growth driving forces $\phi^{(i)}$ are defined as
    \begin{align}
    \phi^{(i)} (\BSigma_\text{e}) &= 
    \begin{cases}
        \Sigma_{\text{e}}^{(i)}, & \text{for} \,\, i=1,2 \\
        \trb{\BSigma_{\text{e}}}, & \text{for} \,\, i=3
    \end{cases}, \nonumber \\
    \Sigma_{\text{e}}^{(i)} &= \BSigma_{\text{e}} : \left( \be^{(i)} \otimes \be^{(i)} \right),
    \end{align}
    where $\be^{(i)}$ denote the eigenvectors of the symmetric part of the Mandel stress tensor $\BSigma_{\text{e}} = \symb{\bC_{\text{e}}\bS_{\text{e}}}$.
    Here, the driving force for growth in the direction of the first two principal directions is assumed to be the Mandel stress in the respective directions, while for the third growth direction, the driving force is assumed to be the trace of the Mandel stress tensor.
    The choice of the driving forces is motivated by \citet{Zah:2021:mog}, where various growth driving forces were compared based on an optimization principle, and the above choice was found to be most physiological.

    Remodeling is accounted for by the evolution of the fiber directions $\ba^{(f)}$. They are assumed to evolve in such a way that they align with the sum and difference of first two principal values of the Mandel stress tensor defined in the intermediate configuration as
    \begin{align}
        \bar{\ba}_\mathrm{targ,int}^{(1)} &= \big\langle\Sigma_\mathrm{e}^{(1)}\big\rangle\,\be^{(1)} + \big\langle\Sigma_\mathrm{e}^{(2)}\big\rangle\,\be^{(2)}, \nonumber \\
        \bar{\ba}_\mathrm{targ,int}^{(2)} &= \big\langle\Sigma_\mathrm{e}^{(1)}\big\rangle\,\be^{(1)} - \big\langle\Sigma_\mathrm{e}^{(2)}\big\rangle\,\be^{(2)},
    \end{align}
    where $\be^{(1)}$ and $\be^{(2)}$ are the first two eigenvectors of the symmetric part of the Mandel stress and $\Sigma_\mathrm{e}^{(1)}$ and $\Sigma_\mathrm{e}^{(2)}$ are the corresponding eigenvalues.
    Then an evolution equation for the angle between the fiber direction $\ba^{(f)}$ and the target fiber direction $\ba_{\text{targ}}^{(f)}$ is defined as
    \begin{align}
        \dot{\eta}^{(f)}=
        - \dfrac{k_{\eta}^{+}\eta^{2}}{\eta^{2} + m_{\eta}^{+}}, \nonumber \\
        \eta^{(f)} = \arccos\left(\ba^{(f)} \cdot \ba_{\text{targ}}^{(f)}\right),
    \end{align}
    where $\ba_{\text{targ}}^{(f)}$ is the target fiber direction in the intermediate configuration, $k_{\eta}^{+}$ is a parameter that determines the rate of angle evolution, and $m_{\eta}^{+}$ controls the shape of the evolution function.
    \subsection{Choice of model and material parameters}
    Usually, parameters are obtained by fitting the model behavior to experimental data.
    However, in case of soft biological tissue, there are limitations to the types of experiments that can be conducted, in particular on human tissue, which consequently results in a lack of experimental data in general.
    For instance, material testing of components of arteries excised from deceased patients can provide sufficient data to determine the parameters for the passive material response due to elastin, collagen, and connective tissue~\cite{HolSomReg:2004:amp}.
    Whereas in the context of active arterial response due to SMC activation, this type of testing is difficult as it requires intact and viable SMCs, meaning the arteries would need to be excised within a few hours post mortem.
    This is not only operationally challenging but also ethically complicated in human subjects.
    In these cases, experiments are usually conducted on animal tissue~\cite{JohElyTakWalWalCol:2009:css} to obtain at least characteristic physiologically relevant data which may serve as basis for the parameter adjustment.
    Ideally all the necessary experiments are performed on the same specimen or, since this is in most cases virtually impossible, at least on a similar arterial tissue.
    However, as is commonly the case in biological modeling, experimental data in the 
    existing literature is often fragmented, with several studies using not only different arteries but also different species in investigating various aspects of arterial behavior.
    Here, we have utilized medial SMC response parameters obtained from experiments on rat middle cerebral arteries.
    Whereas the passive material parameters for the adventitia and fibrous cap were obtained from experiments on human iliac arteries.
    Rat SMCs and human SMCs have similar gross morphology and contractile mechanisms and exhibit similar Ca$^{2+}$ sensitivities~\cite{AalMul:1983:har}.
    In collagen contraction assays conducted on rat SMCs, nine out of ten vasoactive compounds tested elicited a response consistent with known in-vivo human SMC behavior, demonstrating a high predictive value~\cite{TreBieZafPol:2023:sco}.
    Although fundamental differences like cell size and {receptor-specific} responses exist between species, the basic contractile mechanisms are comparable.
    The maximal active stress measured in rats~\citep[ca. 74~kPa][]{ArnMalUve:1984:sam} corresponds well to human measurements~\citep[ca. 50--60~kPa][]{FraBreGioKasHolAma:2022:ros}.
    Despite these similarities, it is important to acknowledge that using parameters from different species and arterial types qualifies the model predictions.
    We attempt to mitigate this limitation by conducting parameter variation studies to investigate the sensitivity of the results to changes in key parameters within physiologically reasonable ranges.
    As will be demonstrated subsequently, the simulations, despite their limitations, provide characteristic insights into the drug-artery interactions under investigation.

    \subsubsection{Structural parameters}
\begin{figure}[htb]
        \centering
        \includegraphics[width=0.475\textwidth]{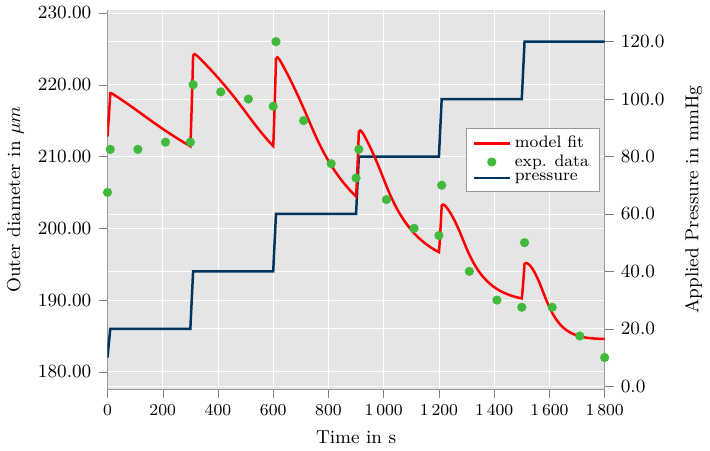}
        \caption{SMC model parameter fit to experimental data \citep[reproduced from][]{UhlBal:2023:cmo}. Green dots represent experimental data from \citet{JohElyTakWalWalCol:2009:css}, whereas the red curve represents the model fit. The blue curve shows the applied pressure step function as part of the experimental protocol. The model is able to reproduce the transient and steady-state response of the artery to changes in intraluminal pressure.}
        \label{fig:smc_parameter_fit}
    \end{figure}

    Material parameters for the active material response that were identified in \citet{UhlBal:2023:cmo} are used for the healthy media.
    Therein, the authors fitted model and material parameters to experiments conducted on the middle cerebral arteries of male Sprague–Dawley rats~\cite{JohElyTakWalWalCol:2009:css}.
    Fig.~\ref{fig:smc_parameter_fit} shows the result of the model parameter fit \cite{UhlBal:2023:cmo} to experimental data obtained from~\citet{JohElyTakWalWalCol:2009:css}.
    To account for the vasodilatory effects of endothelium-produced gasotransmitters~(most importantly nitric oxide) which were not considered in the experimental study, we reduced the value of the model parameter $\kappa$~(corresponding to the maximum achievable active force) by 30\%.
    All other parameters are kept the same as in the original study and can be obtained from Table~2 and 3 in \citet{UhlBal:2023:cmo}.
    The degenerated media in atherosclerotic plaques is in general stiffer than the healthy media and exhibits a pronounced nonlinear response even at low strain levels~\cite{HolSomReg:2004:amp}.
    In addition, the SMCs in the degenerated media are known to switch to synthetic/proliferative phenotypes and produce less active force~\citep[see Section 3.3 in][]{DelTou:2019:tov}.
    To account for these observations, the material parameters of the degenerated media are adjusted accordingly and the model parameter $\kappa$ is reduced by 33\% compared to the healthy media to account for the lesser number of contractile SMCs.
    The lipid pool is modeled as a low stiffness, isotropic rubber-like material with an effective Young's modulus of the order of $1 \, \text{kPa}$, while the calcifications are stiff and isotropic~\cite{YanBacZheNaqWooTenBilTan:2009:ivi,OhaFinGhaHerTraHerRioKotElaPet:2008:nct}.
    
\begin{figure*}[t]
    \centering
    \begin{minipage}[b]{0.47\linewidth}
        \vspace{0pt}
        \includegraphics[width=\linewidth,page=1]{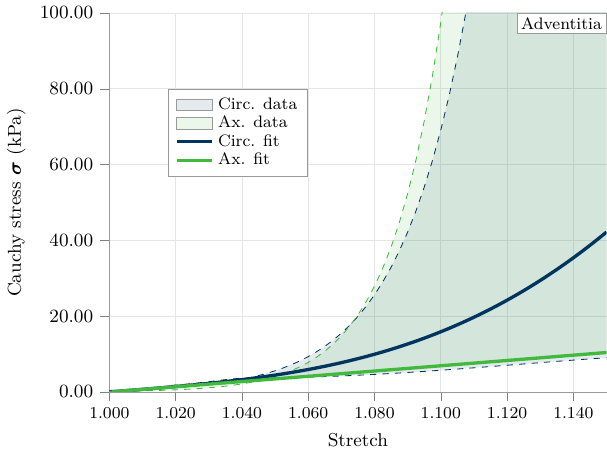}
        \par\vspace{-1.6\baselineskip}\noindent\llap{\footnotesize(a)\hspace{-0.6em}}
        \vspace{0.8\baselineskip}
    \end{minipage}
    \hfill
    \begin{minipage}[b]{0.47\linewidth}
        \vspace{0pt}
        \includegraphics[width=\linewidth,page=1]{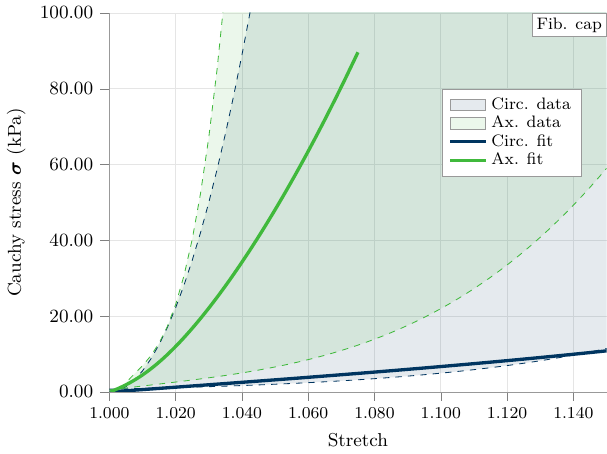}
        \par\vspace{-1.6\baselineskip}\noindent\llap{\footnotesize(b)\hspace{-0.6em}}
        \vspace{0.8\baselineskip}
    \end{minipage}
    \caption{Stress stretch response of model using the material parameters fitted to uniaxial test data from \citet{HolSomReg:2004:amp} for (a) adventitia and (b) fibrous cap. The area between the dashed lines represents the experimental data region covering the range of external iliac arteries from data, whereas the solid lines represent the model fit to artery I from \citet{HolSomReg:2004:amp}. Blue curves correspond to circumferential loading and green curves to axial loading}
    \label{fig:parameter_fit}
    \end{figure*}

    The material parameters $\alpha_{1}$ and $\alpha_{4}$ for the adventitia and fibrous cap were obtained by fitting them to uniaxial tensile test data from \citet{HolSomReg:2004:amp}, whereas $\alpha_{5}$ was assumed to be 4.0 for the adventitia and 2.5 for the fibrous cap.
    The parameter $\alpha_{5}$ controls the nonlinearity of the stress-stretch response, with higher values leading to a stiffer response at higher stretch values.
    The values chosen here offer a reasonable compromise between achieving a good fit to the experimental data and ensuring numerical stability during simulations, cf. the discussions in~\citet{BalHeiKlaRheSch:2023:coa} and the comparison of the model response with experiments in Fig.~\ref{fig:parameter_fit}.
    Since the corresponding fiber angles were not reported in the experimental study, they were also optimized.
    Fitting was done by minimizing the cost function given as
    \begin{align}
    z ={}& \dfrac{1}{N_{\text{circ.}}} \sum_{i=1}^{N_{\text{circ.}}} \left(\dfrac{\sigma_{\text{circ.}}^{\text{exp.}}(\lambda_{i}) - \sigma_{\text{circ.}}^{\text{sim.}}(\lambda_{i})}{\operatorname{Max}(\sigma_{\text{circ.}}^{\text{exp.}})}\right)^{2} \nonumber \\
    &+ \dfrac{1}{N_{\text{ax.}}} \sum_{j=1}^{N_{\text{ax.}}} \left(\dfrac{\sigma_{\text{ax.}}^{\text{exp.}}(\lambda_{j}) - \sigma_{\text{ax.}}^{\text{sim.}}(\lambda_{j})}{\operatorname{Max}(\sigma_{\text{ax.}}^{\text{exp.}})}\right)^{2}.
    \end{align}
    Parameter optimization was carried out using the differential evolution optimization algorithm implemented in Mathematica.
    Differential evolution is a population-based, evolutionary algorithm suitable for nonlinear, multi-dimensional optimization problems where the cost function is not easily analytically differentiable~\cite{StoPri:1997:dea}.
    The solid lines in Fig.~\ref{fig:parameter_fit} show the results of the material parameter fit for the adventitia and fibrous cap regions.
    Material parameters were optimized to fit the stress-stretch response of Artery I from \citet{HolSomReg:2004:amp}, where the authors conducted uniaxial tensile tests on 9 different human iliac artery specimens, and reported a significant variation in their mechanical response, which is depicted by the shaded area between the dashed lines in Fig.~\ref{fig:parameter_fit}.
    Table~\ref{tab:material_parameters} details the chosen values of material parameters
    for the different plaque components.
    \colsr{Additionally, the parameter $\alpha_{3}$ was set to $2.757$ for all components~\cite{UhlBal:2023:cmo} whereas $\alpha_{2}$, a penalty-type parameter was chosen component-wise as shown in Table~\ref{tab:material_parameters} to ensure incompressibility while maintaining numerical stability.
    The fiber dispersion parameter $\kappa_{\text{disp}}$ was chosen to be $0.133$ for all components.
    }
    \begin{table*}[t]
        \caption{Component-wise material parameters adopted \label{tab:material_parameters}}
        \begin{center}
            \begin{tabular}{rcccccc} \toprule
                Parameter & $\alpha_1$(kPa) & \colsr{$\alpha_2$(kPa)} & $\alpha_4$(kPa) & $\alpha_5$ & $\mu_{\text{a}}$(kPa) & $\kappa$(kPa) \\ \midrule
                Media & $11.52507$ & \colsr{$151.738$} & $1.27631$ & $3.088$ & $11.857$ & $103.78$ \\
                Adventitia & $11.1251$ & \colsr{$198.654$} & $370.4226$ & $4.0$ & - & - \\
                Degen. media & $13.8301$ & \colsr{$151.738$} & $3.19078$ & $3.088$ & $7.251$ & $69.53$ \\
                Lipid pool & $1.3$ & \colsr{$500.0$} & - & - & - & - \\
                Calcifications & $11.52 \cdot 10^{3}$ & \colsr{$151.738$} & - & - & - & - \\
                Fibrous cap & $9.2896$ & \colsr{$165.528$} & $73375.455$ & $2.5$ & - & - \\ \bottomrule
            \end{tabular}
        \end{center}
    \end{table*}
    \begin{table}[t]
        \centering
        \caption{\colsr{Component-wise transport parameters} \label{tab:transport_parameters}}
            \begin{tabular}{rcc} \toprule
                Parameter & $D_{\text{aml.}}(\text{cm}^{2}\text{s}^{-1})$ & $\colsr{K}(\text{s}^{-1})$ \\ \midrule
                Media & $7 \cdot 10^{-5}$ & 0 \\
                Adventitia & $7 \cdot 10^{-5}$ & - \\
                Degen. media & $3.5 \cdot 10^{-5}$ & 0 \\
                Lipid pool & $1.45 \cdot 10^{-8}$ & - \\
                Calcifications & $1.45 \cdot 10^{-8}$ & - \\
                Fibrous cap & $6.23 \cdot 10^{-6}$ & - \\ \bottomrule
            \end{tabular}
    \end{table}
    \subsubsection{Diffusivity and drug-tissue interaction parameters}
    To model transmural amlodipine diffusion, we require drug diffusivities in the different plaque components.
    In our literature survey, we unfortunately did not come across studies that have investigated the diffusivity of amlodipine through the arterial wall or its components.
    However, there have been several studies on the diffusion/transport process of antiproliferative drugs like paclitaxel through the healthy wall, in the context of drug-eluting stents. %
    Paclitaxel, a lipophilic drug, exhibits a molecular weight and logP value that are comparable to those of amlodipine, suggesting similar diffusivity.
    \citet{BozChoBar:2014:mtt} and \citet{HosHosBazCalHug:2012:mmo} report that the diffusivity of paclitaxel through the media, fibrous cap, and lipid pool is on the order of $10^{-7}\, \text{cm}^{2}\text{s}^{-1}$, $10^{-8}\, \text{cm}^{2}\text{s}^{-1}$, and $10^{-13}\, \text{cm}^{2}\text{s}^{-1}$, respectively.
    The primary interest of this study lies in the steady-state distribution of amlodipine and its effect on the active tissue response, which depends on the relative values of diffusivities among the various components.
    In other words, we are interested in how much faster or slower amlodipine diffuses through one component compared to another.
    In the interest of minimizing the number of time steps required, we virtually accelerate the diffusion process by uniformly multiplying the diffusivities by a factor of 100 for all components.
    In the absence of reaction terms, this scaling is equivalent to scaling the time step sizes by a factor of 100.
    This accelerated diffusion is necessary not only due to the slow time scale of diffusion compared to SMC activity but also because increasing time step sizes beyond 0.3s leads to convergence issues in the solution of the time-dependent, nonlinear equations of the active SMC model at some Gauss integration points.
    This approach preserves the relative diffusivity differences between components while allowing us to reach steady-state at an earlier time.
    Preliminary simulations confirmed that this scaling does not affect the steady-state distribution of the drug.
    
    Since component specific diffusivity data is unavailable for the adventitia and degenerated media, reasonable assumptions are made based on their composition and density.
    Initially, the reaction term coefficient $K$ is assumed to be zero, indicating no drug absorption by the tissue. %
    Table~\ref{tab:drug_interaction_parameters} lists the values of the drug-tissue interaction parameters $\gamma_{1,\text{min}}$, $\gamma_{3,\text{min}}$, $c_{\text{50}}$, and $p$ that were obtained by fitting them to data from \citet{WanJiaLiLai:2008:eos}, who explored the relationship between amlodipine concentration and L-type calcium channel current inhibition in aortic smooth muscle cells of rats.
    \begin{table*}[htb]
        \caption{Drug-tissue interaction parameters.
        $\gamma_{1,\text{max}}$ and $\gamma_{3,\text{max}}$ were obtained from \citet{UhlBal:2023:cmo}, whereas the other parameters were fitted to experimental data from \citet{WanJiaLiLai:2008:eos}, cf.~Fig.~\ref{fig:drug_interaction}
        \label{tab:drug_interaction_parameters}}
        \begin{center}
            \begin{tabular}{cccccc} \toprule
                $\gamma_{1,\text{min}}(\mu \text{M})$ & $\gamma_{1,\text{max}}(\mu \text{M})$ & $\gamma_{3,\text{min}}(\mu \text{M})$ & $\gamma_{3,\text{max}}(\mu \text{M})$ & $c_{\text{50}}(\mu \text{M})$ & $p(\mu \text{M}^{-1})$\\ \midrule
                $0.07798$ & $0.5131$ & 0.13678 & 0.9 & $0.667345$ & $2.9581$ \\ \bottomrule
            \end{tabular}
        \end{center}
    \end{table*}
    \subsection{Atherosclerotic artery simulations}
    The weak forms of the coupled deformation-diffusion partial differential equations~\eqref{eq:pdes} in the reference configuration are given as
    \eb
    \displaystyle  \int_{\cal B}\half\delta\bC : \bS
    \,\text{d}V
    - \displaystyle \int_{\partial B}\delta\bu \cdot \bT_0\, \text{d}A
    = 0 \, ,
    \ee
    \begin{align}
    \colsr{
    \displaystyle \int_{\cal B}\delta c \left(
    \dfrac{\partial (J\,c)}{\partial t}
    +
    J\,r(c)\right)
    \text{d}V} \nonumber \\
    +  
    \displaystyle
    \int_{\cal B}
    \Grad{\delta c}\cdot\bD\cdot\operatorname{Grad}c\,J\,\text{d}V \nonumber \\
    - \int_{\partial B} \delta c\,\, q^{\cal B}\,J\,\text{d}A = 0\, .
    \end{align}
    Herein, $\Psi$ is the total strain energy density function for each component, $\bu$ is the displacement vector, $\delta\bC$ is the Green-Lagrange deformation tensor evaluated for the virtual displacement vector~$\delta \bu$,
    $\bT_0$ is the traction vector on the boundary $\partial \cal B$,
    $c$ is the drug concentration, $\delta c$ is the variation of drug concentration, $J$ is the determinant of the deformation gradient, $r(c)$ denotes the reaction term, \colsr{$\bD = \bF^{-1} \bD_{\Bt} \bF^{-T}$} is the diffusion tensor, $\operatorname{Grad}$ denotes the gradient operator in the reference configuration, \colsr{the material time derivative is denoted by $\frac{\partial(\cdot)}{\partial t}$}, and $q^{\cal B}=\bD\cdot\Grad{c}\cdot\bn_0$ is the flux across the boundary with $\bn_0$ being the outward normal vector in the reference configuration.
    The weak forms are discretized using standard displacement tetrahedral finite elements with quadratic shape functions for~$\bu$, $c$, $\delta\bu$, and~$\delta c$.
    The resulting system of nonlinear equations is solved with the Newton-Raphson scheme and load stepping with non-equidistant time steps, which are chosen depending on the individual part of the loading protocol.
    Time integration of the diffusion equation along with all evolution equations given in Section~\ref{sec:cont_mech_model} is performed using the implicit Euler method.
    An exception is the evolution of fiber directions, which is integrated explicitly and offline in the sense that the updated fiber directions are used only in the time step after they have been computed.
    Local nonlinear equations at the Gauss points resulting from implicit time integration are solved using a separate Newton algorithm rendering the overall solution scheme a multilevel Newton-Raphson finite element scheme.
    The element level code, including the consistent tangents required in the linearization of the nonlinear system of equations, was generated using the Wolfram Mathematica based package AceGen~\cite{KorWri:2016:aof}.
    This was then integrated into the finite element domain decomposition framework FEDDLib \cite{feddlib}, where nonlinear solvers from Trilinos \cite{trilinos, trilinos:2025:arxiv} along with efficient preconditioners from FROSch \cite{heinlein_parallel_2016a,Heinlein:2020:FRO} were used to monolithically solve the fully coupled system of equations; see also \cite[Sec. 4]{BalHeiKlaKneRamRheSasUhl:2024:acf} for more details. %

    \subsubsection{Atherosclerotic geometries}
\begin{figure}[t]
        \centering
        \begin{minipage}[b]{\linewidth}
            \vspace{0pt}
            \includegraphics[width=\linewidth,page=1]{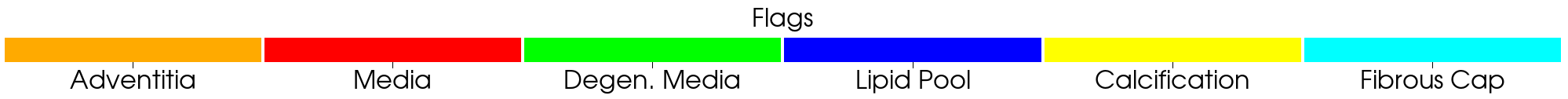}
        \end{minipage}
        \hfill
        \begin{minipage}[b]{0.3\linewidth}
            \vspace{0pt}
            \includegraphics[width=\linewidth,page=1]{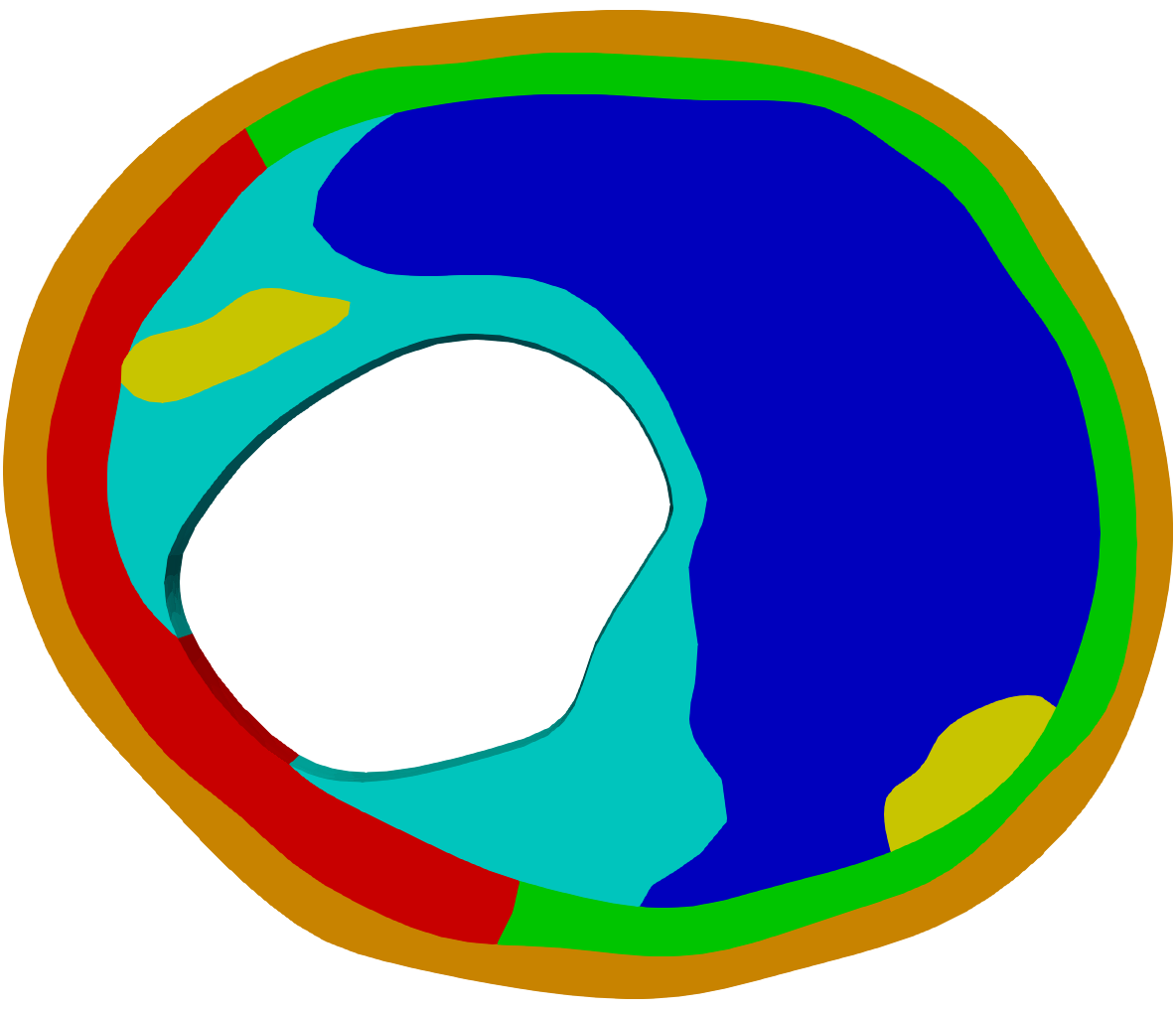}
            \par\vspace{-1.1\baselineskip}\noindent\llap{\footnotesize(a)\hspace{-0.6em}}
            \vspace{0.8\baselineskip} %
        \end{minipage}
        \hfill
        \begin{minipage}[b]{0.3\linewidth}
            \vspace{0pt}
            \includegraphics[width=\linewidth,page=1]{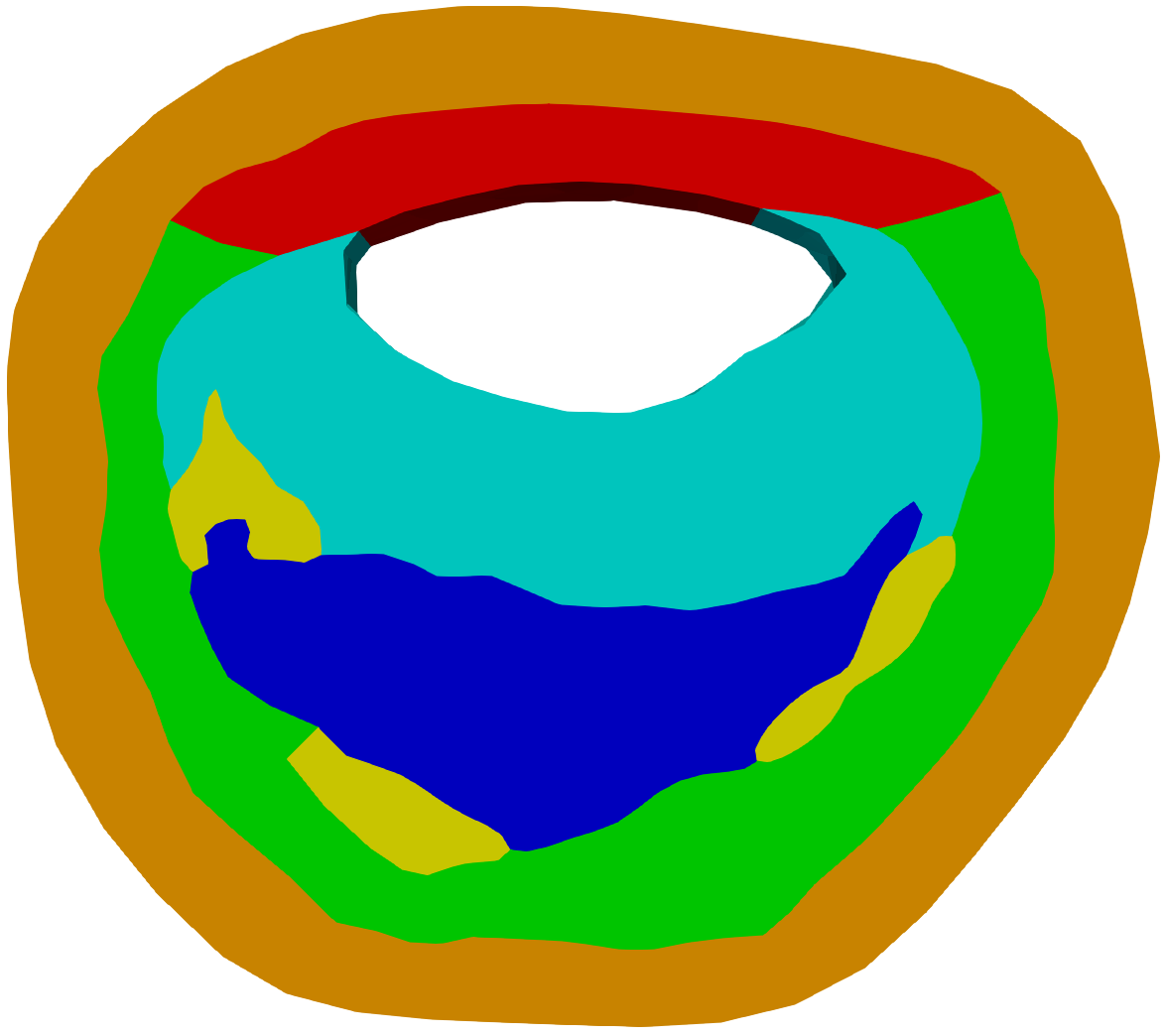}
            \par\vspace{-1.1\baselineskip}\noindent\llap{\footnotesize(b)\hspace{-0.6em}}
            \vspace{0.8\baselineskip}
        \end{minipage}
        \hfill
        \begin{minipage}[b]{0.3\linewidth}
            \vspace{0pt}
            \hspace{0.05\linewidth}
            \includegraphics[width=0.9\linewidth,page=1]{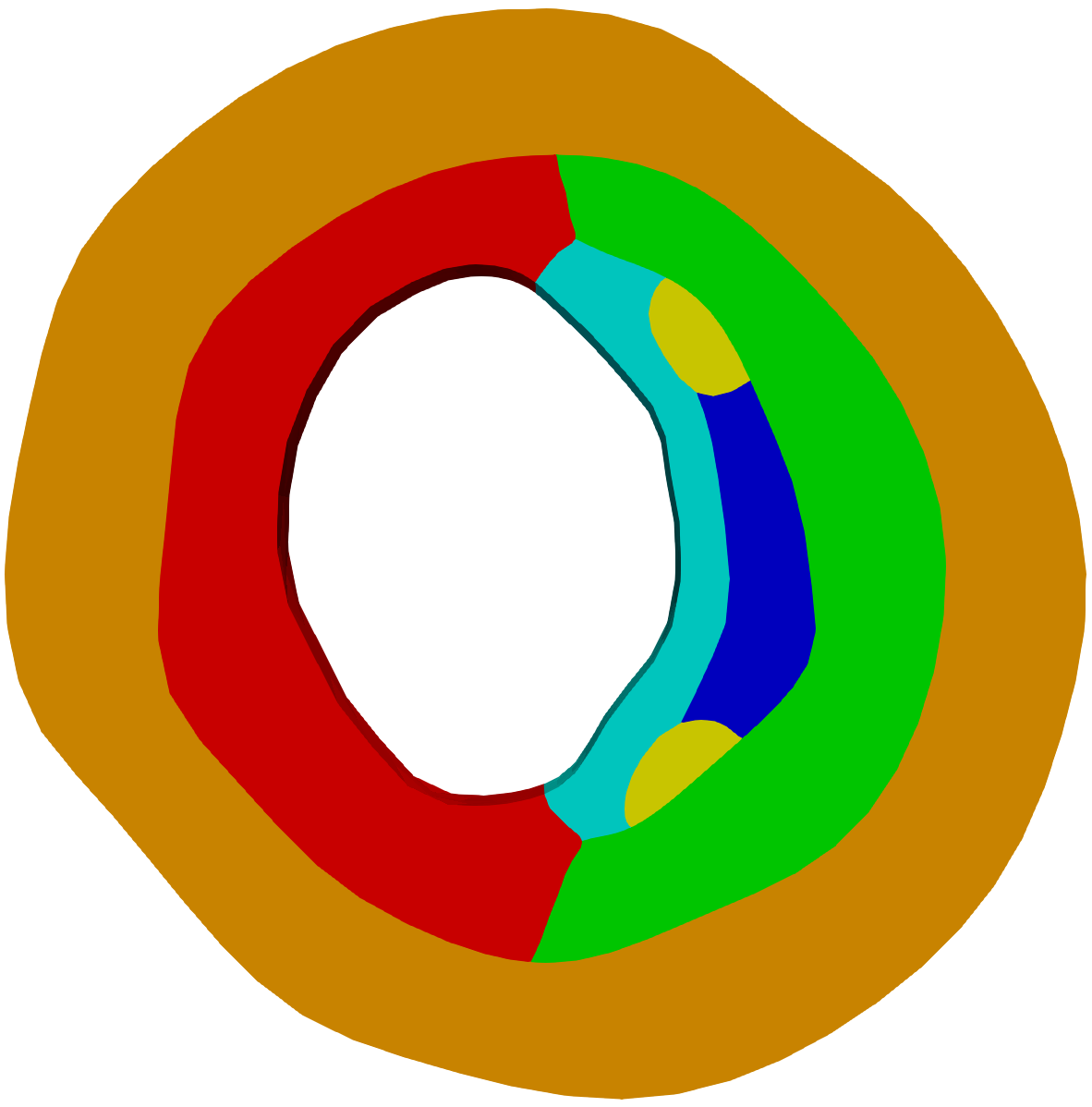}
            \par\vspace{-1.1\baselineskip}\noindent\llap{\footnotesize(c)\hspace{-0.6em}}
            \vspace{0.8\baselineskip}
        \end{minipage}
        \hfill
        \begin{minipage}[b]{0.3\linewidth}
            \vspace{0pt}
            \includegraphics[width=\linewidth,page=1]{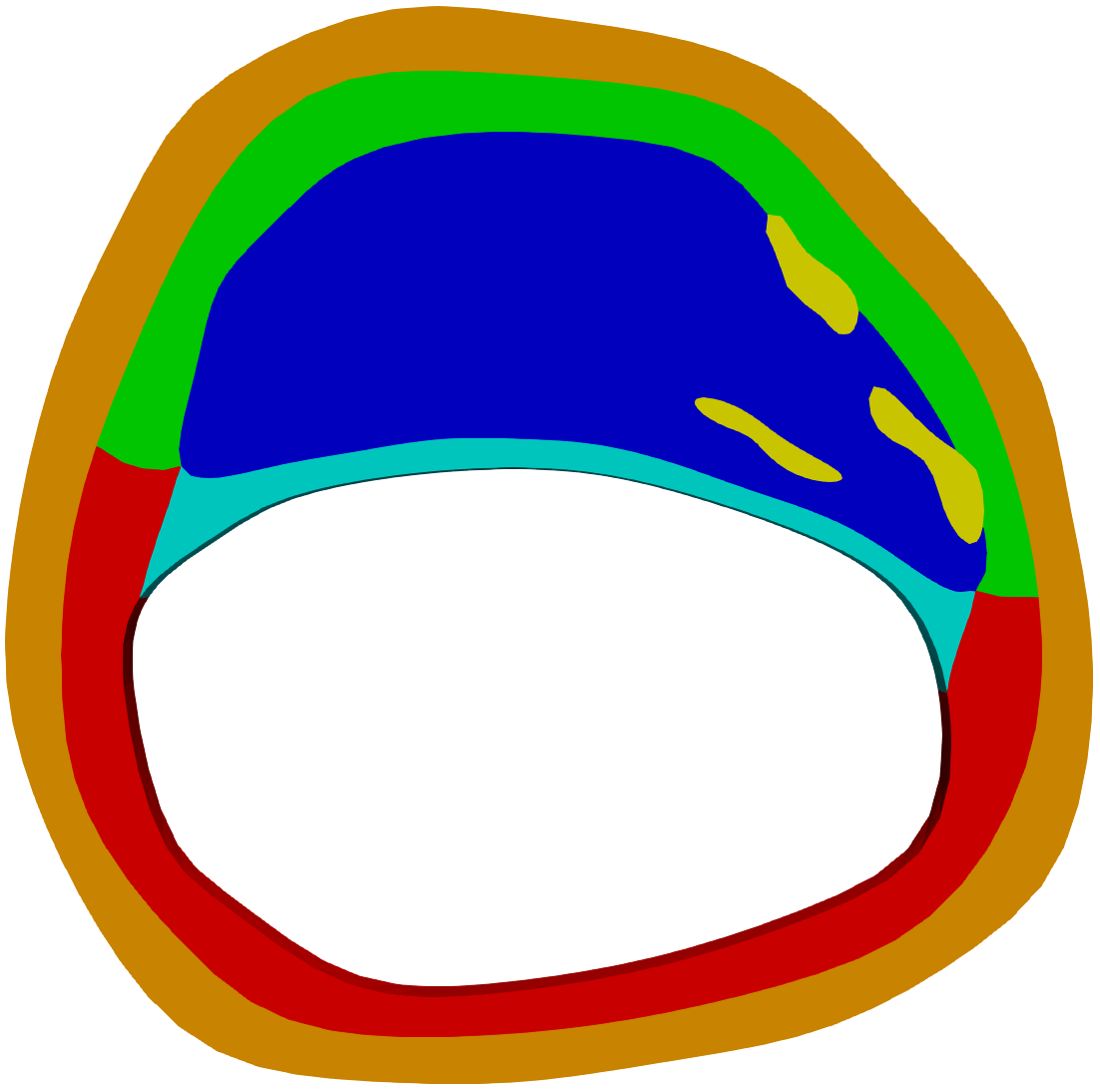}
            \par\vspace{-1.1\baselineskip}\noindent\llap{\footnotesize(d)\hspace{-0.6em}}
            \vspace{0.8\baselineskip} %
        \end{minipage}
        \hfill
        \begin{minipage}[b]{0.3\linewidth}
            \vspace{0pt}
            \includegraphics[width=0.95\linewidth,page=1]{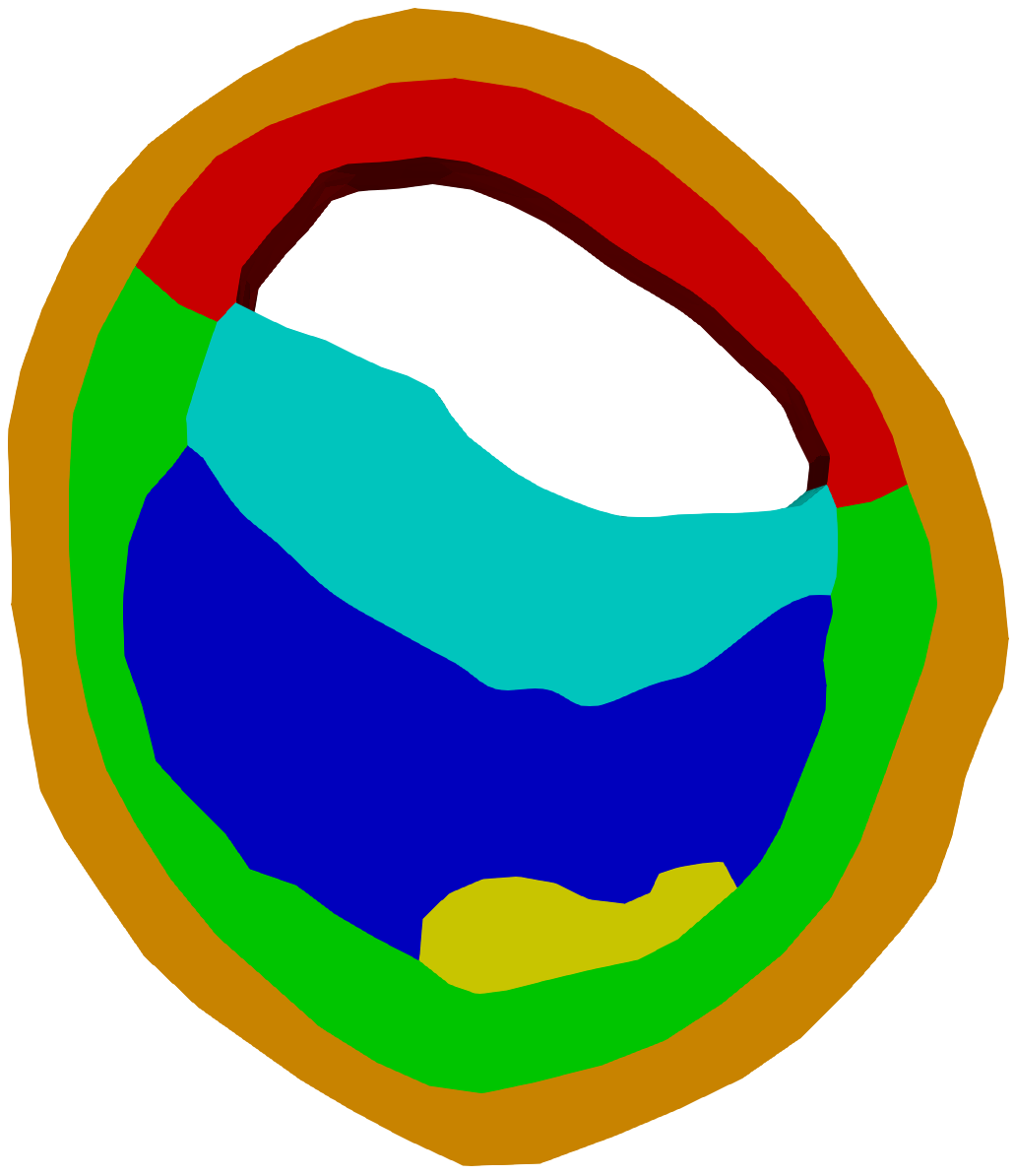}
            \par\vspace{-1.1\baselineskip}\noindent\llap{\footnotesize(e)\hspace{-0.6em}}
            \vspace{0.8\baselineskip}
        \end{minipage}
        \hfill
        \begin{minipage}[b]{0.3\linewidth}
            \vspace{0pt}
            \includegraphics[width=\linewidth,page=1,clip,trim=15cm 1cm 28cm 1cm]{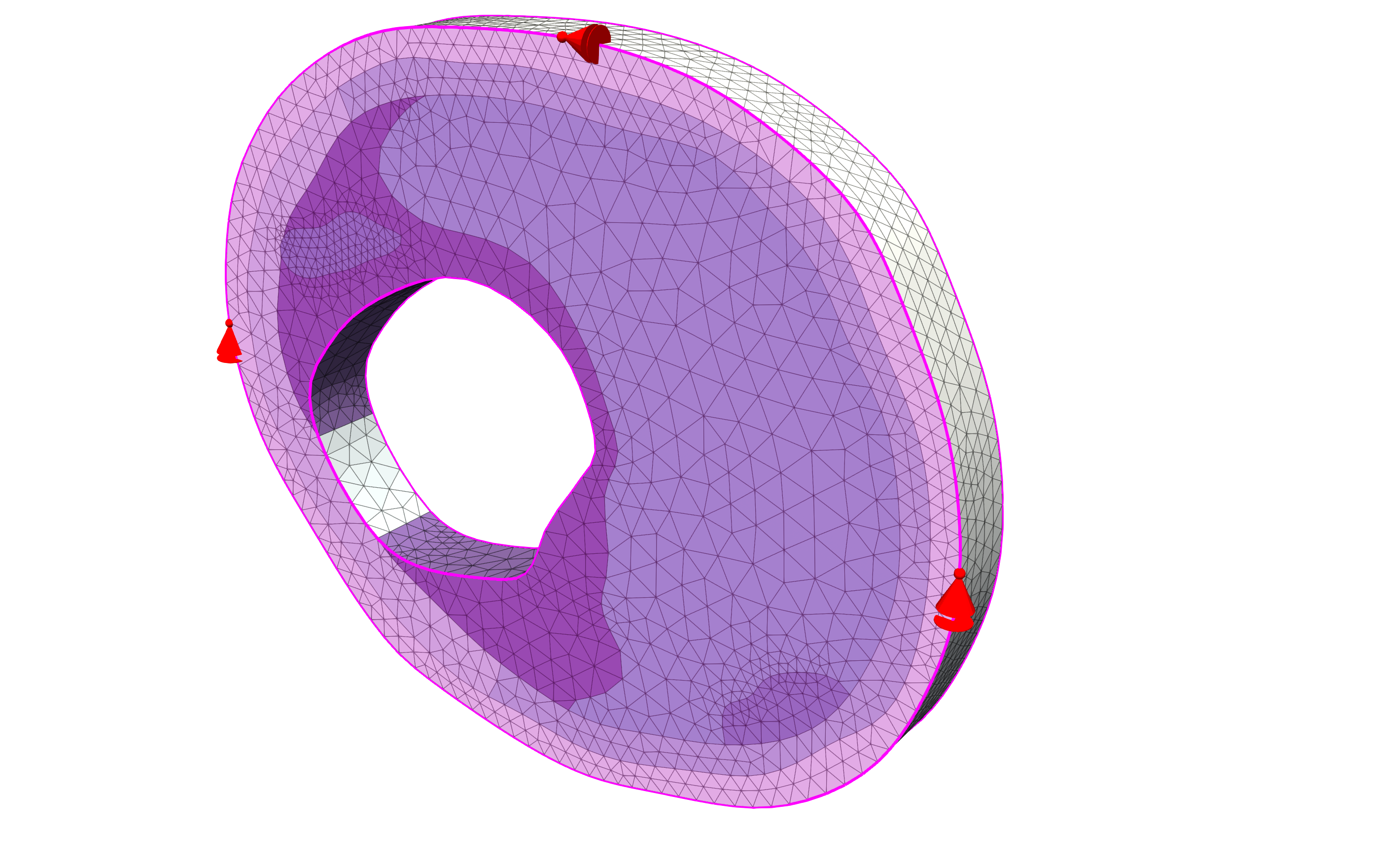}
            \par\vspace{-1.1\baselineskip}\noindent\llap{\footnotesize(f)\hspace{-0.6em}}
            \vspace{0.8\baselineskip}
        \end{minipage}
        \caption{Arterial geometries constructed from histological cross-sections of atherosclerotic human arteries.
        (a) and (d) show coronary thin cap fibroatheromas with moderate stenosis taken from \citet{NarNakVirKolPetNewMalFusFin:2013:hco} and \citet{PhiAndLacLorMakBot:2013:mro}, respectively.
        (b) shows a tibial artery with severe stenosis, calcifications, and a large fibrous cap obtained from \citet{KimGuz:2023:mac}.
        (c) depicts an external iliac artery with eccentric stenosis obtained from \citet{HolSchSta:2000:moa}.
        (e) depicts a thick cap fibrous atheroma with moderate stenosis taken from \citet{PlaHeeDae:2009:pia}.
        Exemplary boundary value problem is shown in (f): three nodes (in red) are fixed in the circumferential direction to prevent rigid body motion, the front and back faces (in pink) are fixed in the axial direction, while pressure is applied to the luminal surface
        \label{fig:arteries_geometry}}
    \end{figure}

    Five arterial geometries were constructed using histological cross-sections of atherosclerotic human arteries available in the literature.
    They were chosen to represent different stages of atherosclerosis, with varying degrees of stenosis and plaque morphology.
    3D geometries were constructed in Autodesk Inventor~\cite{AutodeskInventor} by constant cross-section extrusion of manually segmented regions.
    These geometries were then adjusted for apparent tissue-processing artifacts and meshed using Coreform Cubit~\cite{coreformCubit} \colsr{as well as Gmsh}~\cite{GeuzaineRemacle2009} to produce tetrahedral element meshes \colsr{that were refined locally at the interface and shoulder regions.}
    \colsr{A mesh convergence study, including geometries with rounded fibrous cap shoulders, is presented in Appendix~\ref{sec:app_mesh}.}
    Based on the findings of \citet{HolSomReg:2004:amp}, six distinct regions were considered: adventitia, healthy media, degenerated media, lipid pool, calcifications, and the fibrous cap.
    Fig.~\ref{fig:arteries_geometry} depicts the chosen geometries along with the segmentation of the different regions.
    Fig.~\ref{fig:arteries_geometry}(a) and (d) show coronary thin cap fibroatheromas with moderate stenosis taken from \citet{NarNakVirKolPetNewMalFusFin:2013:hco} and \citet{PhiAndLacLorMakBot:2013:mro}, respectively. Fig.~\ref{fig:arteries_geometry}(b) was obtained from \citet{KimGuz:2023:mac} and shows a tibial artery with severe stenosis, calcifications, and a large fibrous cap.
    Fig.~\ref{fig:arteries_geometry}(c) depicts an external iliac artery with eccentric stenosis obtained from \citet{HolSchSta:2000:moa}.
    Finally, Fig.~\ref{fig:arteries_geometry}(e) depicts a fibrous cap atheroma with moderate stenosis taken from \citet{PlaHeeDae:2009:pia}.
    \subsubsection{Boundary conditions and simulation protocol}\label{sec:boundary_conditions}
    Boundary conditions are designed to mimic in vivo conditions as closely as possible.
    The top and bottom faces of the arterial sections are fixed in the axial direction,
    while three nodes on the outer surface of the adventitia are fixed in the circumferential direction to prevent rigid body motion and to ensure that the system is statically determinate~(see Fig.~\ref{fig:arteries_geometry}(f)).
    The luminal surface is subjected to a mild hypertensive pressure of \colsr{11.33}~kPa, which translates to 85~mmHg, higher (by about 20~mmHg) than the mean arterial pressure experienced by moderate sized arteries in humans~\citep[see Fig.~2.2 in][]{DelTou:2019:tov}.
    Pressure was applied as a deformation-dependent load using surface elements \citep[see][]{Wri:2008:nfe} to ensure load remains normal to the luminal surface.
    
    In the body, amlodipine may reach a point in the arterial wall through two ways depending on the location: near the luminal side through direct diffusion from the bloodstream and the near adventitial side by means of perfusion through the vasa vasorum.
    There is a lack
    of experimental data on tissue/vessel wall concentrations of amlodipine.
    Pharmacokinetic studies typically focus on plasma concentrations and dose-response relationships, which is insufficient for our purpose.
    \citet{MerElli:1992:cpo} found that plasma concentrations of amlodipine in 14-day repeated oral dosing reached 25~ng/mL.
    Assuming a molecular weight of 409.9 g/mol, this corresponds to a concentration of approximately $c = 0.06\,\mathrm{\mu M}$.
    \citet{JepBruNie:1998:ade} determined that in guinea pigs, the myocardial tissue amlodipine accumulation was 600 times that of plasma.
    Furthermore, the large volume distribution of 21 L/kg,
    \colsr{attributed to its strong affinity to cell membranes,} also suggests wide distribution in tissues~\cite{MerElli:1992:cpo}.
    In addition, \citet{WanJiaLiLai:2008:eos} found that concentrations beyond $2.5\,\mathrm{\mu M}$ did not result in any further significant reduction in calcium channel current.
    Consequently, we assumed that tissue amlodipine concentration at the intimal surface (from direct blood contact) and the adventitial surface (through vasa vasorum) to be constant and equal to $c = 2\,\mathrm{\mu M}$.
    
    The artery exists in a complex state where there is an interplay of different phenomena such as residual stresses due to growth and remodeling, stretch induced SMC activation, and pulsatile blood flow, all of which take place at different time scales.
    SMC contractions occur over the course of seconds to minutes, while blood flow pulses are observed even below seconds, whereas growth and remodeling occur over the course of days or weeks.
    In principle, it would be possible to simulate all processes simultaneously, but this would not only require sufficiently low growth and remodeling rates to ensure compatibility with SMC activation but also a fine temporal resolution to capture the fast SMC activation leading to a large number of time steps and thus requiring avoidable computer resources.
    \begin{figure*}[t]
        \centering
             \includegraphics[width=\linewidth]{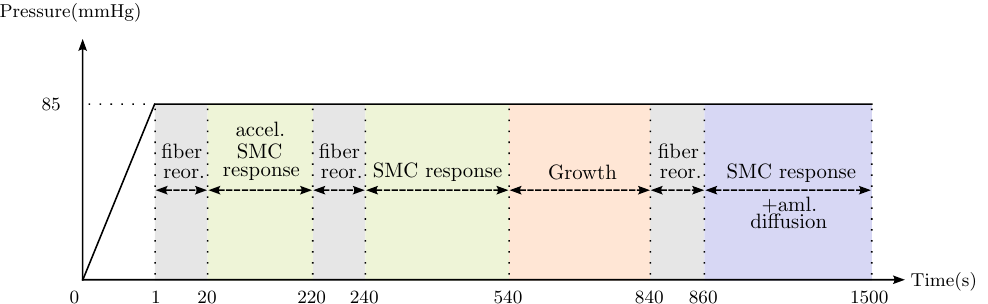}
        \caption{\footnotesize Simulation protocol. 0--1s: loading phase, 1--860s: adaptation phase, 860+s: physiologically relevant phase to be analyzed, virtual drug intake at 860s for effect assessment
        \label{fig:loading_curve}}
    \end{figure*}
    \begin{table}[t]
        \centering
        \caption{\colsr{Upper limits of the growth factors in the three growth directions $i$} \label{tab:growth_parameters}}
        \begin{tabular}{lccc} \toprule
            Parameter & $i=1$ & $i=2$ & $i=3$ \\ \midrule
            $\vartheta^{+}_{(i)}$ & $1.00506$ & $1.0206$ & $1.11992$ \\ \bottomrule
        \end{tabular}
    \end{table}
    We thus adopt the simulation protocol depicted in Fig.~\ref{fig:loading_curve}, where we assume that the growth and remodeling processes are negligible during SMC activation, an assumption supported by their different time scales.
    The simulation is split into three phases: 0--1s loading phase, 1--860s adaptation phase, and 860+s physiological phase.
    Only after the adaptation phase, the artery is in a physiologically relevant state.
    At time zero, we begin with a stress-free reference configuration \colsr{with initial fiber orientations of the two fiber families $\bA^{(f)}$ set to $\pm 30^{\circ}$ to the axial direction in the axial-circumferential plane. Further, it is assumed} that all the cross-bridges in the SMCs are in the free dephosphorylated state.
    Over the course of the adaptation phase, the cross-bridges reach a steady-state distribution, and the artery adapts to the applied pressure through growth and reorientation leading to
    \colsr{one possible} physiological stress state with realistic, distributed fiber arrangement.
    The adaptation phase itself is decomposed in several sub-phases in order to arrive at a \colsr{plausible} physiological state with reduced effort.
    Between 20--540s, SMCs active response is activated to allow cross-bridge cycling and SMC adaptation to the applied pressure.
    During the initial active stage (20--220s) the SMC contraction rate is increased by a factor of 20 (by modifying the rate-dependent parameters) to accelerate the adaptation and reduce the number of time steps required.
    This is followed by a growth phase (540--840s) where stress-driven growth is activated to allow the artery to adapt its geometry and stress state to the applied pressure.
    \colsr{Growth directions are aligned with the principal stress directions at the onset of growth and are subsequently treated as material directions.}
    \colsr{Here, the reference geometries, derived from the histological cross-sections, are considered stress-free, which, in reality, is not the case~\cite{HolSomAueRegOgd:2007:l3r,Sok:2019:rdo}.
    Since layer-specific opening angle data was unavailable for these specimens, the residual stresses have been approximated by the growth-based formulation described in Section~\ref{sec:residual_stresses}, as true residual stresses are unknown.
    }
    The fiber reorientation algorithm is employed after each sub-phase to ensure realistic fiber directions even for complex geometries and stress states.
    The growth and reorientation rates are chosen sufficiently high so that there are no discernible changes after the respective time intervals.
    \colsr{The reorientation rate $k_{\eta}^{+}$ is chosen to be $0.001$s$^{-1}$ and the curvature parameter $m_{\eta}^{+}$ is set to $5.0$.
    The upper limits of the growth factors $\vartheta^{+}_{(i)}$ are listed in Table~\ref{tab:growth_parameters}, whereas the lower limits $\vartheta^{-}_{(i)}$ and the shape parameters $m^{+/-}_{\vartheta,(i)}$ are set to $0.98$ and $3.0$, respectively, in all three growth directions.
    The growth rates $k^{+/-}_{\vartheta,(i)}$ are chosen to be $10^{-5}\,\text{kPa}^{-1}\text{s}^{-1}$ for the arteries in Fig.~\ref{fig:arteries_geometry}(a), (c) and (e), and set to $2.5 \cdot 10^{-6}\,\text{kPa}^{-1}\text{s}^{-1}$ and $10^{-6}\,\text{kPa}^{-1}\text{s}^{-1}$ for those in Fig.~\ref{fig:arteries_geometry}(b) and (d), respectively.}
    After the adaptation phase, i.e., at 860s, when a \colsr{physiologically plausible} state is reached, drug treatment is virtually considered by activating the diffusion process.
    The results are then compared with a control case without drug diffusion to assess the effect of amlodipine on the stress state in the arterial wall.

\section{Results}

    \subsection{External iliac artery}
\begin{figure*}[t]
        \centering
        \begin{minipage}[b]{0.28\linewidth}
            \vspace{0pt}
            \includegraphics[width=\linewidth,page=1]{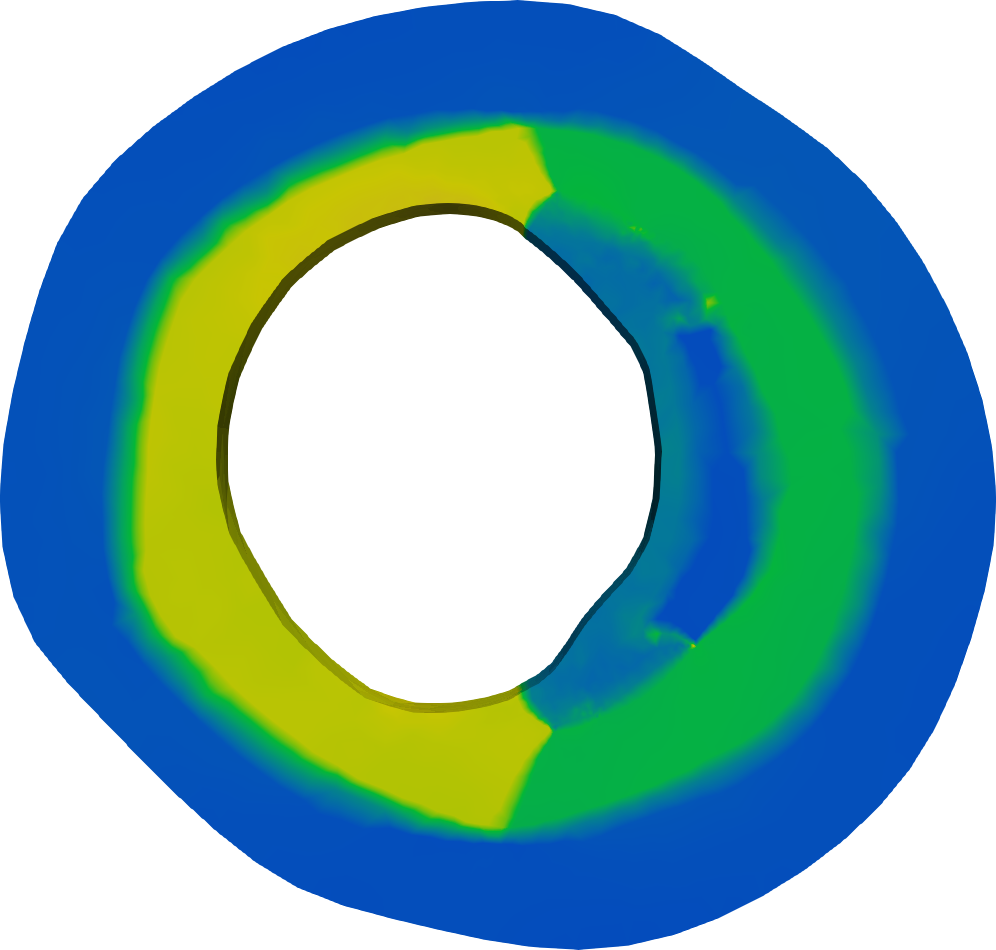}
            \par\vspace{-1.6\baselineskip}\noindent\llap{\footnotesize(a)\hspace{-0.6em}}
            \vspace{0.8\baselineskip}
        \end{minipage}
        \begin{minipage}[b]{0.05\linewidth}
            \vspace{0pt}
            \includegraphics[width=\linewidth,page=1]{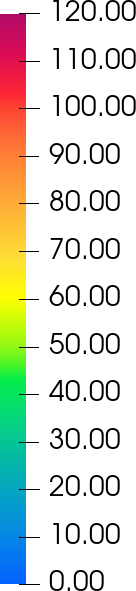}
            \vspace{0.8\baselineskip}
        \end{minipage}
        \hfill
        \begin{minipage}[b]{0.28\linewidth}
            \vspace{0pt}
            \includegraphics[width=\linewidth,page=1]{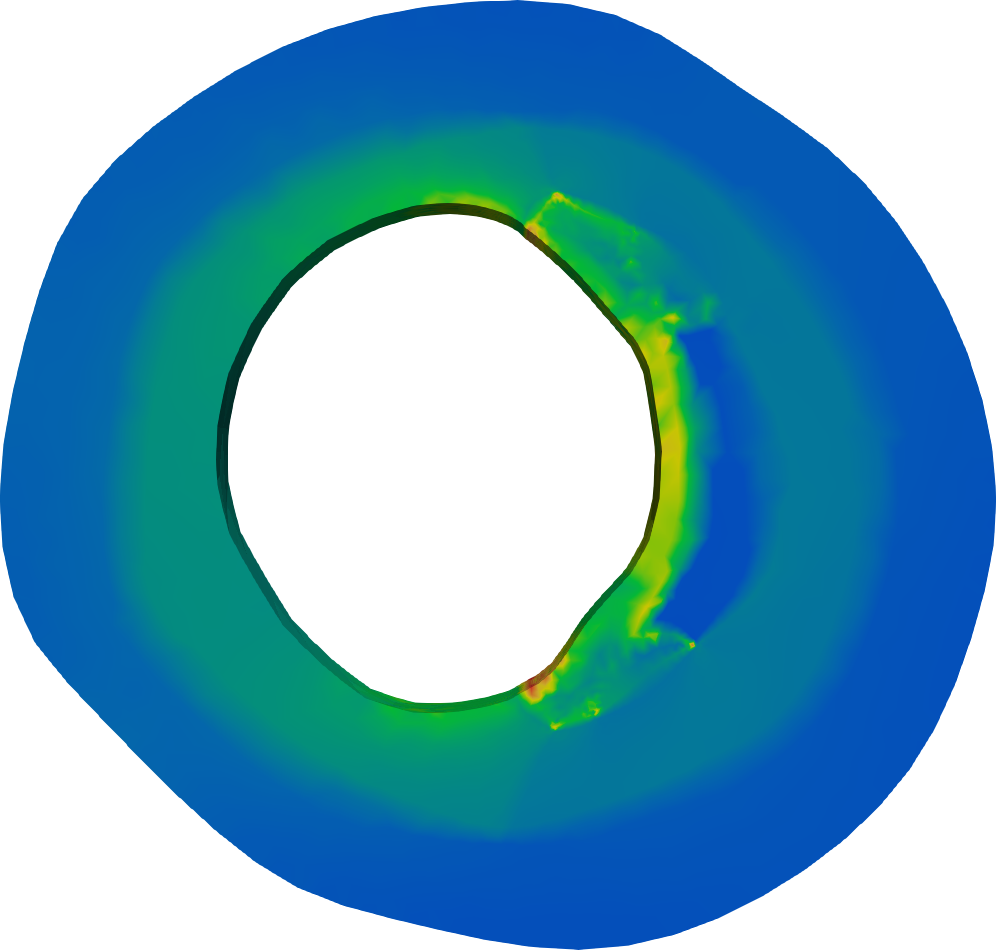}
            \par\vspace{-1.6\baselineskip}\noindent\llap{\footnotesize(b)\hspace{-0.6em}}
            \vspace{0.8\baselineskip}
        \end{minipage}
        \hfill
        \begin{minipage}[b]{0.28\linewidth}
            \vspace{0pt}
            \includegraphics[width=\linewidth,page=1]{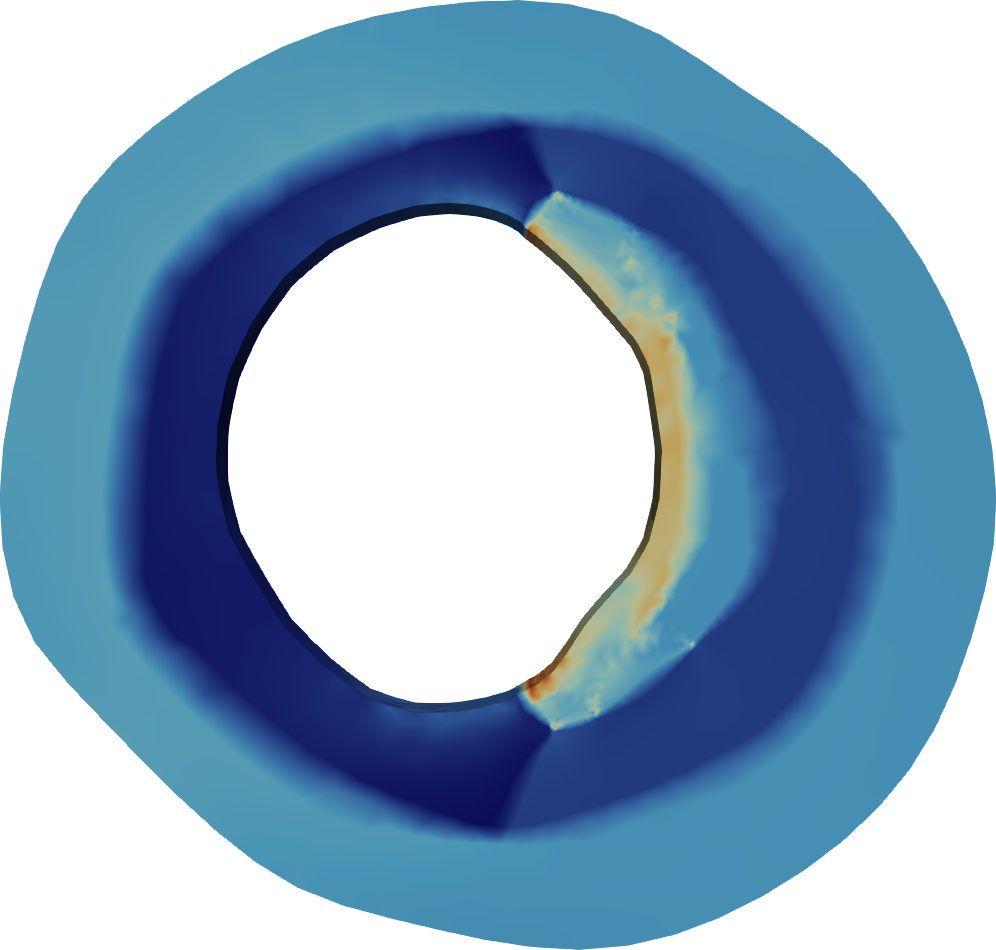}
            \par\vspace{-1.6\baselineskip}\noindent\llap{\footnotesize(c)\hspace{-0.6em}}
            \vspace{0.8\baselineskip}
        \end{minipage}
        \begin{minipage}[b]{0.05\linewidth}
            \vspace{0pt}
            \includegraphics[width=\linewidth,page=1]{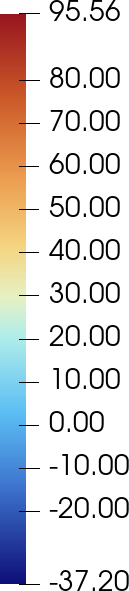}
            \vspace{0.8\baselineskip}
        \end{minipage}
        \hfill
        \caption{Contour plots of von Mises stress~(kPa) in the external iliac artery inflated to 85mmHg.
        (a) shows \colsr{control} stress levels, (b) the stress levels \colsr{with} amlodipine treatment, and (c) shows the difference in von Mises stress with and without amlodipine ($\sigma_\text{aml.} - \sigma_\text{control}$).
        \colsr{As can be seen in (c), in the simulation with drug, von Mises stresses in the medial layers are lower than in the control, whereas those in the fibrous cap are considerably higher}}
        \label{fig:iliac_artery}
    \end{figure*}

    Fig.~\ref{fig:iliac_artery}(a) depicts the von Mises stress distribution in the external iliac artery~(Fig.~\ref{fig:arteries_geometry}(c)) simulated according to the protocol specified in Fig.~\ref{fig:loading_curve} with material and transport parameters defined in Tables~\ref{tab:material_parameters}, \ref{tab:transport_parameters}, and \ref{tab:drug_interaction_parameters}.
    \colsr{In the following, stress values are reported as the volume-weighted median of each component, with the 5th--95th percentile range in parentheses.}
    \colsr{Volume-weighted percentiles are used instead of peak values, since the peak values occur at the sharp edges of the histology-based geometries and increase with mesh refinement, whereas the medians and 5th--95th percentiles are practically mesh-independent~(Appendix~\ref{sec:app_mesh}).}
    \colsr{The median von Mises stress values are $57.6$~kPa~(5th--95th percentile: $41.1$--$60.2$~kPa) in the healthy media, $39.0$~kPa~($27.7$--$40.5$~kPa) in the degenerated media, and $15.2$~kPa~($11.8$--$31.3$~kPa) in the fibrous cap.}
    These values are commensurate with the measured physiological levels reported in the 
   existing literature~\cite{GiuSprWilKhi:2023:tlc,SomRegKolHol:2010:bmp} and correspond well with experimental observations indicating that the media constitutes the most significant load-bearing layer in the physiological state~\cite{HolSomReg:2004:amp,HolGas:2007:csd}.
    \begin{figure*}[t]
    \centering
    \begin{minipage}[b]{0.33\linewidth}
        \vspace{0pt}
        \includegraphics[width=\linewidth,page=1]{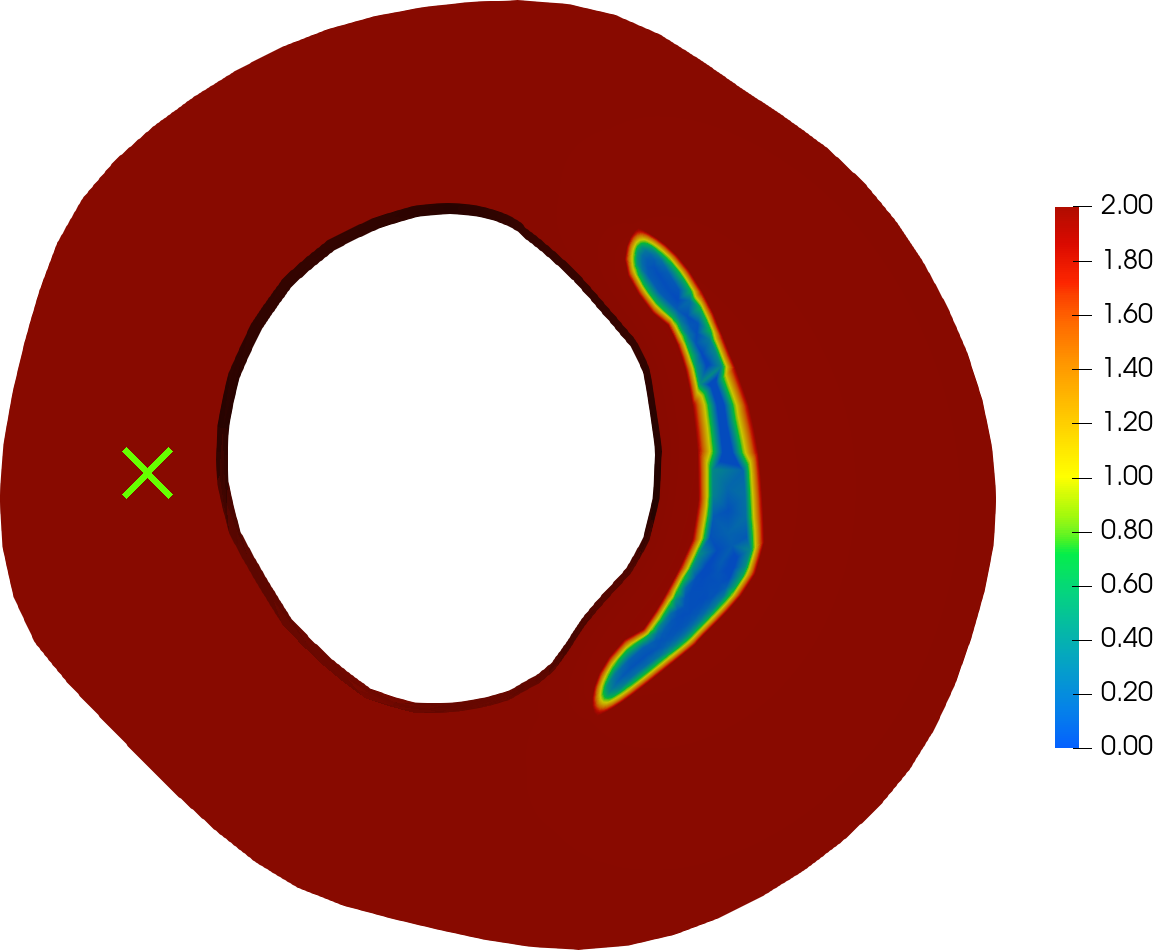}
        \par\vspace{-0.75\baselineskip}\noindent\llap{\footnotesize(a)\hspace{-0.6em}}
        \vspace{0.8\baselineskip}
    \end{minipage}
    \hspace{0.06\linewidth}
    \begin{minipage}[b]{0.47\linewidth}
        \vspace{0pt}
        \includegraphics[width=\linewidth,page=1]{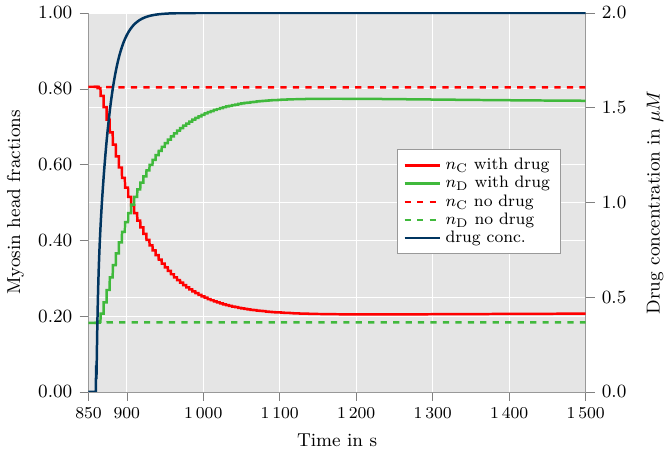}
        \par\vspace{-1.6\baselineskip}\noindent\llap{\footnotesize(b)\hspace{-0.6em}}
        \vspace{0.8\baselineskip}
    \end{minipage}
    \caption{Amlodipine diffusion profile and its effects on myosin head fractions in the external iliac artery.
    (a) depicts the steady-state amlodipine concentration ($\mathrm{\mu M}$) profile.
    (b) shows {the} evolution of myosin head fractions (dashed: no drug; solid: with drug) for the attached, phosphorylated state $n_\text{C}$(red) and attached, dephosphorylated state $n_\text{D}$(green) at the material point marked in (a) by a cross.
    Corresponding drug concentration profile plotted in blue.
    Post-treatment, $n_\text{C}$ reduces from 0.82 to 0.21, significantly reducing contraction 
    }
    \label{fig:drug_effects}
    \end{figure*}
    Here, the high fraction of attached, phosphorylated myosin heads~(see Fig.~\ref{fig:drug_effects}(b)) in the media contributes to increased contraction and thereby elevated stress.
    Fig.~\ref{fig:drug_effects}(a) depicts the amlodipine concentration in the artery at steady-state for the simulation with drug effects.
    For the specified parameters, robust drug penetration into all components except the necrotic core and calcifications can be observed.
    The effect of amlodipine on the myosin head fractions in the attached state can be seen in Fig.~\ref{fig:drug_effects}(b), where the evolution of the values over time are plotted for the material point marked by a cross in Fig.~\ref{fig:drug_effects}(a).
    With increasing drug concentration, the fraction of attached, phosphorylated myosin heads~($n_\text{C}$) reduces from 0.82 to 0.21, resulting in a significant reduction in smooth muscle activation and hence contractile force.
    \colsr{Consequently, stresses in both the healthy and degenerated media are lower than in the control~(see Fig.~\ref{fig:iliac_artery}(a) and (b)).}
    In stark contrast to the medial layers, \colsr{the stresses in the fibrous cap are higher than in the control by a median of $29$~kPa~($8$--$37$~kPa), with local peaks reaching up to about $95$~kPa.}
    In Fig.~\ref{fig:iliac_artery}(c), which illustrates the difference between von Mises stress with and without drug effects, it can be seen \colsr{that the neck areas of the fibrous cap experience the largest stress increase.

    To analyze the impact of several influencing parameters, including the SMC contractility, the drug availability and response, and axial strains, further investigations on the external iliac artery are given in Appendix~\ref{sec:analysis}.
    It is found that the drug-induced increase in fibrous cap stress is smaller
    \begin{itemize}
     \item for decreased contractility (cf.~\ref{sec:NO_effects}), either by
     \begin{itemize}
       \item increased NO-based reduction in SMC contractility,
       \item reduced basal tone (reduced maximum rate of phosphorylation),
       \item or, rather counterintuitively, increased degree of disease in the degenerated media, however, in this case the absolute magnitude of fibrous cap stress after drug-intake is larger and thus, potentially more problematic,
     \end{itemize}
     \item for lower drug dose and less sensitive drug response~(cf.~\ref{sec:transport_param_effects}),
     \item for decreased tissue diffusivity and increasing drug absorption rate (cf.~\ref{sec:transport_param_effects}),
     \item for increased axial strain (cf.~\ref{sec:axial_stretch})
    \end{itemize}
    However, although all parameter changes lead to changes in the drug-induced increase in fibrous cap stress, the stress gains remain relevant in all scenarios.
    } %
    \subsection{Role of plaque morphology}
\begin{figure*}[t]
    \centering
    \begin{minipage}[b]{0.235\linewidth}
        \vspace{0pt}
        \includegraphics[width=\linewidth,page=1]{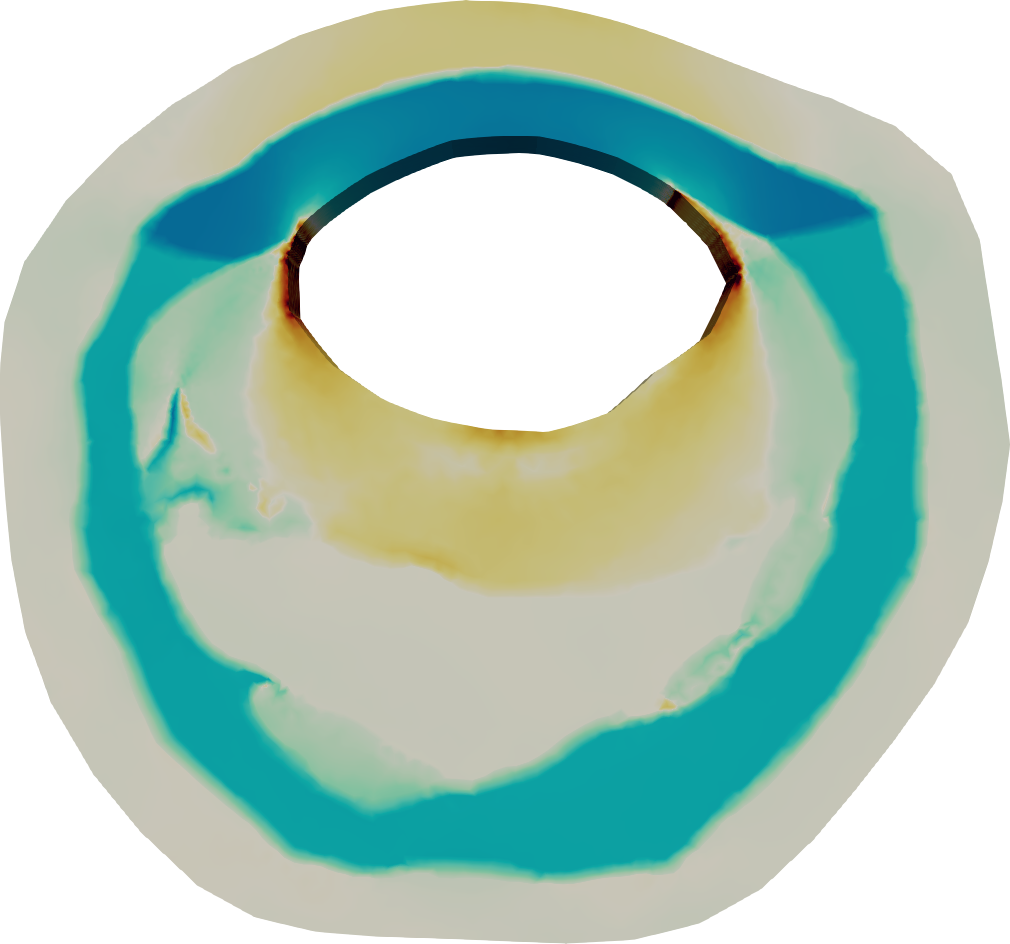}
        \par\vspace{-1.1\baselineskip}\noindent\llap{\footnotesize(a)\hspace{-0.6em}}
        \vspace{0.8\baselineskip} %
    \end{minipage}
    \hfill
    \begin{minipage}[b]{0.24\linewidth}
        \vspace{0pt}
        \includegraphics[width=\linewidth,page=1]{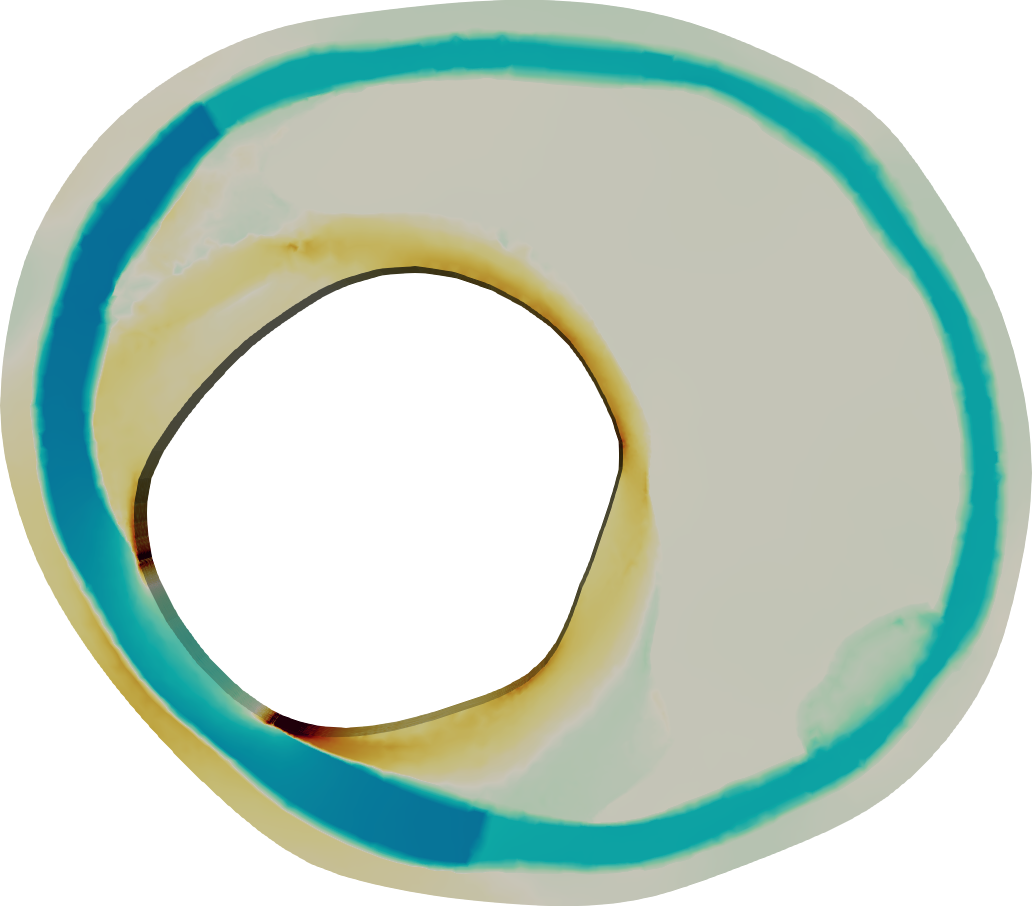}
        \par\vspace{-1.1\baselineskip}\noindent\llap{\footnotesize(b)\hspace{-0.6em}}
        \vspace{0.8\baselineskip}
    \end{minipage}
    \hfill
    \begin{minipage}[b]{0.21\linewidth}
        \vspace{0pt}
        \includegraphics[width=\linewidth,page=1]{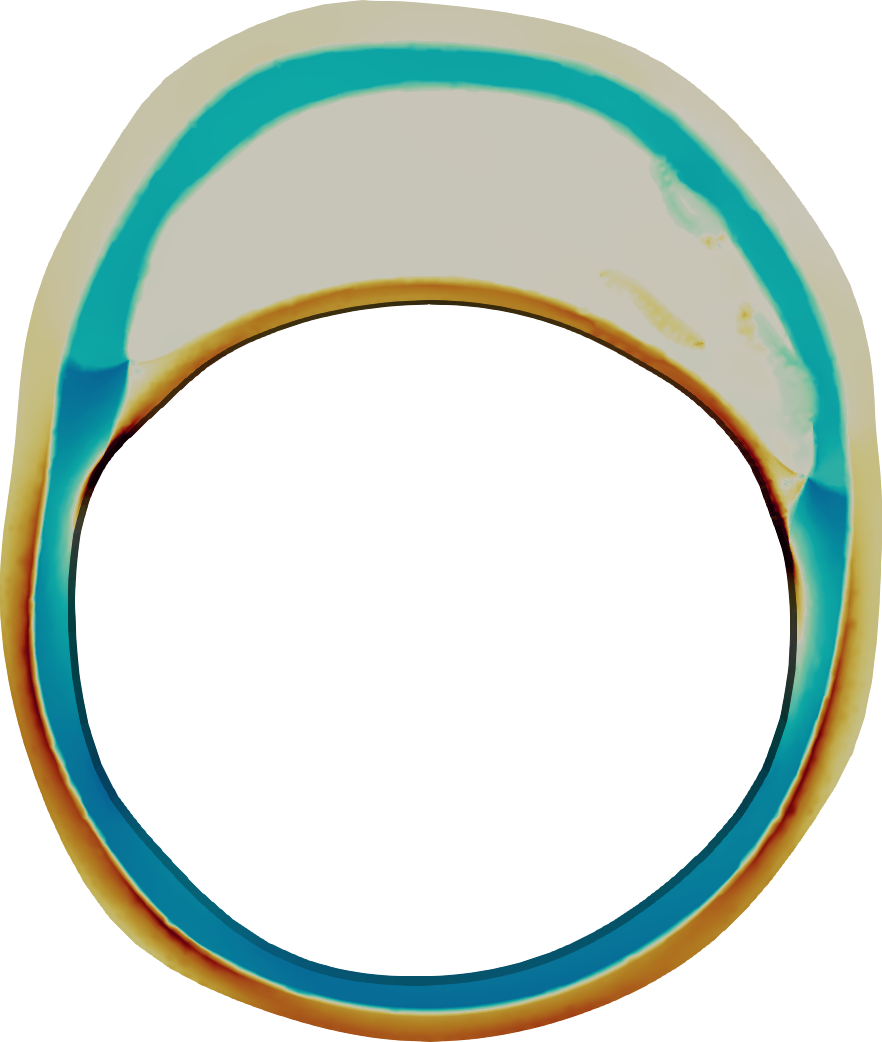}
        \par\vspace{-1.1\baselineskip}\noindent\llap{\footnotesize(c)\hspace{-0.6em}}
        \vspace{0.8\baselineskip}
    \end{minipage}
    \hfill
    \begin{minipage}[b]{0.21\linewidth}
        \vspace{0pt}
        \includegraphics[width=\linewidth,page=1]{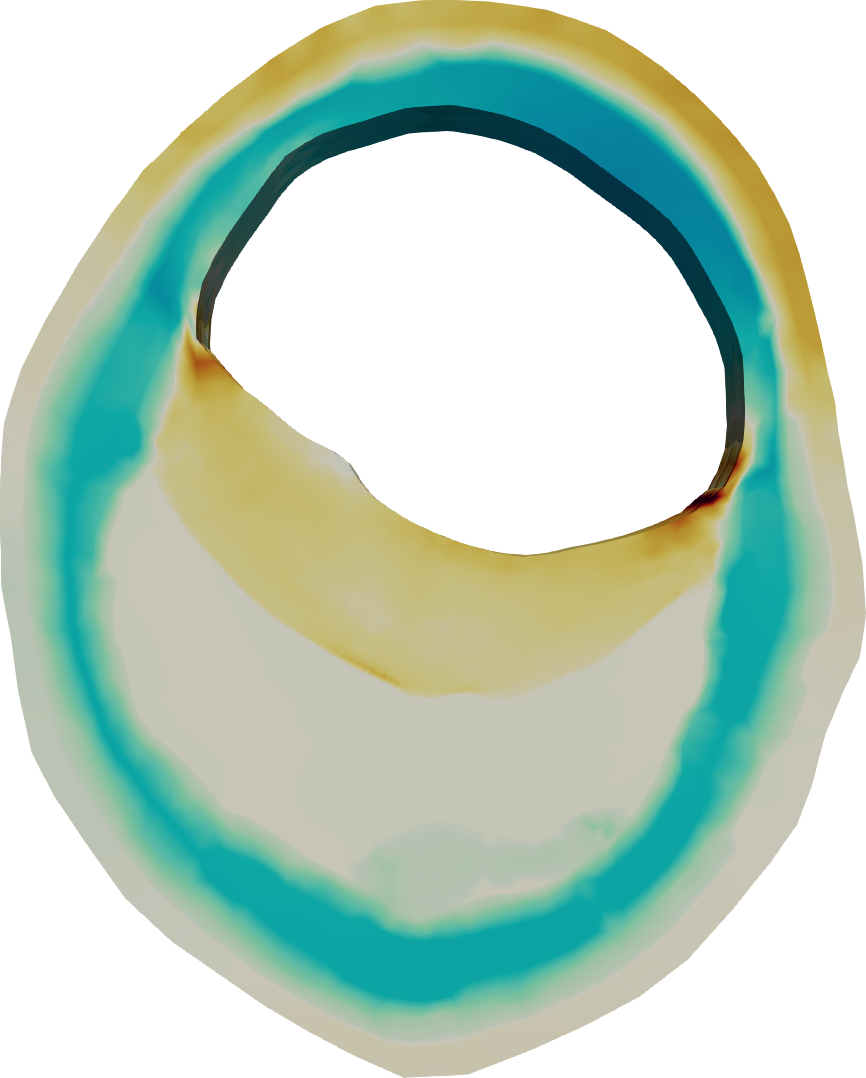}
        \par\vspace{-1.1\baselineskip}\noindent\llap{\footnotesize(d)\hspace{-0.6em}}
        \vspace{0.8\baselineskip}
    \end{minipage}
    \hfill
    \begin{minipage}[b]{0.045\linewidth}
        \vspace{0pt}
        \includegraphics[width=\linewidth,page=1]{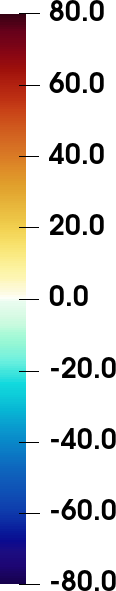}
        \vspace{0.1\baselineskip}
    \end{minipage}
    \caption{%
        Change in von Mises stress in kPa post-drug treatment for different plaque morphologies~(see Fig.~\ref{fig:arteries_geometry}), represented by the difference between post-drug and baseline von Mises stress distributions ($\sigma_\text{aml.} - \sigma_\text{baseline}$).
        \colsr{The fib. cap stress increases by a volume-weighted median of $5.4$~kPa (5th--95th percentile: $-8.7$ to $13.9$~kPa) in the tibial artery in (a), with the highest increases at high curvature neck regions, and by $3.9$~kPa ($-5.6$ to $18.1$~kPa) in the TCFA in (b).
        The fib. cap of the TCFA depicted in (c) shows the largest increase, with a median of $28.4$~kPa ($2.2$ to $44.3$~kPa) and about $70$~kPa or more in $1\%$ of its volume at the neck regions.
        The moderately stenosed artery in (d) has a median fib. cap stress increase of $7.4$~kPa ($3.2$ to $14.1$~kPa), with some adventitial areas seeing an increase of up to about $20$~kPa.
        The stress reduction in the healthy media is similar in all arteries~(median $-21.4$ to $-33.8$~kPa).}
        Summarizing, as a result of drug-induced relaxation of smooth muscle, an unexpected increase of stress in parts of the atherosclerotic plaque is observed in all arterial morphologies
    }
    \label{fig:plaque0}
    \end{figure*}

    Plaque morphology plays an important role in the stress distribution within atherosclerotic arteries.
    To investigate the effect of different plaque types on the change in stress levels post-drug treatment, we simulated four additional artery geometries depicted in Fig.~\ref{fig:arteries_geometry}(b,c,d,e), while maintaining the same material and transport parameters as in the external iliac artery case.
    In all simulations, the region around the interface between the fibrous cap and the media exhibits the highest stress values.
    This phenomenon may be attributed to numerical challenges linked with the
    sharp edges resulting from the
    histology-based geometry creation~\cite{ChaChaShiMacIftTeaKamBouKaa:2004:mao}.
    \colsr{These peak values are confined to the immediate vicinity of the sharp edges and do not affect the volume-weighted medians reported here~(Appendix~\ref{sec:app_mesh}).}
    However, the main interest in our study lies in the change in stresses post-drug treatment.
    To this end, Fig.~\ref{fig:plaque0} shows these changes in von Mises stress comparing the states before and after simulated drug intake.
    \colsr{In all cases, the stress in most of the fibrous cap increases post-drug treatment, with areas of high curvature and neck regions experiencing the highest increases.}
    The medial layers exhibit reduced stresses due to drug-induced reduction in SMC contraction.
    \colsr{The largest increase occurs in the thin cap fibrous atheroma~(TCFA) in Fig.~\ref{fig:plaque0}(c), with a median of $28.4$~kPa (5th--95th percentile: $2.2$ to $44.3$~kPa), comparable to the external iliac artery~(median $28.9$~kPa), and an increase of about $70$~kPa or more in $1\%$ of the fibrous cap volume~(99th percentile $69.6$~kPa).
    This TCFA also carries by far the highest fibrous cap stress, with a median of $121.3$~kPa without and $150.1$~kPa with amlodipine.
    In the other three arteries, the median increase in the fibrous cap is considerably smaller, with $5.4$~kPa ($-8.7$ to $13.9$~kPa) in the tibial artery~(Fig.~\ref{fig:plaque0}(a)), $3.9$~kPa ($-5.6$ to $18.1$~kPa) in the second TCFA~(Fig.~\ref{fig:plaque0}(b)), and $7.4$~kPa ($3.2$ to $14.1$~kPa) in the moderately stenosed artery~(Fig.~\ref{fig:plaque0}(d)).
    Hence, a thin fibrous cap alone does not determine the magnitude of the drug-induced stress increase.} %
    \colsr{Arterial morphologies with thinner walls~(Fig.~\ref{fig:plaque0}(c) and (d)) show some adventitial stress increases as well~(95th percentile $40.9$~kPa and $15.8$~kPa, respectively).}
    \colsr{The median stress reduction in the healthy media ranges from $21.4$~kPa~(Fig.~\ref{fig:plaque0}(d)) to $33.8$~kPa~(Fig.~\ref{fig:plaque0}(a)), compared to $29.3$~kPa in the external iliac artery.} %
    \subsection{Influence of blood pressure reduction}\label{sec:bp_reduction}
\begin{figure*}
    \centering
    \resizebox{0.928\textwidth}{!}{%
    \begin{minipage}{\textwidth}
    \begin{minipage}[b]{0.27\linewidth}
        \includegraphics[width=\linewidth,page=1]{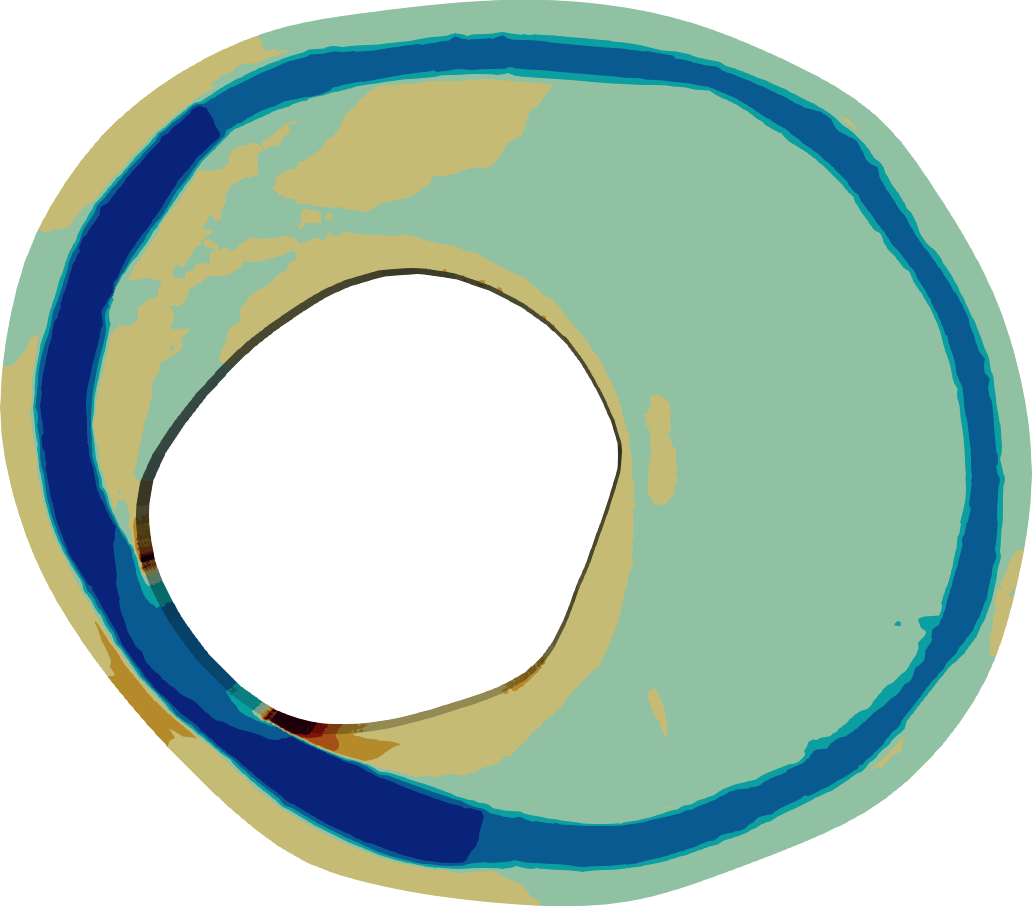}
        \par\vspace{-1.1\baselineskip}\noindent\llap{\footnotesize(a-1)\hspace{-0.6em}}
    \end{minipage}
    \hfill
    \begin{minipage}[b]{0.27\linewidth}
        \includegraphics[width=\linewidth,page=1]{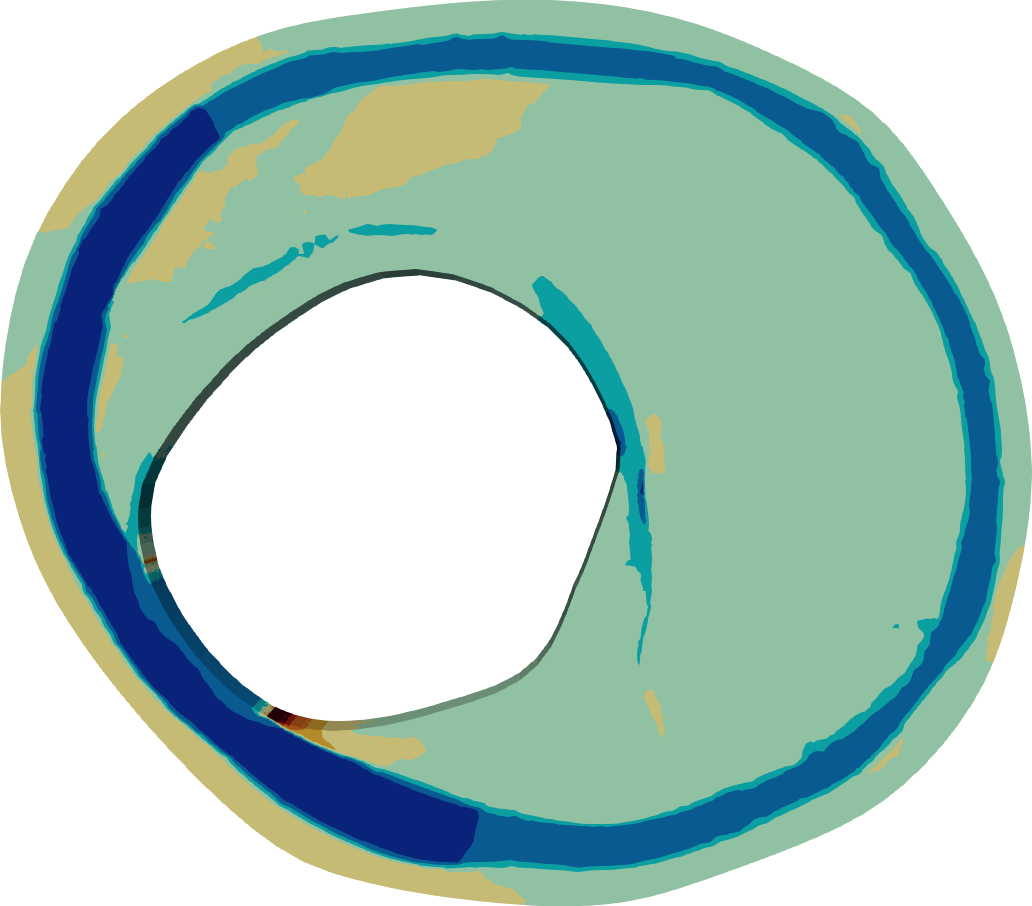}
        \par\vspace{-1.1\baselineskip}\noindent\llap{\footnotesize(a-2)\hspace{-0.6em}}
    \end{minipage}
    \hfill
    \begin{minipage}[b]{0.27\linewidth}
        \includegraphics[width=\linewidth,page=1]{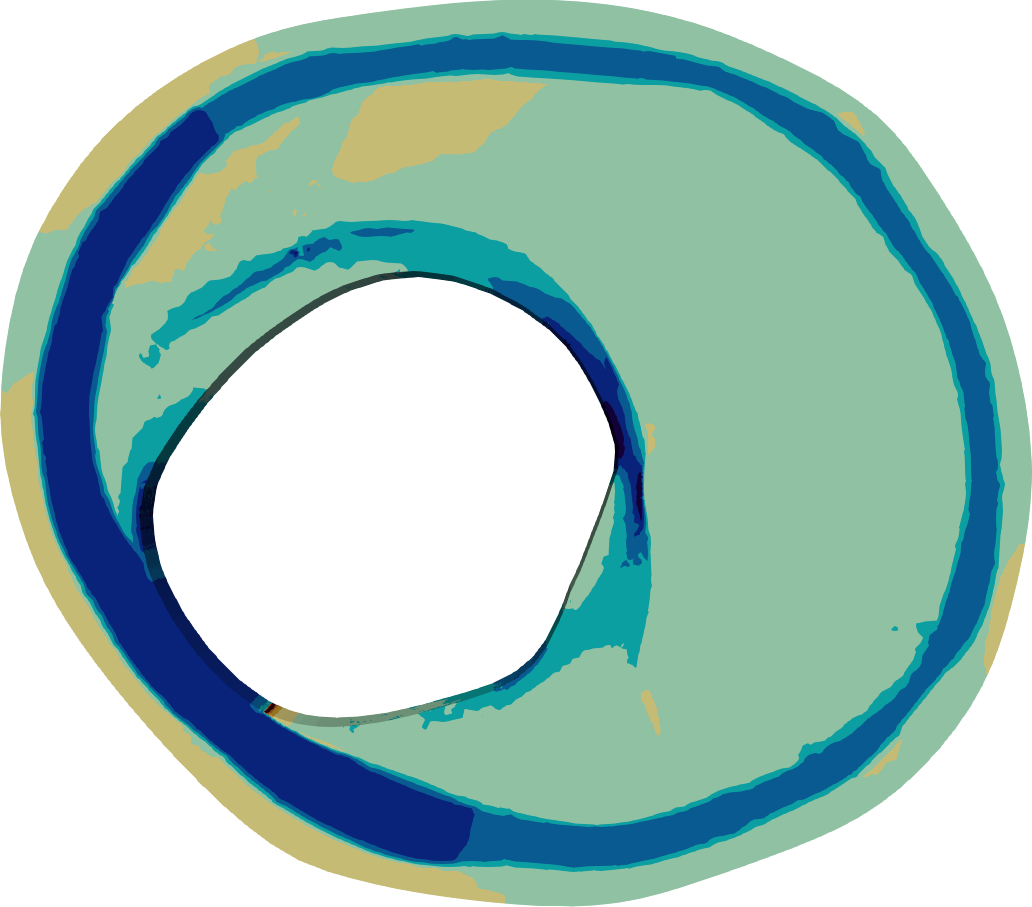}
        \par\vspace{-1.1\baselineskip}\noindent\llap{\footnotesize(a-3)\hspace{-0.6em}}
    \end{minipage}
    \hfill
    \begin{minipage}[b]{0.04\linewidth}
        \includegraphics[width=\linewidth,page=1]{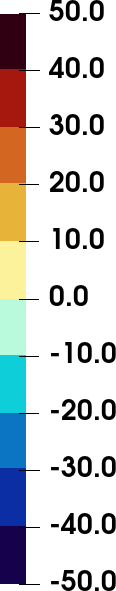}
    \end{minipage}
    \par
    \begin{minipage}[b]{0.26\linewidth}
        \includegraphics[width=\linewidth,page=1]{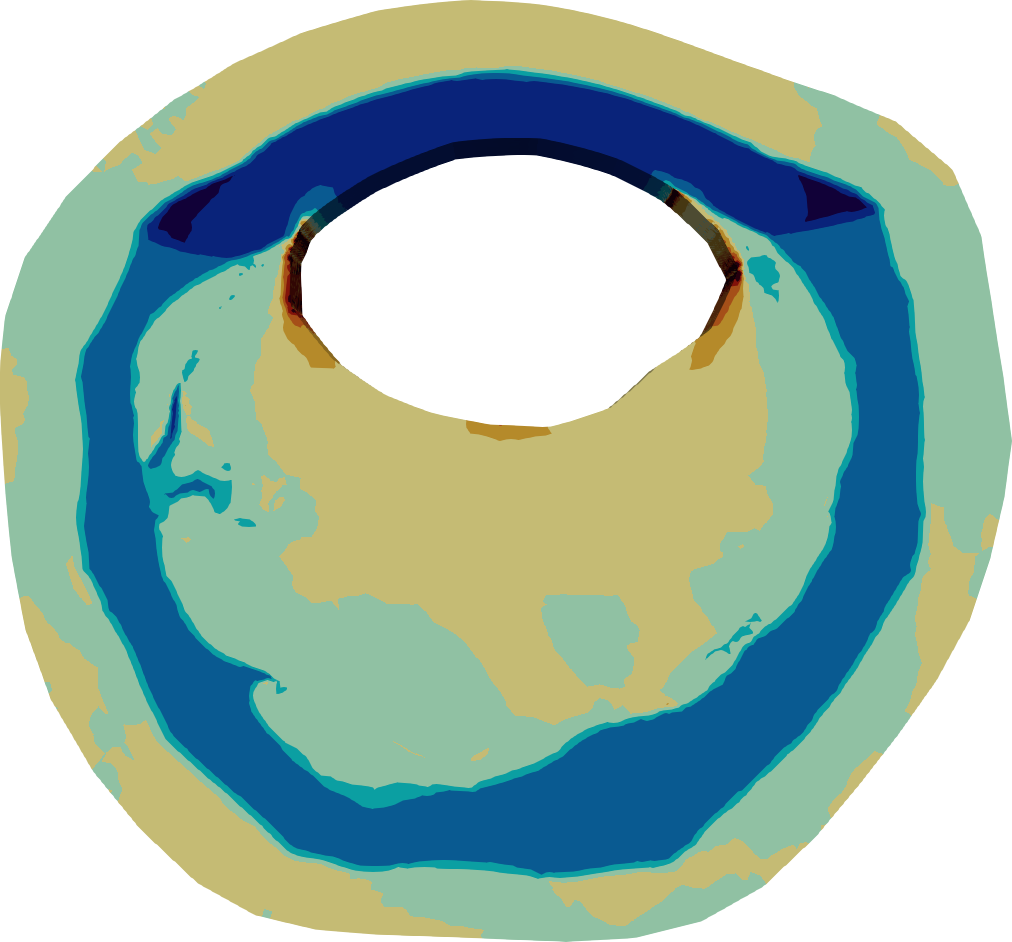}
        \par\vspace{-1.1\baselineskip}\noindent\llap{\footnotesize(b-1)\hspace{-0.6em}}
    \end{minipage}
    \hfill
    \begin{minipage}[b]{0.26\linewidth}
        \includegraphics[width=\linewidth,page=1]{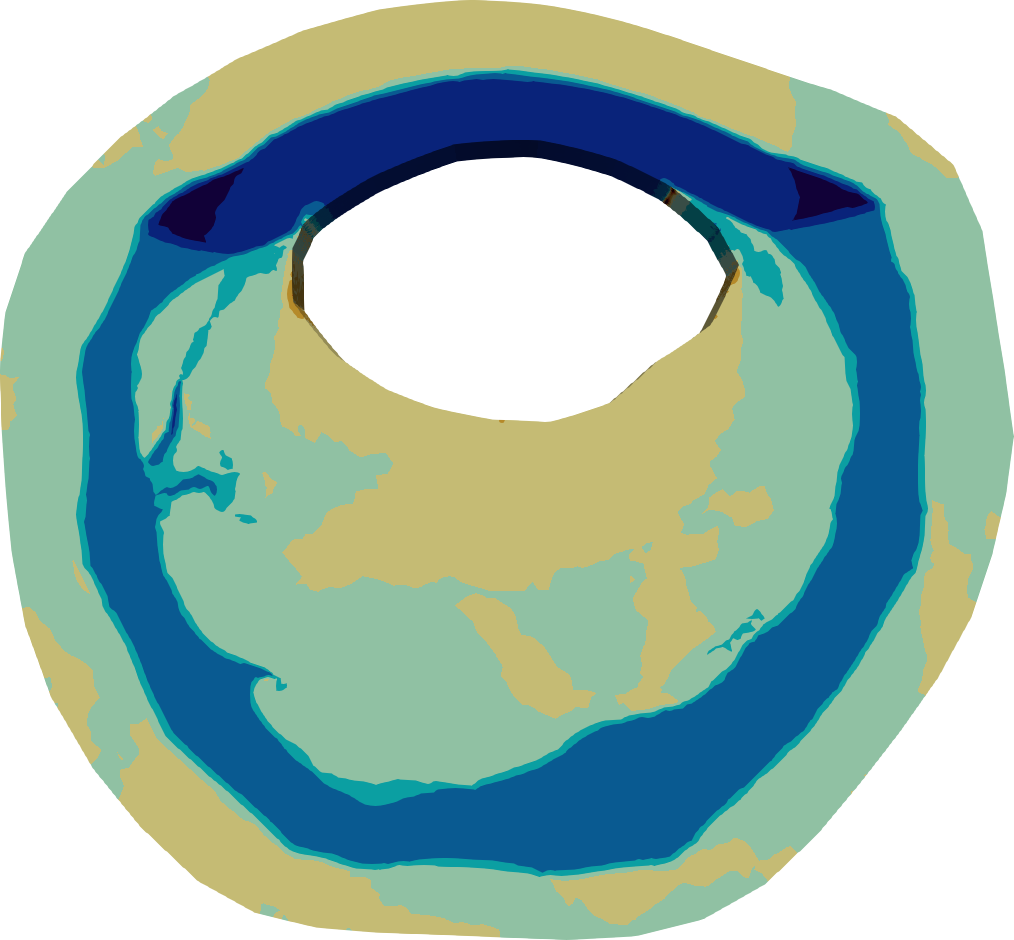}
        \par\vspace{-1.1\baselineskip}\noindent\llap{\footnotesize(b-2)\hspace{-0.6em}}
    \end{minipage}
    \hfill
    \begin{minipage}[b]{0.26\linewidth}
        \includegraphics[width=\linewidth,page=1]{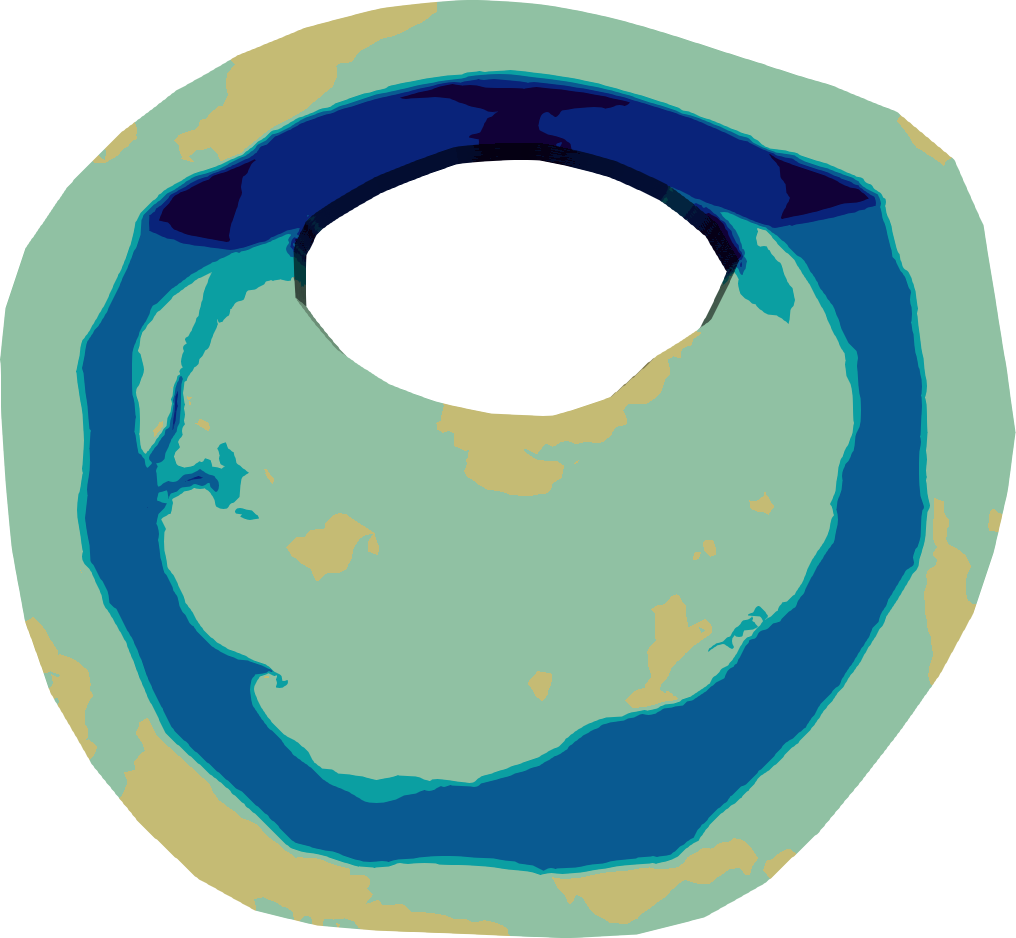}
        \par\vspace{-1.1\baselineskip}\noindent\llap{\footnotesize(b-3)\hspace{-0.6em}}
    \end{minipage}
    \hfill
    \begin{minipage}[b]{0.04\linewidth}
        \includegraphics[width=\linewidth,page=1]{legend50.png}
    \end{minipage}
    \par
    \begin{minipage}[b]{0.27\linewidth}
        \includegraphics[width=\linewidth,page=1]{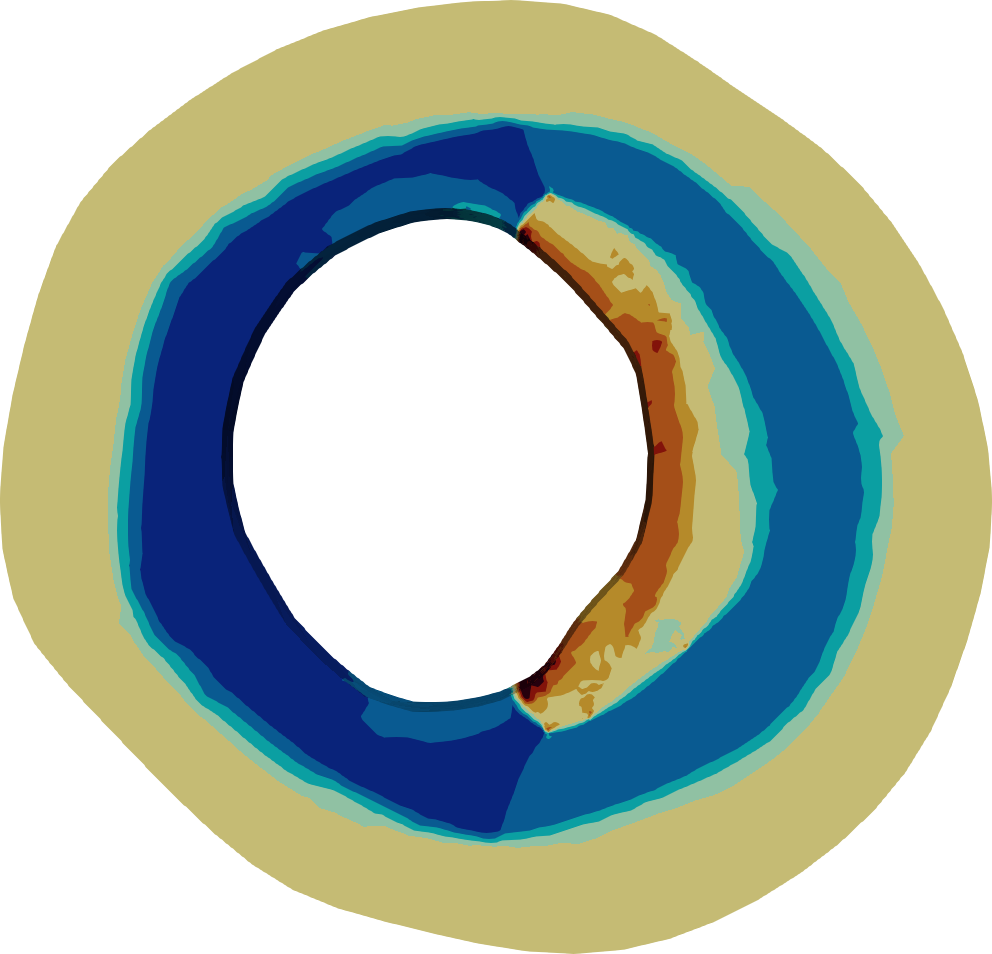}
        \par\vspace{-1.1\baselineskip}\noindent\llap{\footnotesize(c-1)\hspace{-0.6em}}
    \end{minipage}
    \hfill
    \begin{minipage}[b]{0.27\linewidth}
        \includegraphics[width=\linewidth,page=1]{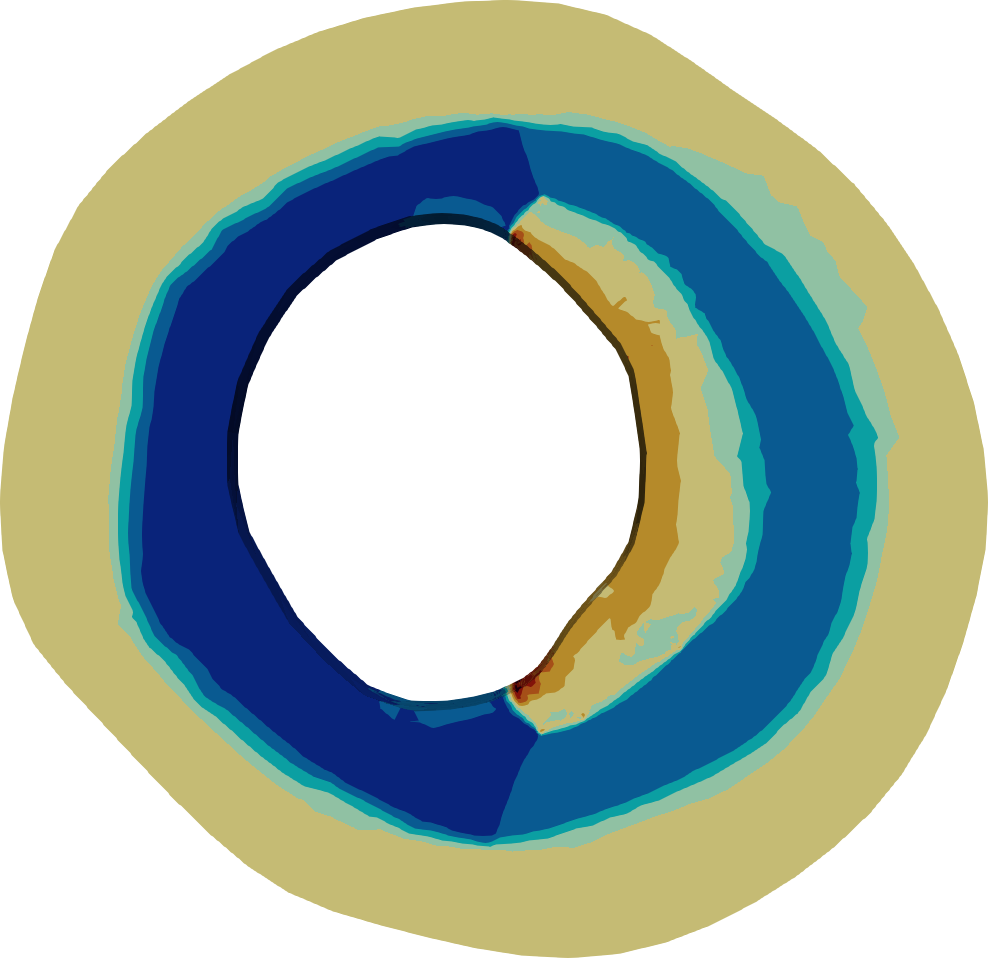}
        \par\vspace{-1.1\baselineskip}\noindent\llap{\footnotesize(c-2)\hspace{-0.6em}}
    \end{minipage}
    \hfill
    \begin{minipage}[b]{0.27\linewidth}
        \includegraphics[width=\linewidth,page=1]{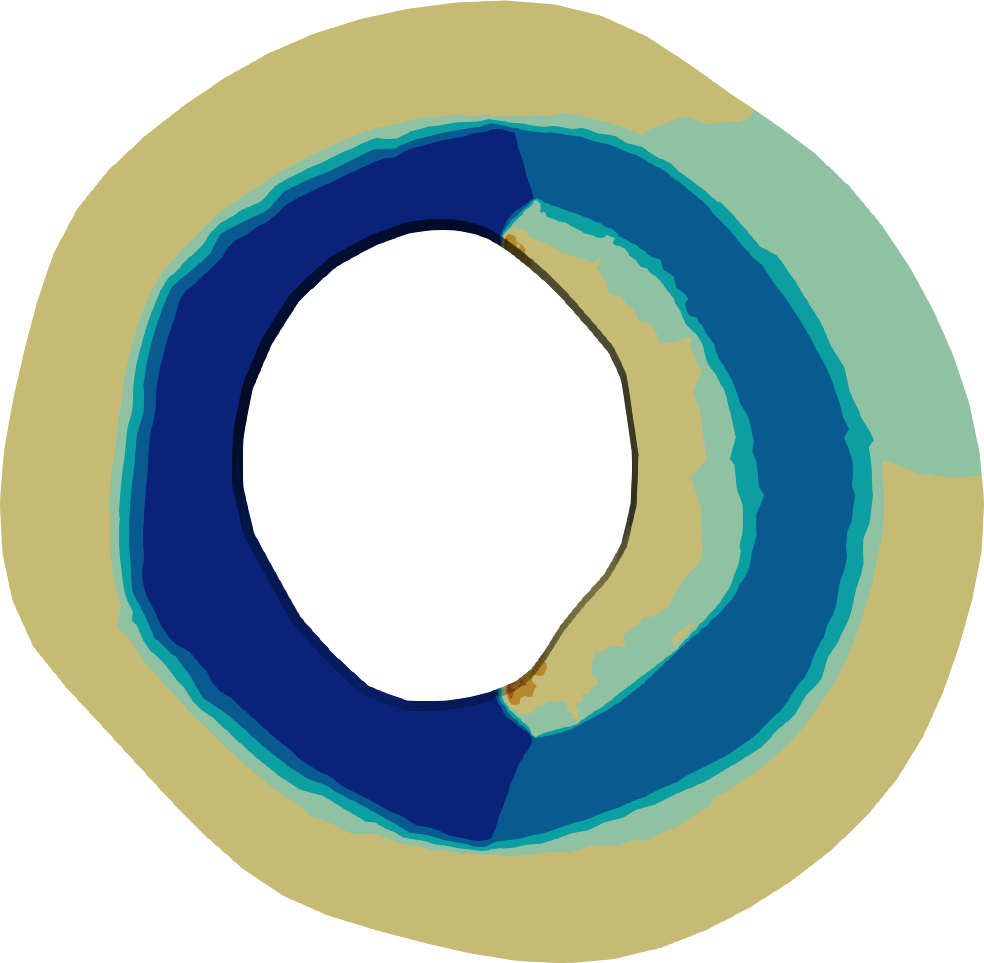}
        \par\vspace{-1.1\baselineskip}\noindent\llap{\footnotesize(c-3)\hspace{-0.6em}}
    \end{minipage}
    \hfill
    \begin{minipage}[b]{0.04\linewidth}
        \includegraphics[width=\linewidth,page=1]{legend50.png}
    \end{minipage}
    \par\medskip
    \begin{minipage}[b]{0.235\linewidth}
        \includegraphics[width=\linewidth,page=1]{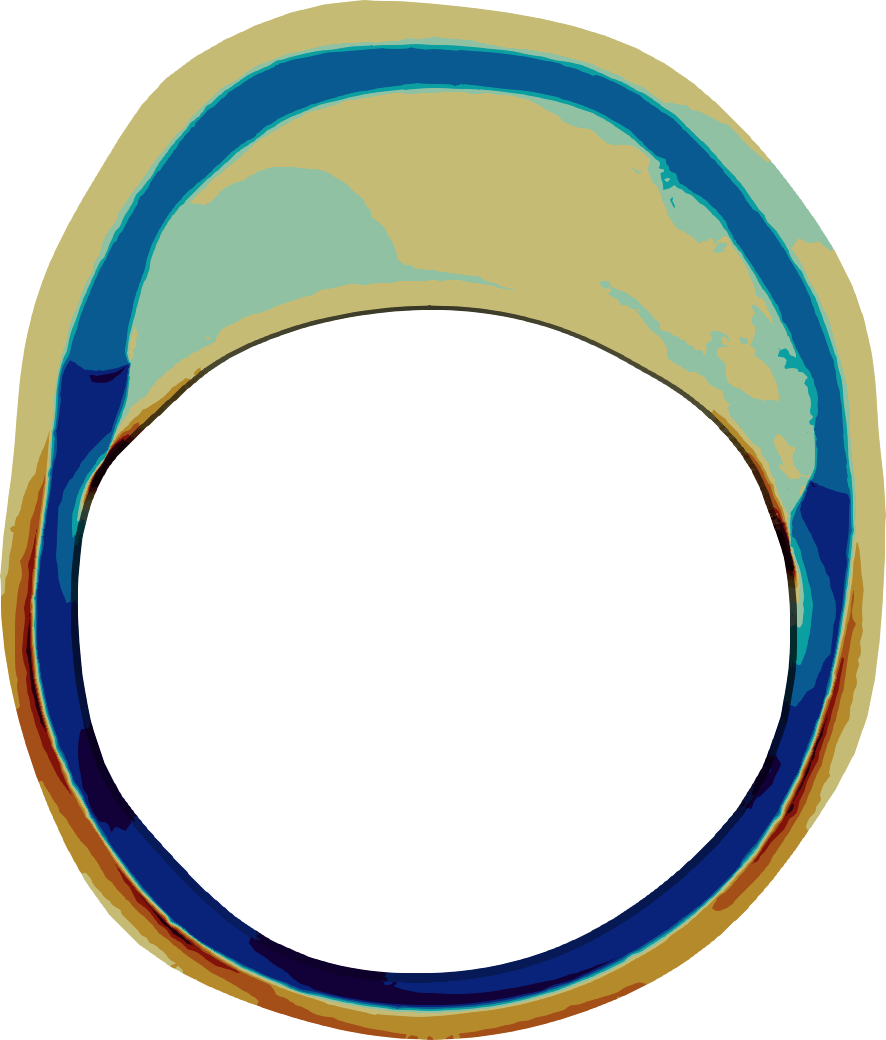}
        \par\vspace{-1.1\baselineskip}\noindent\llap{\footnotesize(d-1)\hspace{-0.6em}}
    \end{minipage}
    \hfill
    \begin{minipage}[b]{0.235\linewidth}
        \includegraphics[width=\linewidth,page=1]{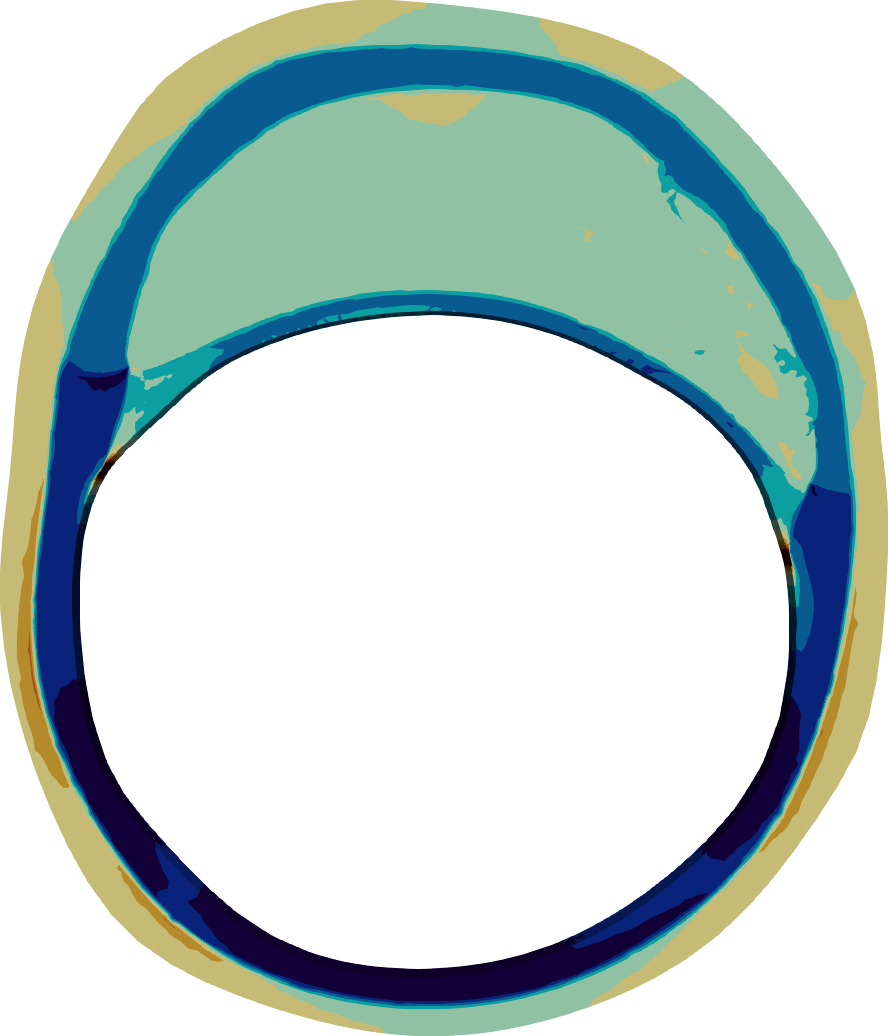}
        \par\vspace{-1.1\baselineskip}\noindent\llap{\footnotesize(d-2)\hspace{-0.6em}}
    \end{minipage}
    \hfill
    \begin{minipage}[b]{0.235\linewidth}
        \includegraphics[width=\linewidth,page=1]{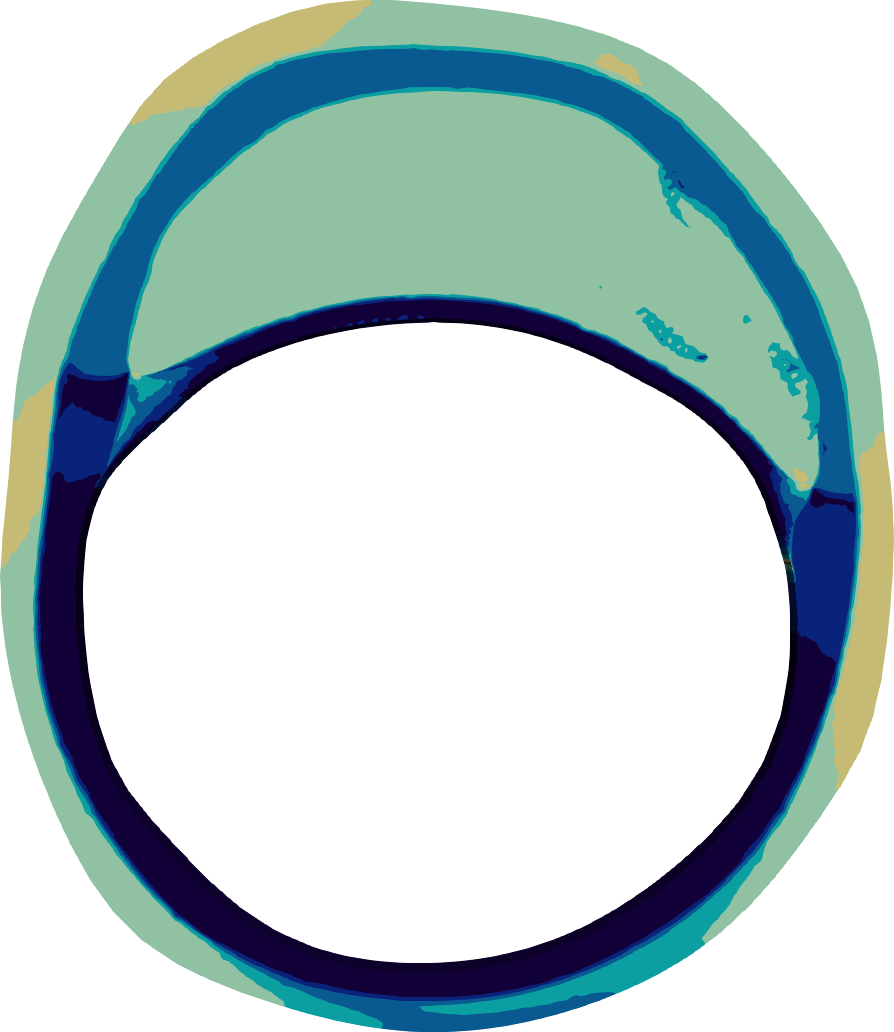}
        \par\vspace{-1.1\baselineskip}\noindent\llap{\footnotesize(d-3)\hspace{-0.6em}}
    \end{minipage}
    \hfill
    \begin{minipage}[b]{0.04\linewidth}
        \includegraphics[width=\linewidth,page=1]{legend50.png}
    \end{minipage}
    \par\medskip
    \begin{minipage}[b]{0.232\linewidth}
        \includegraphics[width=\linewidth,page=1]{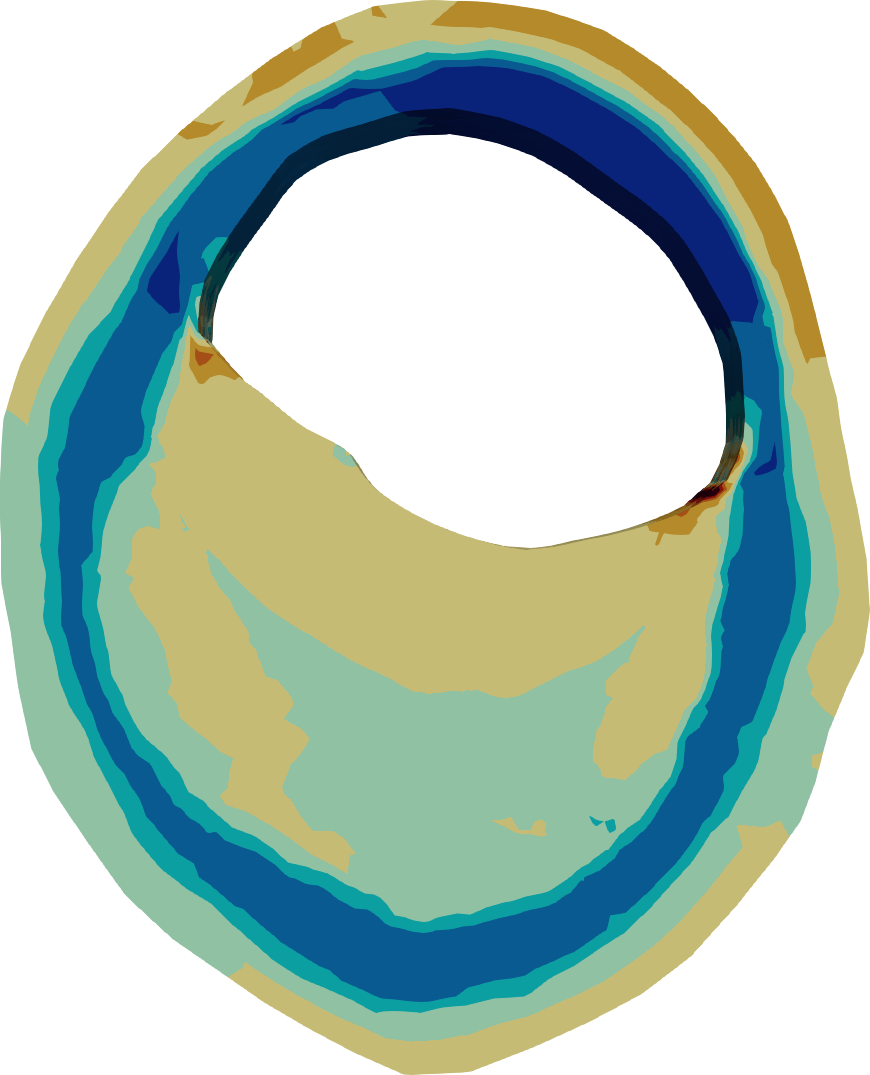}
        \par\vspace{-1.1\baselineskip}\noindent\llap{\footnotesize(e-1)\hspace{-0.6em}}
    \end{minipage}
    \hfill
    \begin{minipage}[b]{0.232\linewidth}
        \includegraphics[width=\linewidth,page=1]{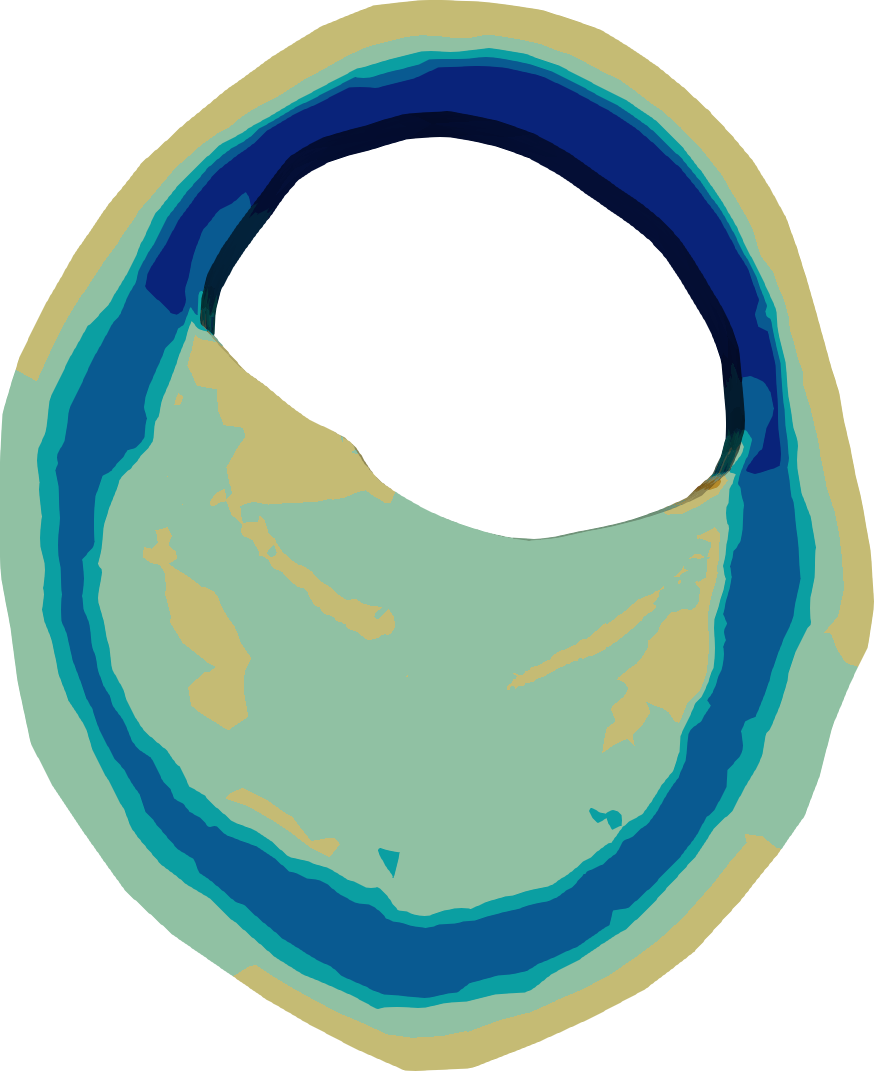}
        \par\vspace{-1.1\baselineskip}\noindent\llap{\footnotesize(e-2)\hspace{-0.6em}}
    \end{minipage}
    \hfill
    \begin{minipage}[b]{0.232\linewidth}
        \includegraphics[width=\linewidth,page=1]{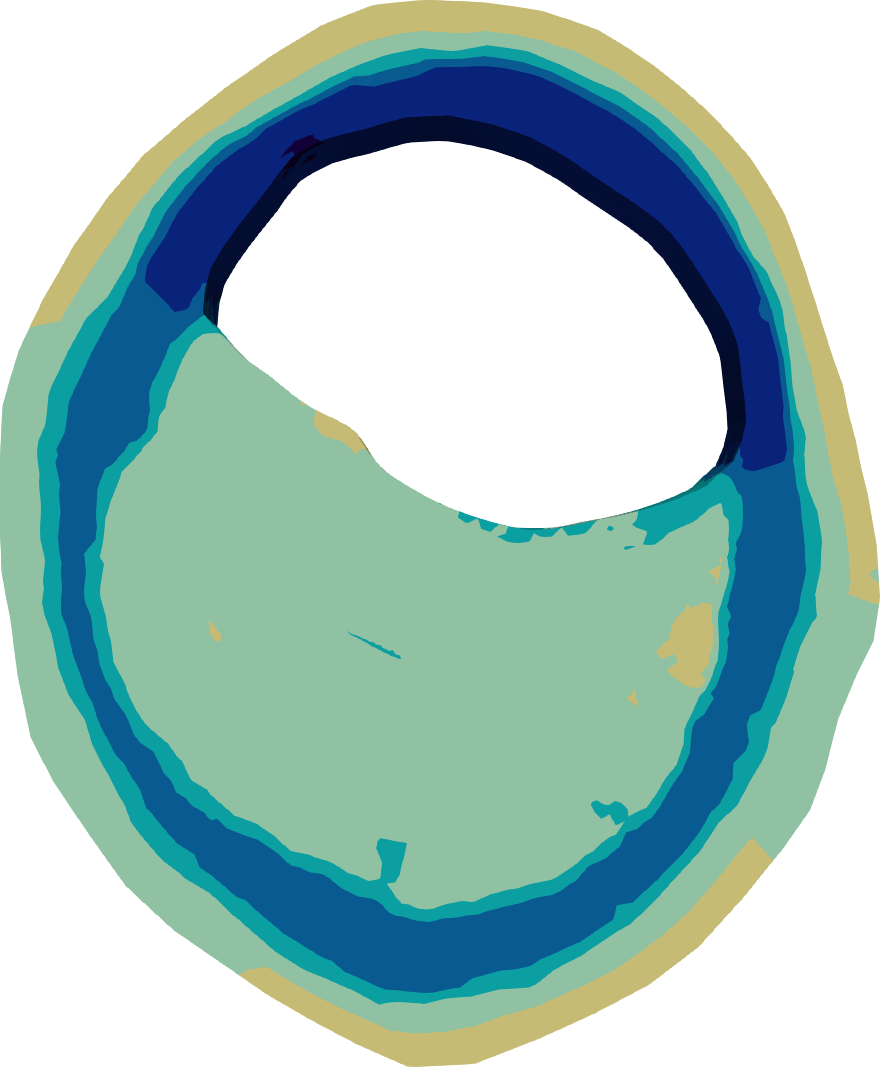}
        \par\vspace{-1.1\baselineskip}\noindent\llap{\footnotesize(e-3)\hspace{-0.6em}}
    \end{minipage}
    \hfill
    \begin{minipage}[b]{0.04\linewidth}
        \includegraphics[width=\linewidth,page=1]{legend50.png}
    \end{minipage}
    \end{minipage}%
    }
    \caption{%
        Blood pressure drop effects on change in von Mises stress in kPa post-drug treatment for different plaque morphologies. (a) to (e) represent the different artery geometries as depicted in Fig.~\ref{fig:arteries_geometry}, with (-1),(-2), and (-3) representing simulations with 10~mmHg, 20~mmHg, and 30~mmHg pressure drops respectively.
        \colsr{As expected, an increasing pressure drop reduces the stress increase in the fib. cap, with a median decrease at 30~mmHg drop in all geometries except the external iliac artery in (c)~(median $5.1$~kPa).
        Relevant stress increases for small and moderate blood pressure reduction remain mainly in the external iliac artery~(median $20.5$~kPa and $12.6$~kPa at 10 and 20~mmHg drop), whereas the median increase in the other geometries is at most $3.8$~kPa already at 10~mmHg drop}
    }
    \label{fig:plaque_bp_effects}
    \end{figure*}

    In the simulations presented above, we focused on direct effects of the drugs on the tissue response and did not account for the systemic impact.
    Therefore, we considered the blood pressure to remain constant at 85 mmHg post-drug treatment.
    However, antihypertensive drugs like amlodipine are primarily prescribed to reduce high blood pressure levels in patients and they do so by lowering the systemic vascular resistance.
    Since a systemic reduction in blood pressure would also affect the local hemodynamic loading conditions in the arteries, it is important to study how such a pressure reduction would influence the stress increases in the fibrous cap observed in the above simulations.
    To mimic the post-drug blood pressure reduction, we conducted simulations where the applied pressure was linearly reduced by fixed amounts between 1000--1200s in the loading protocol depicted in Fig.~\ref{fig:loading_curve}.
    To study different levels of success in
    blood pressure reduction, we
    studied pressure drops of 10~mmHg, 20~mmHg, and 30~mmHg for all artery geometries considered in Fig.~\ref{fig:arteries_geometry}.
    The post-drug change in von Mises stress for the different plaque morphologies and pressure drops is shown in Fig.~\ref{fig:plaque_bp_effects}.
    In all plaque morphologies, the post-drug fibrous cap stress increase reduces with increasing pressure reduction.
    \colsr{At a pressure drop of 10~mmHg, the fibrous cap of the external iliac artery still sees a widespread increase~(Fig.~\ref{fig:plaque_bp_effects}(c-1): median $20.5$~kPa, 5th--95th percentile: $2.2$ to $26.5$~kPa, compared to $28.9$~kPa without pressure drop, whereas the median increase in the other geometries drops to between $-0.1$~kPa~(Fig.~\ref{fig:plaque_bp_effects}(a-1)) and $3.8$~kPa~(Fig.~\ref{fig:plaque_bp_effects}(d-1)), with isolated regions experiencing higher increases~(99th percentile $11.8$~kPa in Fig.~\ref{fig:plaque_bp_effects}(e-1) to $33.2$~kPa in Fig.~\ref{fig:plaque_bp_effects}(d-1)).}
    However, with increasing pressure reduction, i.e., success of the drug intake, these regions reduce in size and magnitude.
    \colsr{At 20~mmHg reduction, the fibrous cap stress of both TCFA geometries decreases on median~($-4.1$~kPa in Fig.~\ref{fig:plaque_bp_effects}(a-2) and $-19.3$~kPa in Fig.~\ref{fig:plaque_bp_effects}(d-2)), with increases remaining only in isolated neck regions at the boundary of the fibrous cap and media~(99th percentile $3.4$~kPa and $-0.7$~kPa), and the median change in the tibial and the moderately stenosed artery is close to zero~($0.2$~kPa in Fig.~\ref{fig:plaque_bp_effects}(b-2) and $-0.6$~kPa in Fig.~\ref{fig:plaque_bp_effects}(e-2)), while the external iliac artery still shows a median increase of $12.6$~kPa~(Fig.~\ref{fig:plaque_bp_effects}(c-2)).
    With a 30~mmHg drop, the fibrous cap stress decreases on median in all geometries~($-2.4$~kPa in Fig.~\ref{fig:plaque_bp_effects}(b-3) to $-41.7$~kPa in Fig.~\ref{fig:plaque_bp_effects}(d-3)) with the exception of the external iliac artery~(Fig.~\ref{fig:plaque_bp_effects}(c-3)), which still shows a small increase~(median $5.1$~kPa, 5th--95th percentile: $-9.5$ to $7.3$~kPa) that reaches about $10$~kPa in some regions near the lumen~(99th percentile $11.4$~kPa).}

    \subsection{\colsr{Axial pre-stretch and blood pressure reduction}}\label{sec:stretch_bp}
\begin{figure*}[t]
    \centering
    \newcommand{\stbptitle}[1]{\makebox[\linewidth][c]{\hspace{0.19\linewidth}\fontfamily{lmr}\fontsize{8.43}{10}\selectfont #1}\par\vspace{-0.45\baselineskip}}
    \begin{minipage}[b]{0.32\linewidth}
        \stbptitle{TCFA~(Fig.~\ref*{fig:arteries_geometry}(a))}
        \includegraphics[width=\linewidth]{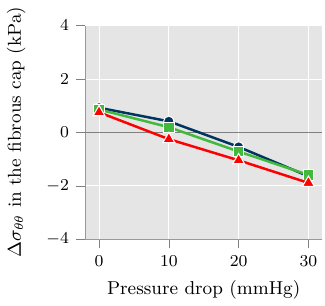}
        \par\vspace{-1.6\baselineskip}\noindent\llap{\footnotesize(a)\hspace{-0.6em}}
        \vspace{0.6\baselineskip}
    \end{minipage}%
    \hfill
    \begin{minipage}[b]{0.32\linewidth}
        \stbptitle{Tibial artery}
        \includegraphics[width=\linewidth]{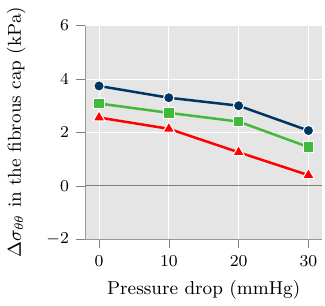}
        \par\vspace{-1.6\baselineskip}\noindent\llap{\footnotesize(b)\hspace{-0.6em}}
        \vspace{0.6\baselineskip}
    \end{minipage}%
    \hfill
    \begin{minipage}[b]{0.32\linewidth}
        \stbptitle{External iliac artery}
        \includegraphics[width=\linewidth]{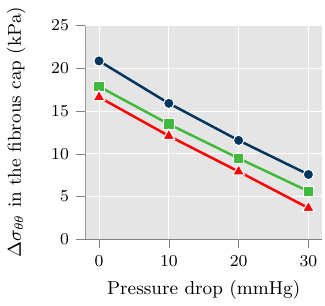}
        \par\vspace{-1.6\baselineskip}\noindent\llap{\footnotesize(c)\hspace{-0.6em}}
        \vspace{0.6\baselineskip}
    \end{minipage}
    \par
    \begin{minipage}[b]{0.32\linewidth}
        \stbptitle{TCFA~(Fig.~\ref*{fig:arteries_geometry}(d))}
        \includegraphics[width=\linewidth]{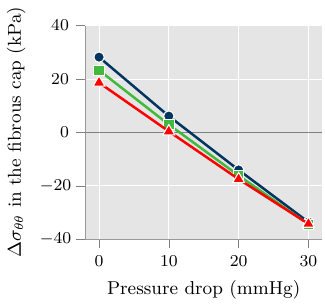}
        \par\vspace{-1.6\baselineskip}\noindent\llap{\footnotesize(d)\hspace{-0.6em}}
        \vspace{0.6\baselineskip}
    \end{minipage}%
    \hfill
    \begin{minipage}[b]{0.32\linewidth}
        \stbptitle{Moderately stenosed artery}
        \includegraphics[width=\linewidth]{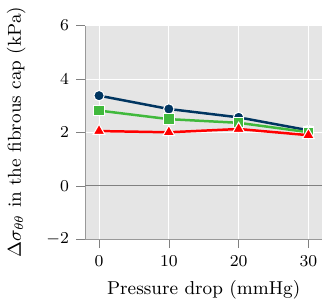}
        \par\vspace{-1.6\baselineskip}\noindent\llap{\footnotesize(e)\hspace{-0.6em}}
        \vspace{0.6\baselineskip}
    \end{minipage}%
    \hfill
    \begin{minipage}[b]{0.32\linewidth}
        \includegraphics[width=\linewidth]{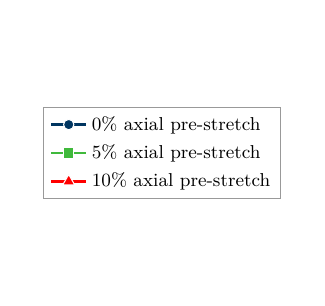}
        \vspace{0.6\baselineskip}
    \end{minipage}
    \caption{\colsr{Effect of axial pre-stretch and blood pressure reduction on the drug-induced change in circumferential stress in the fibrous cap~($\Delta\sigma_{\theta\theta} = \sigma_{\theta\theta,\text{aml.}} - \sigma_{\theta\theta,\text{control}}$ at $t=1500 \, \text{s}$) for 0\%, 5\%, and 10\% axial pre-stretch.
    (a) to (e) represent the artery geometries in the order of Fig.~\ref{fig:plaque_bp_effects}; each simulation with drug effects is compared with the simulation without drug effects at full pressure and the same pre-stretch.
    The lines indicate the volume-weighted medians; note the different stress ranges of the panels.
    Pre-stretch reduces the increase only moderately, whereas an increasing pressure drop reduces it in all geometries except the moderately stenosed artery with pre-stretch in (e); a clear decrease occurs only in the TCFA in (d) from 20~mmHg on}
    \label{fig:stretch_bp}}
    \end{figure*}

    \colsr{In Section~\ref{sec:axial_stretch}, an axial pre-stretch of the external iliac artery reduced the drug-induced increase in the von Mises stress of the fibrous cap.
    The decomposition of the stress tensor in Table~\ref{table:axial_stretch_components}, however, shows that this reduction stems mainly from the axial stress. With 10\% pre-stretch, the drug-induced increase in the axial stress of the fibrous cap drops from a median of $34.0$~kPa to $21.6$~kPa, whereas the increase in the circumferential stress drops only from $20.9$~kPa to $16.6$~kPa.
    The pre-stretch itself also raises mainly the axial stress, e.g., the median axial stress in the fibrous cap before treatment from $5.6$~kPa to $105.7$~kPa, while the median circumferential stress remains below $10$~kPa.
    With pre-stretch, the von Mises stress in the fibrous cap is therefore dominated by the axial stress imposed by the boundary condition rather than by the luminal pressure.
    The circumferential stress, in contrast, is governed by the luminal pressure and directly reflects the load that the relaxing medial layers shift onto the fibrous cap.
    To study the combined effect of axial pre-stretch and blood pressure reduction in all five geometries, we therefore evaluate the drug-induced change in the circumferential stress~$\sigma_{\theta\theta}$ of the fibrous cap.
    To this end, the simulations with and without drug effects were repeated for all geometries with 5\% and 10\% axial pre-stretch, and the simulations with drug effects additionally with the pressure drops of Section~\ref{sec:bp_reduction}.

    Fig.~\ref{fig:stretch_bp} shows the resulting changes.
    Without pressure drop and pre-stretch, the drug-induced increase in the circumferential stress of the fibrous cap is largest in the external iliac artery and the TCFA in Fig.~\ref{fig:stretch_bp}(d)~(median $20.9$~kPa and $28.2$~kPa), and small in the other three geometries~(median $0.9$ to $3.7$~kPa).
    In all geometries, the circumferential load released by the medial layers is taken up by the fibrous cap.
    An axial pre-stretch of 10\% reduces the increase only moderately, to a median of $16.6$~kPa and $18.6$~kPa in Fig.~\ref{fig:stretch_bp}(c) and (d), respectively, and to at most $2.6$~kPa in the other geometries.
    With increasing pressure drop, the circumferential stress increase in the fibrous cap decreases in all geometries, except for the moderately stenosed artery with pre-stretch, where it remains at about $2$~kPa.
    In the TCFA in Fig.~\ref{fig:stretch_bp}(d), the increase is almost offset at a pressure drop of 10~mmHg~(median $6.1$~kPa without and $0.4$~kPa with 10\% pre-stretch) and turns into a decrease at 20~mmHg and above~(median $-14.1$ to $-34.4$~kPa).
    In the external iliac artery, in contrast, an increase remains even at 30~mmHg~(median $7.6$~kPa without and $3.6$~kPa with 10\% pre-stretch), and in the tibial and the moderately stenosed artery, the small increase remains positive~(median $0.4$ to $2.1$~kPa at 30~mmHg).
    Hence, axial pre-stretch does not change the qualitative findings of the previous sections. 
    Amlodipine still shifts circumferential load from the medial layers onto the fibrous cap, and both pre-stretch and blood pressure reduction lower the resulting stress increase.}

\section{Discussion}
Calcium channel blockers like amlodipine are important and widely used antihypertensive drugs that reduce blood pressure through the systemic dilation of arteries.
Since they directly affect the active response of smooth muscle and thus the compliance of the arterial wall, they have the potential to significantly influence the mechanical stresses in the arterial wall.
This becomes particularly important in the presence of atherosclerosis, where spatially heterogeneous plaque formations exist in a subtle equilibrium with the surrounding tissue.
Due to the desired drug-induced reduction of blood pressure, stresses in the arteries are expected to decrease and thus reduce the risk of complications.
However, due to the complex morphology of atherosclerotic plaques which includes pathological components like calcifications, lipid pools, or fibrous caps, exhibiting a reduced drug-tissue interaction,
stresses may locally increase and thus lead to higher risk of adverse events such as plaque rupture.
To test this hypothesis, we developed a computational framework that couples drug transport, SMC contraction, anisotropic nonlinear arterial material behavior, and residual stresses to simulate the effects of calcium channel blockers on local stress distributions in atherosclerotic arteries.
\colsr{Simulations of the atherosclerotic external iliac artery showed that although most of the artery experiences a post-drug reduction in stress, the stress in the fibrous cap increases by a median of $28.9$~kPa.}
\colsr{This increase becomes smaller with decreasing SMC contractility before treatment.
It drops to a median of $22.4$~kPa for the highest NO activity considered~(45\% reduction of $\kappa$; see Fig.~\ref{fig:kappa_effects}), to $10.0$~kPa without contractility of the degenerated media~(see Fig.~\ref{fig:degen_kappa}), and to $12.0$~kPa for the lowest basal tone considered~($\eta$ at 10\% of its reference value; see Fig.~\ref{fig:basal_tone}).
In all three cases, the fibrous cap stress after drug treatment remains similar.
Hence, the drug-induced stress increase in the fibrous cap is mainly determined by the load that the contracted SMCs carry before treatment.}
\colsr{This is of particular relevance, since} atherosclerotic arteries exhibit endothelial dysfunction and reduced endothelial Nitric Oxide synthase (eNOS) activity~\cite{WanHe:2024:edm} causing reduced NO bioavailability in the medial layers.
\colsr{Our results indicate that with lower NO activity, the drug-induced increases in fibrous cap stresses are higher.}

\colsr{The magnitude of the increase further depends not only on how much amlodipine reaches the SMCs and how strongly they respond to it but also on the axial pre-stretch.}
\colsr{Lower diffusivities in the degenerated media and adventitia reduced the fibrous cap stress increase to a median of $20.6$~kPa for the lowest combination~(see Fig.~\ref{fig:transport_param_effects}(a)), whereas drug reaction rates above $10^{-2}\text{s}^{-1}$ reduced it strongly, to a median of $3.5$~kPa at $3\cdot 10^{-1}\,\text{s}^{-1}$~(see Fig.~\ref{fig:transport_param_effects}(b)), due to reduced drug availability.}
However, it is unlikely that such high reaction rates are physiologically relevant for amlodipine, given its slow clearance and long half-life~\cite{BurMoi:1994:upp,NayGua:1991:abt}.
\colsr{The stress increase also requires drug concentrations in the arterial wall to be well above the plasma level.
At the plasma level of about $0.06\,\mathrm{\mu M}$, the median increase in the fibrous cap was only $0.1$~kPa, whereas it was close to saturation at the assumed boundary concentration of $2\,\mathrm{\mu M}$~(see Fig.~\ref{fig:dose_effects}).
How strongly the SMCs respond to a given concentration had a smaller influence, since a less sensitive or flatter drug-tissue interaction reduced the increase by only about a quarter to a third~(see Fig.~\ref{fig:drug_response}).
Apart from the drug itself, the axial pre-stretch of the artery also lowered the increase, from a median of $28.9$~kPa without pre-stretch to $23.6$~kPa and $14.9$~kPa at a pre-stretch of 5\% and 10\%, respectively~(see Fig.~\ref{fig:axial_stretch}).
This reduction stems mainly from a smaller increase in the axial stress, whereas the increase in the circumferential stress of the fibrous cap, which reflects the load shifted from the medial layers, was reduced only moderately by the pre-stretch in all five geometries~(see Fig.~\ref{fig:stretch_bp}).
All other simulations were performed without axial pre-stretch, and the stress increases reported for them can therefore be regarded as an upper estimate in this respect.
In all these variations at the assumed boundary concentration, however, the drug-induced shift of load from the medial layers to the fibrous cap persisted.}

\colsr{Further simulations on four different atherosclerotic arterial morphologies showed that this increase in fibrous cap stress after drug treatment occurs in all of them, although its magnitude varies strongly.
The median increase was largest in one of the thin cap fibroatheromas~($28.4$~kPa, see Fig.~\ref{fig:plaque0}(c)), but only $3.9$ to $7.4$~kPa in the other three arteries, including the second thin cap fibroatheroma.
Hence, a thin fibrous cap alone does not determine the magnitude of the drug-induced stress increase.}

\colsr{Amlodipine has been shown to reduce systemic mean arterial pressure in hypertensive patients by an average of 10--20~mmHg~\cite{JefBhaRob:2015:ibp,GaoZhoYan:2021:eoa}.}
To ascertain the influence of systemic blood pressure reduction on the observed stress increases, we performed simulations with varying pressure drops up to 30~mmHg.
\colsr{Already at a pressure reduction of 10~mmHg, the median increase in the fibrous cap dropped to at most $3.8$~kPa in the four plaque morphologies, and at 20~mmHg and above, the median change was close to zero or negative in all geometries except the external iliac artery~(median $12.6$~kPa and $5.1$~kPa at 20~mmHg and 30~mmHg, respectively; see Fig.~\ref{fig:plaque_bp_effects}), indicating that sufficient blood pressure reduction offsets the drug-induced increase of local stresses.
The circumferential stress in the fibrous cap shows the same trends with and without axial pre-stretch~(see Fig.~\ref{fig:stretch_bp}).
At a pressure reduction of 30~mmHg, however, the circumferential stress still increases in the external iliac artery~(median $7.6$~kPa without pre-stretch) and slightly in the tibial and the moderately stenosed artery~(median of at most $2.1$~kPa).}

It must be noted that an observed increase in fibrous cap stress in and of itself does not necessarily imply an increased risk of plaque rupture.
This would depend on the quantitative values of increased stresses in relation to
the material strength of the fibrous cap which has high variability ranging from 72.1--259.3~kPa with a median of 158.3~kPa~\cite{TenFenZhaSutHuaBroJinGil:2015:aus}, which is often plaque specific.
\colsr{In our study, the fibrous cap of the thin cap fibroatheroma in Fig.~\ref{fig:plaque0}(c) already carried stresses within this range before treatment~(median $121.3$~kPa, 95th percentile $165.8$~kPa), which rose to a median of $150.1$~kPa~(95th percentile $210.4$~kPa) with amlodipine, close to the median strength.
At a pressure reduction of 10~mmHg, they remained close to the level before treatment~(median $125.8$~kPa, 95th percentile $175.0$~kPa).
In the second thin cap fibroatheroma, the fibrous cap stress remained well below the median strength~(95th percentile $96.4$~kPa after drug treatment).
A drug-induced stress increase in a fibrous cap that is already loaded close to its strength may indicate an increased rupture risk.}

\colsr{Clinical trials have not shown more coronary events with amlodipine, even where it lowered the blood pressure only slightly.
In CAMELOT~\cite{NisTuzLibThoGhaGarBerShiBueTopCAM:2004:eoa}, amlodipine lowered the blood pressure of normotensive patients with coronary artery disease by only $4.8/2.5$~mmHg, yet reduced cardiovascular events by 31\%, mainly through fewer coronary revascularizations and hospitalizations for angina, while nonfatal myocardial infarctions were not more frequent than with placebo~(14 vs.\ 19 patients).
Likewise, PREVENT~\cite{PitByiFurHunManMilRilFo:2000:eoa} showed no difference in myocardial infarctions~(19 vs.\ 20 patients) at a blood pressure reduction of about $8/4$~mmHg, and in ALLHAT~\cite{The:2002:moi}, where amlodipine and chlorthalidone lowered the blood pressure to almost the same level, the rate of fatal coronary heart disease and nonfatal myocardial infarction was the same.
In large trials of hypertensive patients with a considerable blood pressure reduction, amlodipine-based regimens did not result in more coronary events than the comparators either~\cite{JulKjeWebBruEkmHanHuaLarMcIMitPlaSchSmiZan:2004:oih,DahSevPouWedBeeCauColKjeKriMcIMehNieOBOst:2005:poc}, which agrees with our finding that a sufficient blood pressure reduction offsets the stress increase in the fibrous cap.
Adverse outcomes have only been reported for short-acting calcium channel blockers, i.e., immediate-release formulations that act within minutes and only for a few hours.
Their use was associated with an increased risk of myocardial infarction at comparable blood pressure levels in a case-control study~\cite{PsaHecKoeSisRagWeiRosLemSmiWahWagFur:1995:tro}, and in a meta-analysis of randomized trials in patients with coronary heart disease, immediate-release nifedipine increased the total mortality in a dose-dependent manner, up to a risk ratio of $2.83$ at a daily dose of 80~mg~\cite{FurPsaMey:1995:ndr}.
Neither ALLHAT for amlodipine, which is long-acting owing to its slow elimination, nor the ACTION trial for an extended-release formulation of nifedipine~\cite{PooLubKirVanWagDanJusFoxPocClaMotParBouDarHilHjaKraMolOttSeaSolWeb:2004:eol} showed such an increase.
These observations limit the clinical relevance of our predictions, but they do not rule them out.
The blood pressure reductions in CAMELOT and PREVENT correspond to a reduction in mean arterial pressure of only about 3 to 5~mmHg, i.e., below the smallest pressure drop simulated here, so our model would still predict a stress increase in the fibrous cap, which would be considerable in the external iliac artery and in the thin cap fibroatheroma in Fig.~\ref{fig:plaque0}(c), if amlodipine accumulated in the arterial wall to the assumed concentration.
At the plasma concentration, however, the predicted increase was negligible~(see Fig.~\ref{fig:dose_effects}), and the concentration that amlodipine reaches at the medial SMCs in human atherosclerotic arteries is not known.
Furthermore, neither trial was powered to detect differences in myocardial infarction, and their unselected patient populations are unlikely to contain many plaques whose fibrous caps are already loaded close to their strength, which is where our results indicate a potential risk.
The clinical endpoints also reflect the anti-ischemic and possibly blood pressure-independent effects of amlodipine~\cite{NisTuzLibThoGhaGarBerShiBueTopCAM:2004:eoa,BreIvaPolCheTuzHuFriNis:2006:ata}, which are not included in our model.
Finally, the absence of significant atheroma progression in the intravascular ultrasound substudy of CAMELOT concerns the growth of the plaque over two years, not the stress in the fibrous cap, and therefore does not test our prediction directly.}

\colsr{Whether the predicted load shift becomes clinically relevant thus depends on several factors that differ between patients and drugs.
Our simulations identified three of them, namely the drug concentration in the arterial wall, the blood pressure reduction and the axial pre-stretch.
The axial pre-stretch in particular differs between patients, since it decreases with age and in the presence of cardiovascular risk factors~\cite{KamSeaBowDeeDesBohPouMac:2016:isl,HorAdaKul:2017:aco}, and the stress increase may therefore be larger in older patients.
A further factor, which differs between drugs and was not addressed by our simulations, is the time course of the drug action.
Amlodipine relaxes the SMCs only gradually, since its plasma concentration peaks 6--8~h after oral administration~\cite{MerElli:1992:cpo}, and both its binding to the calcium channel~\cite{BurDodGar:1989:ppo} and its partitioning into cell membranes~\cite{MasCamWanHer:1989:col} are slow processes.
How fast a calcium channel blocker acts can, however, change the response of the circulation, as shown for nifedipine.
In healthy volunteers, a gradual rise of the plasma concentration of nifedipine lowered the diastolic blood pressure without reflex tachycardia, whereas a rapid rise to the same concentration caused tachycardia without lowering the diastolic blood pressure~\cite{KleVanDanFabUrqBre:1987:roi}, and in patients with coronary artery disease, immediate-release nifedipine activated the sympathetic nervous system earlier and more strongly than its extended-release formulation~\cite{ParDIoFloToa:2020:cos}.
A rapidly acting calcium channel blocker may thus relax the SMCs without a corresponding reduction in blood pressure, which is the scenario in which our model predicts the largest stress increase in the fibrous cap.
A gradual onset may further give the arterial wall time to adapt, since arteries remodel towards a preferred mechanical state~\cite{Hum:2008:vaa} and isolated arterioles remodel within days in response to sustained changes in SMC tone~\cite{BakVanVanEveSpaVan:2002:irf}.
This capacity may, however, be limited in the fibrous caps of rupture-prone plaques, which contain few SMCs~\cite{DavRicWooKatMan:1993:rot}.
The combination of these factors may explain why clinical trials of long-acting calcium channel blockers have not shown more coronary events, whereas the dose-dependent excess mortality with immediate-release nifedipine would be consistent with the mechanism described here, although other explanations have been proposed for it~\cite{FurPsaMey:1995:ndr,PsaHecKoeSisRagWeiRosLemSmiWahWagFur:1995:tro}.}

There are certain limitations to this study that need to be addressed in future work.
Firstly, we did not have access to patient-specific data to generate arterial geometries and fit material and SMC model parameters to patient- and tissue-specific data.
Our SMC model parameters were fit to data from rat arteries, which may not accurately represent human SMC behavior.
Further, we relied on histology-based reconstruction of geometries which is known to introduce artifacts due to tissue processing and 2D sectioning.
IVUS based virtual histology~\cite{WarYooMeyMolSamHay:2022:afe}, optical coherence tomography~\cite{ChaChaShiMacIftTeaKamBouKaa:2004:mao} or computed tomography angiography~\cite{CurRosMazNanConMat:2023:3dp} based reconstructions may provide more realistic geometries in future studies.
\colsr{Moreover, the material and SMC model parameters were assumed to be constant within each component and therefore change abruptly at the interfaces between the components, in particular between the fibrous cap and the healthy and degenerated media.
In reality, the tissue composition, and with it the material properties, are likely to change gradually across such transition zones.
Spatially varying volume fractions of the individual constituents, obtained, for instance, from a constrained mixture-based approach, could capture these gradients and would likely also reduce the stress concentrations at the sharp interfaces.
Therefore, to make meaningful comparisons, we use volume-weighted percentiles rather than peak values~(see Appendix~\ref{sec:app_mesh}).}
Additionally, we did not consider the effect of cyclic loading due to pulsatile blood flow.
Although we expect the qualitative trends to remain similar, the absolute stress values may change under pulsatile loading.
\colsr{Furthermore, the drug was applied at a constant boundary concentration and the blood pressure reduction was prescribed independently of it, hence the time course of the drug action, reflex responses of the circulation and an adaptation of the arterial wall after treatment were not captured.}
We considered a simplified drug transport model with an isotropic diffusion tensor and simple first-order reaction kinetics.
Although amlodipine transport is diffusion dominated, considering advection and more complex reaction kinetics, provided relevant data is available, may improve the accuracy of the model.
We considered NO effects to be homogeneous and time-independent throughout the artery due to lack of data.
A fluid-structure interaction based coupled model \cite{BalHeiKlaRheSch:2023:coa,BalDepFauForHeiKlaQuaRheSch:2015:nmo} extended to integrate blood flow based, shear stress induced endothelial NO production, and intramural NO transport may provide additional insights.

\colsr{The coupled model contains a large number of parameters for the passive and active material behavior, the drug transport, the drug-tissue interaction and the growth model, which were identified from experiments on different arteries and species.
A simpler model with fewer parameters, each of which could be measured on the simulated tissue, would be easier to calibrate and to adapt to individual patients, and would therefore be preferable for clinical application, provided that it still captures the mechanisms that determine the prediction.
The drug-induced load shift, however, results from the interaction of several mechanisms, some of which are complex by nature.
First, the active stress of the SMCs depends on their stretch~\cite{UhlBal:2023:cmo}, which differs between the healthy and the degenerated media and changes with the blood pressure.
Since the stress increase in the fibrous cap is mainly determined by the load that the contracted SMCs carry before treatment~(see Figs.~\ref{fig:kappa_effects}, \ref{fig:degen_kappa} and \ref{fig:basal_tone}), a constant active stress or a passive model with a drug-dependent stiffness would not capture how the shifted load depends on the SMC tone before treatment and on the blood pressure reduction.
Second, calcium channel blockers act on the calcium influx, and a model with an explicit intracellular calcium concentration allows the drug-tissue interaction to be identified directly at this level, e.g., from experimental data~\cite{WanJiaLiLai:2008:eos}~(see Fig.~\ref{fig:drug_interaction}), and to be combined with the stretch-dependent activation.
Third, the drug effect depends on the plaque morphology, since the drug diffuses through plaque components whose diffusivities differ by several orders of magnitude~\cite{BozChoBar:2014:mtt,HosHosBazCalHug:2012:mmo}, and the stress state before treatment depends on the residual stresses, which strongly affect the stress in the fibrous cap~\cite{OhaDubTraLeRioChaThiPetFin:2007:ior}.
To limit the consequences of the large number of parameters, the passive, active and drug-tissue interaction parameters were each identified on dedicated experiments, and the sensitivity of the results to the least certain parameters was examined in parameter studies~(see Sections~\ref{sec:NO_effects}, \ref{sec:transport_param_effects} and \ref{sec:axial_stretch}).
These showed that the magnitude of the stress increase depends strongly on the parameters, whereas the load shift itself persisted in all variations at the assumed drug concentration. The model predictions should therefore currently be regarded as indicative of the mechanism and its sensitivities rather than as quantitative predictions for individual patients.
At the same time, the parameter studies identify the quantities that dominate the result, i.e., the load carried by the SMCs before treatment, the drug concentration in the arterial wall and the blood pressure reduction. A reduced model that retains these and is calibrated against the full model and the experiments outlined in the following may thus be a suitable basis for patient-specific predictions.}

\colsr{A conclusive validation requires that all model parameters are identified on the same type of artery that is simulated, and that the model predictions are compared with measurements that were not used for the identification.
The geometry of each specimen would be acquired non-destructively, for example by optical coherence tomography at physiological pressure and axial stretch. After completion of all mechanical tests, the segment would be fixed under the same loading and sectioned, which allows the histologically identified components to be registered to the imaged geometry.
Polarized light microscopy
\citep[see][]{SchZeiPieRegHol:2012:dot}
could yield the layer-specific collagen fiber angles and their dispersion, and immunostaining of contractile markers could quantify the fraction of contractile smooth muscle cells in the healthy and degenerated media~\cite{HaoGabCamBacVirBoc:2006:pmo}, which is known to decrease in atherosclerotic plaques~\cite{BenSinOwe:2016:vsm}.
The passive material parameters could be identified from extension-inflation tests of the same segment in a Ca$^{2+}$-free solution at several axial stretches.
The parameters of the smooth muscle model used here were identified from pressure step experiments on rat middle cerebral arteries~\cite{JohElyTakWalWalCol:2009:css}.
These experiments would need to be repeated on the chosen artery.
In addition to the diameter response to stepwise increases in pressure, the phosphorylation of the regulatory myosin light chain could be measured in segments that are snap-frozen at defined pressures and times, as in \citet{JohElyTakWalWalCol:2009:css}, and the intracellular Ca$^{2+}$ concentration in the wall could be recorded simultaneously with the diameter, as in \citet{KnoNel:1998:roa}.
In addition, a depolarization-induced contraction with potassium chloride could be applied to induce maximal contraction, and the axial force could be measured at several pressures and axial stretches.
These data would allow the stretch-dependent calcium dynamics, the cross-bridge kinetics and the maximum active stress to be identified for the healthy and the degenerated media separately.

On the same segments, a cumulative concentration-response experiment with amlodipine applied to both the luminal and the adventitial surface could identify the parameters of the drug-calcium interaction in the target tissue.
Amlodipine fully reverses potassium-induced contraction of human internal mammary arteries~\cite{BaiLiuJinYanTanHe:2010:eoa}, but its strong partitioning into cell membranes leads to a slow onset of action~\cite{MasCamWanHer:1989:col}, which would require sufficiently long equilibration times at each concentration.
Measuring the diameter, the axial force and the myosin light chain phosphorylation at concentrations between the plasma level and the saturation level, followed by stepwise reductions of the luminal pressure, would reproduce the scenarios simulated in this study experimentally.

The transport parameters could be identified from the transmural drug distribution.
\citet{TzaVukKolAstEde:2010:lcd} incubated human and rabbit arteries in drug solutions, sectioned them en face into thin slices and determined the transmural drug content, from which they obtained layer-specific partition coefficients and uptake kinetics, and showed that the lesion composition strongly affects drug deposition, cf.~also~\citet{LevVukHwaEde:2004:sbt}.
Applying this protocol to amlodipine at several exposure times could provide the diffusivities and the reaction parameters of the individual plaque components as well as the partition coefficient between the surrounding fluid and the tissue, which would test the assumed boundary concentration.

The residual stresses generated by the growth model could be validated by opening angle experiments~\cite{ChuFun:1986:rsi} on rings adjacent to the tested segment.
Since a single opening angle cannot capture the layer-specific residual deformations~\cite{HolSomAueRegOgd:2007:l3r}, the full displacement field of the cut ring would be recorded and compared with the simulated unloaded and cut configuration.
Such experiments have been performed on human vulnerable coronary plaques, where the residual stresses were found to strongly affect the peak stress in the fibrous cap~\cite{OhaDubTraLeRioChaThiPetFin:2007:ior}.

Finally, the predicted changes in cap strain could be assessed in vivo.
Intravascular ultrasound elastography and palpography could identify high-strain regions in coronary plaques~\cite{SchDeMasStrPasBoeSerVan:2003:cvp,SchRegMasMcFSaiDisDeDeVanSer:2004:ioh}, and compound ultrasound strain imaging enables a noninvasive assessment of carotid plaques~\cite{HanDeBotMolPasDe:2016:von}.
Serial strain measurements before and after initiation of amlodipine therapy, combined with patient-specific simulations, would test the predicted change in cap strain, as serial intravascular imaging under amlodipine has been shown to be feasible~\cite{NisTuzLibThoGhaGarBerShiBueTopCAM:2004:eoa}.
Agreement between the predicted and measured quantities would validate the coupled model, while discrepancies would indicate which of the submodels, i.e., drug transport, smooth muscle activation or residual stress, requires further development.}

\colsr{
\section{Conclusion}
In summary, our results suggest that calcium channel blockers do not only reduce the stresses in atherosclerotic arteries, but can also shift load from the relaxed medial SMCs to the fibrous cap.
The hypothesis that this shift may increase the rupture risk of vulnerable plaques whose fibrous caps are already loaded close to their strength, in particular when the blood pressure reduction is insufficient, is neither contradicted nor supported by the available clinical trials.
Testing it comprehensively requires the experimental validation outlined above, patient-specific simulations, and clinical studies that relate serial intravascular imaging of the fibrous cap to the blood pressure response in patients with vulnerable plaques.
If confirmed, simulations as shown in this paper could contribute to developing improved treatment protocols, for example by emphasizing a sufficient blood pressure reduction and a gradual onset of drug action in patients with vulnerable plaques.}

\section*{Acknowledgement}
Financial funding from the Deutsche Forschungsgemeinschaft (DFG) through the Priority Program 2311 “Robust coupling of continuum-biomechanical in silico models to establish active biological system models for later use in clinical applications—Co-design of modeling, numerics and usability,” project ID 465228106, is greatly appreciated.
Computations for this publication were performed on the HPC cluster Elysium of Ruhr University Bochum, subsidized by the DFG (INST 213/1055-1).
The authors also gratefully acknowledge the scientific support and HPC resources provided by the Erlangen National High Performance Computing Center (NHR@FAU) of the Friedrich-Alexander-Universität Erlangen-Nürnberg (FAU) under the NHR project k105be.
NHR funding is provided by federal and Bavarian state authorities.
NHR@FAU hardware is partially funded by the German Research Foundation (DFG)—440719683.
Open Access funding enabled and organized by Project DEAL.

\section*{Declarations}
\subsection*{Competing interests}
The authors have no relevant financial or non-financial interests to disclose.
\subsection*{Data availability}
Simulations were performed using the open-source Finite Element and Domain Decomposition Library (FEDDLib) \cite{feddlib}.
The version of FEDDLib along with the required Interface2 library used for this study is archived at \url{https://doi.org/10.60517/173ad37d-1490-467c-94cb-7769483bb732}.
It contains all meshes, input and parameter files required to reproduce the results of this study.
\subsection*{Author Contributions}
Sharan N. Ramesh: Formal analysis, Investigation, Methodology, Software, Validation, Visualization, Writing\textemdash original draft, Writing\textemdash review and editing;
Lea Sa{\ss}mannshausen: Software, Writing\textemdash review and editing;
Axel Klawonn: Conceptualization, Funding acquisition, Supervision, Writing\textemdash review and editing;
Daniel Balzani: Conceptualization, Funding acquisition, Methodology, Supervision, Writing\textemdash review and editing.
\subsection*{Funding}
This work was supported by the Deutsche Forschungsgemeinschaft (DFG) through Priority Program 2311, project ID 465228106.

\begin{appendices}
    \section{Mesh convergence}\label{sec:app_mesh}
    \begin{table*}[htb]
        \caption{\colsr{Mesh convergence: drug-induced change in von Mises stress~(kPa, $\sigma_\text{aml.} - \sigma_\text{control}$ at $t=1500 \, \text{s}$) in the two thin cap fibroatheromas of Fig.~\ref{fig:arteries_geometry}(a) and (d), for the original and the smoothed geometry~(cap shoulders rounded with a fillet radius of $0.1$~mm) at four refinement levels with local element size $h$ and number of degrees of freedom~(DOFs; displacement and concentration unknowns of the quadratic tetrahedral discretization). Volume-weighted percentiles in the fibrous cap, largest nodal value in the fibrous cap, and volume-weighted medians in the healthy and degenerated media. The main text uses the original geometries at L2}\label{tab:mesh_convergence}}
        \begin{center}\footnotesize  %
            \setlength{\tabcolsep}{4pt}  %
            \begin{tabular}{llccccccccc} \toprule
                & & & & \multicolumn{5}{c}{Fibrous cap} & Healthy media & Degen. media \\ \cmidrule(lr){5-9}
                Geometry & Level & $h$~(mm) & DOFs & 5th & 50th & 95th & 99th & Max. & 50th & 50th \\ \midrule
                \multicolumn{11}{l}{Thin cap fibroatheroma, Fig.~\ref{fig:arteries_geometry}(a)} \\
                Original & L0 & $0.12$  & $249468$  & $-5.9$ & $3.9$ & $16.1$ & $22.1$ & $212$ & $-29.6$ & $-20.5$ \\
                         & L1 & $0.06$  & $319468$  & $-5.5$ & $3.8$ & $17.2$ & $24.6$ & $154$ & $-30.3$ & $-20.7$ \\
                         & L2 & $0.03$  & $583496$  & $-5.6$ & $3.9$ & $18.0$ & $28.0$ & $396$ & $-31.0$ & $-21.0$ \\
                         & L3 & $0.015$ & $1846420$ & $-5.2$ & $3.8$ & $17.7$ & $29.4$ & $577$ & $-30.8$ & $-20.7$ \\
                Smoothed & L0 & $0.12$  & $251588$  & $-5.9$ & $3.8$ & $16.0$ & $22.0$ & $104$ & $-29.7$ & $-20.5$ \\
                         & L1 & $0.06$  & $324716$  & $-5.7$ & $3.8$ & $17.4$ & $24.8$ & $64$  & $-30.5$ & $-20.6$ \\
                         & L2 & $0.03$  & $588052$  & $-5.5$ & $3.7$ & $17.6$ & $27.9$ & $76$  & $-30.8$ & $-20.6$ \\
                         & L3 & $0.015$ & $1849844$ & $-5.4$ & $3.8$ & $17.9$ & $30.0$ & $157$ & $-30.9$ & $-20.8$ \\ \midrule
                \multicolumn{11}{l}{Thin cap fibroatheroma, Fig.~\ref{fig:arteries_geometry}(d)} \\
                Original & L0 & $0.12$  & $242044$  & $2.9$ & $27.7$ & $41.6$ & $60.7$ & $492$  & $-29.9$ & $-23.1$ \\
                         & L1 & $0.06$  & $274764$  & $2.6$ & $27.9$ & $42.5$ & $61.0$ & $392$  & $-29.9$ & $-23.2$ \\
                         & L2 & $0.03$  & $376864$  & $2.2$ & $27.9$ & $43.4$ & $68.1$ & $1064$ & $-30.0$ & $-23.2$ \\
                         & L3 & $0.015$ & $882076$ & $2.1$ & $27.9$ & $44.8$ & $67.4$ & $1497$ & $-30.1$ & $-23.3$ \\
                Smoothed & L0 & $0.12$  & $243372$  & $2.5$ & $27.8$ & $42.0$ & $59.9$ & $169$  & $-29.9$ & $-23.2$ \\
                         & L1 & $0.06$  & $277080$  & $2.4$ & $27.8$ & $42.8$ & $61.5$ & $161$  & $-30.1$ & $-23.2$ \\
                         & L2 & $0.03$  & $378096$  & $2.1$ & $27.9$ & $43.9$ & $66.6$ & $354$  & $-30.1$ & $-23.2$ \\
                         & L3 & $0.015$ & $885824$ & $2.1$ & $27.9$ & $45.5$ & $69.7$ & $522$  & $-30.2$ & $-23.2$ \\ \bottomrule
            \end{tabular}
        \end{center}
    \end{table*}
    \colsr{At its shoulders, the fibrous cap tapers to a sharp edge, where the boundary between the cap and the media meets the lumen surface at a small angle~($15$--$23^\circ$ in the two thin cap fibroatheromas).
    These sharp edges result from the manual segmentation of the histological cross-sections.
    At such sharp edges, the stress in an elastic solid is singular, and its peak value therefore does not converge with mesh refinement.
    We therefore examined the convergence of the volume-weighted percentiles used throughout this paper for the two thin cap fibroatheromas in Fig.~\ref{fig:arteries_geometry}(a) and (d), in which the fibrous cap is thinnest.
    Each geometry was meshed with Gmsh at four refinement levels L0--L3, with local element sizes of $0.12$, $0.06$, $0.03$ and $0.015$~mm at the cap shoulders and other stress concentrations and a global element size of $0.12$~mm and $0.15$~mm, respectively.
    In addition, a smoothed variant of each geometry was created, in which the cap shoulders were rounded with a fillet radius of $0.1$~mm, which increases the angle at the edge to about $44^\circ$ while the minimum cap thickness remains unchanged.
    The simulations in the main text use the original geometries at level L2.
    The two thin cap fibroatheromas were chosen since their thin fibrous caps and sharp cap shoulders make them the most demanding cases for the discretization.
    The meshes of the other three arteries have a similar or smaller element size in the fibrous cap~(median edge length $0.041$--$0.054$~mm, compared with $0.047$--$0.048$~mm at L2).
    The meshes used in the main text have $583\,496$, $375\,988$, $177\,200$, $376\,864$ and $144\,496$ degrees of freedom for the arteries in Fig.~\ref{fig:arteries_geometry}(a) to (e), respectively.

    Table~\ref{tab:mesh_convergence} lists the drug-induced change in von Mises stress for all meshes.
    Between L2 and L3, the medians and the 5th and 95th percentiles in the fibrous cap as well as the medians in the healthy and degenerated media change by at most $1.6$~kPa, and they differ by less than $1$~kPa between the original and the smoothed geometries at L3.
    The 99th percentile in the fibrous cap, which is more strongly influenced by the regions near the sharp edges, changes by up to $3.1$~kPa~(8\%) between L2 and L3.
    In contrast, the largest nodal value in the fibrous cap increases with refinement from $212$ to $577$~kPa and from $492$ to $1497$~kPa between L0 and L3, and it remains considerably lower in the smoothed geometries, although it rises there as well~($157$ and $522$~kPa at L3).
    These peak values thus reflect the singularity at the sharp edges rather than a converged stress concentration.
    In the thin cap fibroatheroma in Fig.~\ref{fig:arteries_geometry}(d), the elements within $0.3$~mm of the two sharp edges hold about 2.5\% of the fibrous cap volume but almost all of the volume above the 99th percentile; leaving them out does not change the median~($27.9$ and $27.8$~kPa at L2).
    The volume-weighted medians and 5th--95th percentile ranges reported in this paper are therefore practically independent of the mesh and of the local shape of the cap shoulders.}

    \section{Contractility, drug response, and axial stretch}\label{sec:analysis}
    \setcounter{table}{0}
    \subsection{\colsr{Effect of smooth muscle contractility}} \label{sec:NO_effects}
    \colsr{The drug-induced shift of load from the medial layers to the fibrous cap depends on how much SMC contraction can be attenuated by amlodipine.
    We therefore investigated three factors that determine the contractile state of the SMCs before drug treatment, namely NO-mediated vasodilation, the contractility of the degenerated media, and the basal tone.}

\begin{figure}[t]
        \centering
        \includegraphics[width=0.45\textwidth]{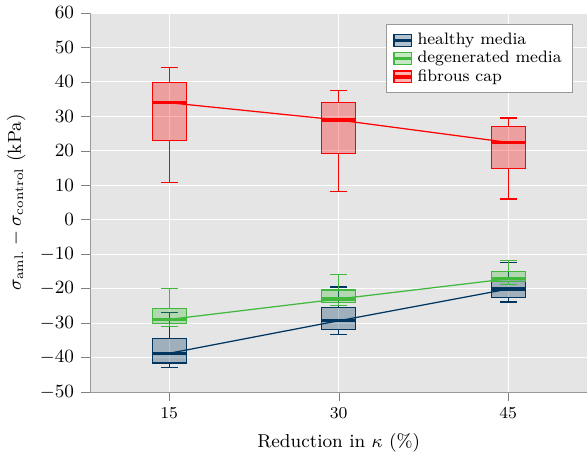}
        \caption{\colsr{Effect of NO-mediated vasodilation on drug-induced changes in von Mises stress in the external iliac artery, for 15\%, 30\%, and 45\% reduction of the maximum active force $\kappa$ in the healthy and degenerated media~(Table~\ref{table:kappa_values}).
        The change denotes the difference in von Mises stress between the simulations with and without drug effects at $t=1500 \, \text{s}$.
        Each box spans the 25th to 75th percentile of the component's volume, with the median indicated by a thick horizontal line, whereas the whiskers represent the 5th and 95th percentiles, so $90\%$ of the component's volume lies between them.
        With increasing $\kappa$ reduction, i.e., increasing NO activity, the drug induced stress reduction in the medial layers becomes smaller~(healthy media: median $-39$ to $-20$~kPa), and so does the stress increase in the fibrous cap~(median $34$ to $22$~kPa)} \label{fig:kappa_effects}}
    \end{figure}

    \colsr{We first consider NO, which} is a potent vasodilator and reduces contraction levels in arteries. %
    In atherosclerotic arteries, endothelial dysfunction often results in reduced NO bioavailability in the medial layers~\cite{WanHe:2024:edm}.
    The magnitude of this reduction could vary significantly in different plaques and even locally within a plaque, based on the extent of endothelial damage and other plaque specific factors like geometry and composition.
    \begin{table}[htb]
        \caption{$\kappa$~(maximum achievable active force in kPa) values for different reduction levels for the healthy and degenerated media components.
        Degenerated media is assumed to have 33\% less maximum contractile force than healthy media.
        Decreasing $\kappa$ values signify increasing NO influence\label{table:kappa_values}}
        \centering
        \begin{tabular}{rccc} \toprule
            Component & 15\% red. & 30\% red. & 45\% red. \\ \midrule
            Healthy media & $126.023$  & $103.783$ & $81.5441$ \\
            Degen. media & $84.4352$ & $69.5349$ & $54.6345$ \\ \bottomrule
        \end{tabular}
    \end{table}
    In our model, we account for the influence of NO by reducing the value of $\kappa$, the parameter that controls the maximum achievable active force in SMCs.
    In general, the reduction in $\kappa$ may be different at every material point and time instant~(based on local NO levels).
    Due to a lack of experimental data, we assume a constant value of $\kappa$ reduction in space and time, which is tantamount to assuming a homogeneous, time-independent NO effect throughout the artery.
    To investigate the influence of varying NO levels on
    \colsr{stress distributions in drug treated arteries}, we conducted simulations with three levels of NO activity.
    This was effected by choosing three values of $\kappa$ reduction: 15\%, 30\%, and 45\%~(see Table~\ref{table:kappa_values}).
    In all three simulations, the stress distributions remained similar to Fig.~\ref{fig:iliac_artery}(b) which depicts the 30\% case.
    \colsr{Fig.~\ref{fig:kappa_effects} shows the distribution of the stress difference~(drug~$-$~control) in each component for the three reduction levels.}
    With increasing reduction of $\kappa$~(increasing NO effect), the magnitude of stress lowered in the medial layers reduces.
    For an artery with already reduced SMC contractility due to NO, the effect of amlodipine is less pronounced, since both NO and amlodipine essentially take on the same function.
    Accordingly, the post-drug stress increase in the fibrous cap \colsr{also} reduces with increasing NO activity.

\begin{figure*}[t]
    \centering
    \begin{minipage}[b]{0.32\linewidth}
        \vspace{0pt}
        \includegraphics[width=\linewidth]{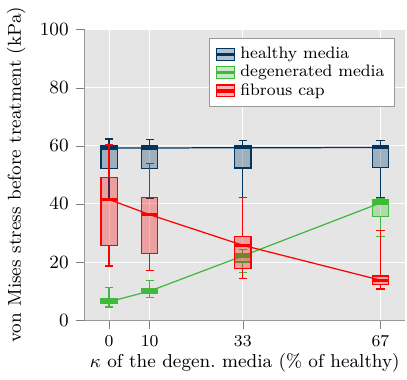}
        \par\vspace{-1.6\baselineskip}\noindent\llap{\footnotesize(a)\hspace{-0.6em}}
        \vspace{0.8\baselineskip}
    \end{minipage}%
    \hfill
    \begin{minipage}[b]{0.32\linewidth}
        \vspace{0pt}
        \includegraphics[width=\linewidth]{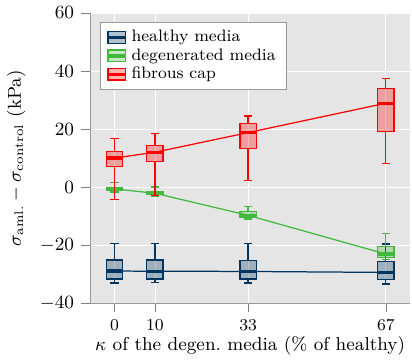}
        \par\vspace{-1.6\baselineskip}\noindent\llap{\footnotesize(b)\hspace{-0.6em}}
        \vspace{0.8\baselineskip}
    \end{minipage}%
    \hfill
    \begin{minipage}[b]{0.32\linewidth}
        \vspace{0pt}
        \includegraphics[width=\linewidth]{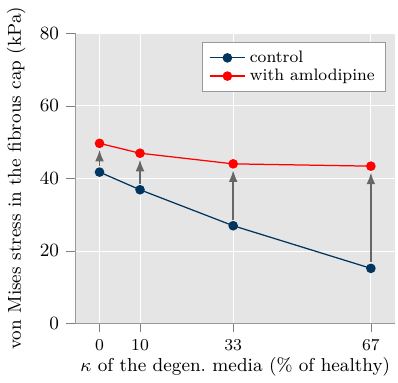}
        \par\vspace{-1.6\baselineskip}\noindent\llap{\footnotesize(c)\hspace{-0.6em}}
        \vspace{0.8\baselineskip}
    \end{minipage}
        \caption{\colsr{Effect of the contractility of the degenerated media on the von Mises stress in the external iliac artery, for a maximum active force~$\kappa$ of the degenerated media of 67\%~(base case), 33\%, 10\%, and 0\% of the healthy media value.
        (a) shows the von Mises stress before drug treatment~($t=860 \, \text{s}$), (b) depicts the difference in von Mises stress between the simulations with and without drug effects~($\sigma_\text{aml.} - \sigma_\text{control}$) at $t=1500 \, \text{s}$, with boxes depicting the $25$--$75$th percentiles and whiskers the $5$--$95$th percentiles.
        (c) shows the median von Mises stress in the fibrous cap at $t=1500 \, \text{s}$ without~(control) and with amlodipine, the arrows indicating the drug-induced change.
        With decreasing $\kappa$ in the degenerated media, i.e., with the degenerated media behaving less like the healthy media, the fibrous cap stress before drug intake increases significantly, but the gain in fibrous cap stress induced by drug-intake decreases.
        However, the total stress after drug-intake shows larger values for higher degenerations which might actually present more problematic scenarios.
        } \label{fig:degen_kappa}}
    \end{figure*}

    \colsr{As the second factor, we consider the contractility of the degenerated media, which we have hitherto assumed to be 33\% lower than that of the healthy media.
    \citet{Sch:1994:pva} reported a reduced density of SMCs in the atherosclerotic media.
    Further, \citet{MarLomManBacFerChiAlfFog:2021:cbf} found in their semi-quantitative confocal microscopy analysis of marker expressions that the smooth-muscle myosin heavy chain and $\alpha$-smooth muscle actin distribution in the altered media was rare, heterogeneous, and scattered.
    This supports the notion that SMCs are sparsely distributed and with reduced contractile machinery in the degenerated media.
    However, this reduced contractility of the degenerated media can vary significantly between plaques.
    To investigate the mechanical effect that this variation has on the drug-induced changes in the component-wise stress distributions, we conducted simulations with varying amounts of maximum achievable active force~($\kappa$) in the degenerated media between 0\% and 67\% of the healthy media value, while keeping all other parameters unchanged.
    It must be noted, though, that proliferative or synthetic phenotype SMCs and their mechanical response may not be captured adequately by the current model, since they may have different cellular dynamics such as altered cross-bridge cycling, membrane depolarization kinetics, mechanosensitivity, and contractile machinery when compared to regular contractile phenotype SMCs.
    Nevertheless, this assumption still provides an approximation of the purely mechanical aspect of altered SMC behavior in the degenerated media.
    From Fig.~\ref{fig:degen_kappa}(a), which shows the volume-weighted median, 25th--75th and 5th--95th percentiles of the von Mises stresses before treatment for the four simulated values of $\kappa$, it can be seen that with lower $\kappa$ values, the degenerated media carries reduced load and this in turn gets shifted onto the fibrous cap.
    The drug-induced increase in fibrous cap von Mises stress changes from a median of $28.9$~kPa to $10.0$~kPa when $\kappa$ reduces from 67\% to 0\% of the healthy media value~(Fig.~\ref{fig:degen_kappa}(b)).
    In contrast, the von Mises stress in the fibrous cap after drug treatment remains relatively unchanged, with a median of $43.4$~kPa for the 67\% case and $49.6$~kPa for the 0\% case~(Fig.~\ref{fig:degen_kappa}(c)).
    Even when there is no contraction in the degenerated media~(0\% case), there is a median increase of $10.0$~kPa, with a $5$th--$95$th percentile range of $-4.2$~kPa to $16.8$~kPa.
    This can be attributed to the unloading of the healthy media, which remains unchanged at a median of $-29$~kPa.
    }

\begin{figure}[t]
        \centering
        \includegraphics[width=0.45\textwidth]{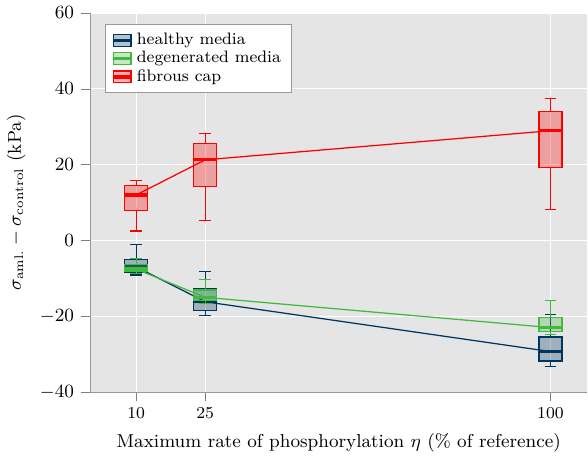}
        \caption{\colsr{Effect of the basal tone on drug-induced changes in von Mises stress in the external iliac artery, for a maximum rate of phosphorylation $\eta$ of 100\%~(base case), 25\%, and 10\% of its reference value.
        The change denotes the difference in von Mises stress between the simulations with and without drug effects~($\sigma_\text{aml.} - \sigma_\text{control}$) at $t=1500 \, \text{s}$, with boxes depicting the 25th--75th percentiles and whiskers the 5th--95th percentiles.
        With decreasing basal tone, the drug-induced stress reduction in the healthy media becomes smaller~(median $-29$ to $-7$~kPa), and so does the stress increase in the fibrous cap~(median $29$ to $12$~kPa)} \label{fig:basal_tone}}
    \end{figure}

    \colsr{Finally, the basal tone, i.e., the contraction of the SMCs before drug intake, is governed in our model by the maximum rate of phosphorylation $\eta$ in Eq.~\eqref{eq:Ca}.
    Since the basal tone of atherosclerotic arteries in vivo is uncertain and may be lower than assumed here, we repeated the simulations with and without drug effects with $\eta$ reduced to 25\% and 10\% of its reference value, while keeping all other parameters unchanged.
    The lower basal tone reduces the median von Mises stress in the healthy media before treatment from $59.4$~kPa to $45.8$~kPa and $37.2$~kPa, respectively, while the fibrous cap carries more load~(median $13.7$~kPa to $22.0$~kPa and $29.3$~kPa).
    Accordingly, the drug-induced stress reduction in the healthy media drops from a median of $-29.3$~kPa to $-16.2$~kPa and $-7.0$~kPa, and the stress increase in the fibrous cap from $28.9$~kPa to $21.3$~kPa~(5th--95th percentile: $5.2$ to $28.1$~kPa) and $12.0$~kPa~($2.5$ to $15.9$~kPa)~(see Fig.~\ref{fig:basal_tone}).}

    \colsr{In all three studies, the drug-induced stress increase in the fibrous cap becomes smaller with decreasing SMC contractility, whereas the fibrous cap stress after drug treatment remains similar, with medians of $42.8$--$45.1$~kPa for the three NO levels, $43.4$--$49.6$~kPa for the degenerated media contractility, and $43.4$--$47.1$~kPa for the basal tone.
    Before treatment, in contrast, the fibrous cap stress differs considerably between the cases~(medians of $10.5$--$20.3$~kPa, $13.7$--$41.5$~kPa, and $13.7$--$29.3$~kPa, respectively).
    The drug-induced increase in fibrous cap stress is therefore mainly determined by the load that the contracted SMCs carry before treatment.
    The lower the contractility, the more load the fibrous cap already carries and the smaller the increase.}
    \subsection{\colsr{Effect of drug availability and response}}\label{sec:transport_param_effects}
    \colsr{The drug-induced stress changes further depend on how much amlodipine reaches the SMCs and how strongly they respond to it.
    We therefore varied the drug transport parameters, the drug concentration at the arterial boundaries, and the drug-tissue interaction parameters.}

\begin{figure*}[t]
    \centering
    \begin{minipage}[b]{0.40\linewidth}
        \vspace{0pt}
        \includegraphics[width=\linewidth,page=1]{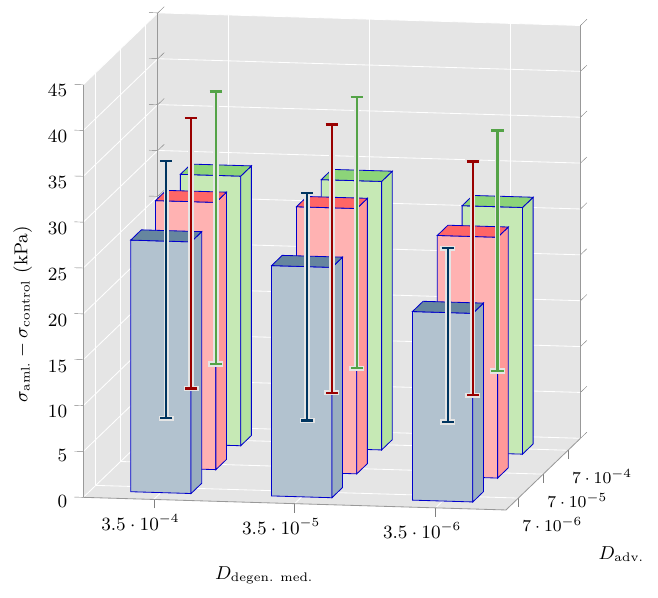}
        \par\vspace{-0.75\baselineskip}\noindent\llap{\footnotesize(a)\hspace{-0.6em}}
        \vspace{0.8\baselineskip}
    \end{minipage}
    \hspace{0.05\linewidth}
    \begin{minipage}[b]{0.47\linewidth}
        \vspace{0pt}
        \includegraphics[width=\linewidth,page=1]{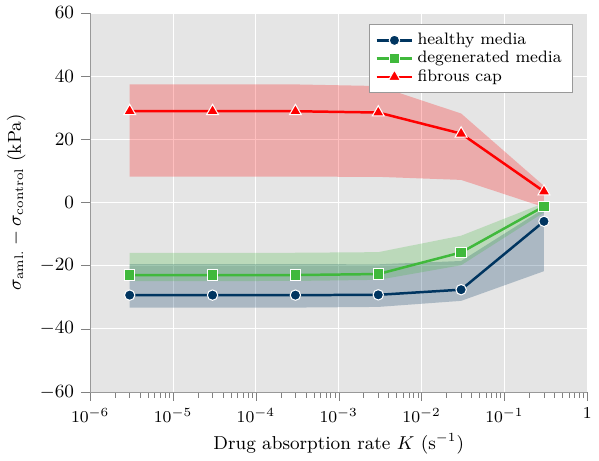}
        \par\vspace{-1.6\baselineskip}\noindent\llap{\footnotesize(b)\hspace{-0.6em}}
        \vspace{0.8\baselineskip}
    \end{minipage}
    \caption{\colsr{Effect of drug transport parameters on change in stress levels post-drug treatment.
    (a) shows the effect of varying diffusivities in the adventitia and degenerated media on the stress increase in the fibrous cap.
    The columns indicate the volume-weighted median values, whereas the whiskers indicate the 5th and 95th percentiles.
    (b) depicts the effect of varying drug reaction rates in the healthy and degenerated media on the change in component stress.
    The points indicate the median values, and the shaded area indicates the 5th and 95th percentiles.
    Reaction rates in the degenerated media assumed to be half that of healthy media.
    Stress increase drops moderately with decreasing diffusivities, and drops significantly with increasing reaction rates, indicating a strong dependence on drug availability in the medial layers}
    }
    \label{fig:transport_param_effects}
    \end{figure*}

    Since not all drug transport parameters for the various components are available in the literature, we first investigated the effect that variations in the unknown assumed parameters have on the stress increase in the fibrous cap post-drug treatment.
    \colsr{We varied the drug diffusivities in the adventitia and degenerated media~(Table~\ref{tab:transport_parameters}) by one order of magnitude in both directions to obtain in total 9 combinations of parameter values.}
    Fig.~\ref{fig:transport_param_effects}(a) shows the \colsr{median and the 5th--95th percentile range of the drug-induced} stress increase in the fibrous cap for these combinations.
    \colsr{Increasing either diffusivity beyond the reference values does not strongly change the increase in the fibrous cap~(median $28.9$ to $29.4$~kPa), whereas decreasing them lowers it to a median of $20.6$~kPa~(5th--95th percentile: $8.1$--$27.0$~kPa) for the lowest combination.
    As less drug reaches the degenerated media, its relaxation is reduced~(median stress change of $-12.2$~kPa instead of $-23.0$~kPa), hence less load is shifted over to the fibrous cap.}

    Furthermore, the drug reaction rates in the healthy and degenerated media were assumed to be zero so far.
    Direct experimental measurements for these binding rates are unfortunately not available in the literature.
    Several studies however have used the second Damk{\"o}hler number $Da_{2} = K L^{2}/D$ ($L$ being the characteristic length, here 1 mm) to estimate these rates~\cite{JaiMcgPonZho:2022:tmd,McgWheMckMcc:2015:das,QiuNeo:2023:tpw}.
    It represents how dominant reaction effects are compared to diffusive effects: $Da_{2} << 1$ means diffusion dominated, whereas $Da_{2} \approx 1$ indicates an appreciable influence of reaction effects and $Da_{2} >> 1$ suggests a dominance of reactive effects.
    Given the unique pharmacokinetic properties of amlodipine, which exhibits not only a long half-life of 35--50 hours but also a more gradual onset, taking several hours to reach maximum effect and also slow offset after washout~\cite{BurMoi:1994:upp,NayGua:1991:abt}, it appears reasonable to assume slow reaction kinetics, rendering diffusion-dominated transport likely.
    Thus, to study this further, we varied the order of reaction rate in the healthy media between $3\cdot 10^{-6} \, \text{s}^{-1}$ to $3\cdot 10^{-1} \, \text{s}^{-1}$($Da_{2}$ range: $10^{-4}$--$10^{1}$) while maintaining all other parameters as indicated in \colsr{Tables~\ref{tab:material_parameters} and \ref{tab:transport_parameters}.}
    We assumed the degenerated media reaction rate to be half that of healthy media in each case.
    The effect these variations have on post-drug fibrous cap
    median stress increase is shown in Fig.~\ref{fig:transport_param_effects}(b), where a reduction in stress increase can be observed with increasing reaction rates, especially above $10^{-2} \, \text{s}^{-1}$~($Da_{2}\approx 1.43$) where reaction dominates and drug availability is significantly impacted, \colsr{with the median stress increase dropping to $21.8$~kPa at $3\cdot 10^{-2} \, \text{s}^{-1}$ and to $3.5$~kPa at $3\cdot 10^{-1} \, \text{s}^{-1}$.}
    \colsr{The simulations of the main text therefore use zero reaction rates~(Table~\ref{tab:transport_parameters}), since the fibrous cap stress increase for reaction rates of up to $10^{-3} \, \text{s}^{-1}$~($Da_{2}\approx 0.143$) does not differ significantly from the zero reaction rate case.}

\begin{figure}[t]
        \centering
        \includegraphics[width=0.45\textwidth]{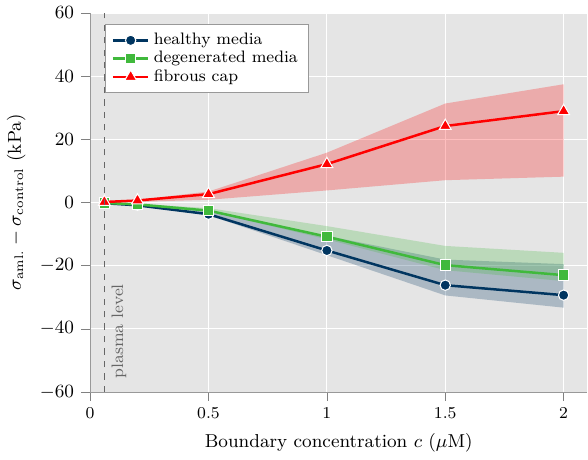}
        \caption{\colsr{Effect of the amlodipine concentration at the intimal and adventitial surfaces on the drug-induced change in von Mises stress in the external iliac artery~($\sigma_\text{aml.} - \sigma_\text{control}$ at $t=1500 \, \text{s}$).
        The lines indicate the volume-weighted medians and the shaded areas the 5th to 95th percentiles of each component; the dashed line marks the plasma level of $0.06 \, \mathrm{\mu M}$.
        The stress changes follow the sigmoidal drug-tissue interaction and are close to saturation at the assumed concentration of $2 \, \mathrm{\mu M}$~(fibrous cap: median $0.1$~kPa at the plasma level and $28.9$~kPa at $2 \, \mathrm{\mu M}$)} \label{fig:dose_effects}}
    \end{figure}

    \colsr{Besides the transport parameters, the drug availability is determined by the amlodipine concentration at the intimal and adventitial surfaces.
    In the simulations of the main text, it was set to $2\,\mathrm{\mu M}$~(Section~\ref{sec:boundary_conditions}), well above the plasma level of about $0.06\,\mathrm{\mu M}$.
    We therefore repeated the drug simulation with boundary concentrations of $0.06$, $0.2$, $0.5$, $1$, and $1.5\,\mathrm{\mu M}$.
    In all cases, the concentration in the arterial wall reaches the boundary value by $t = 1500 \, \text{s}$.
    Fig.~\ref{fig:dose_effects} shows that the stress changes follow the sigmoidal drug-tissue interaction~(Eq.~\eqref{eq:gamma}).
    At the plasma level, amlodipine has practically no mechanical effect~(median stress increase in the fibrous cap of $0.1$~kPa), and the increase remains small up to $0.5\,\mathrm{\mu M}$~(median $2.6$~kPa).
    Between $0.5$ and $1.5\,\mathrm{\mu M}$,
    it rises steeply to a median of $24.2$~kPa, whereas a further increase to $2\,\mathrm{\mu M}$ adds only $4.7$~kPa~(median $28.9$~kPa). The assumed concentration thus lies close to saturation.
    The drug-induced stress increase in the fibrous cap thus requires tissue concentrations well above the plasma level, such as those resulting from the tissue accumulation reported by \citet{JepBruNie:1998:ade}.}

\begin{figure*}[t]
    \centering
    \begin{minipage}[b]{0.40\linewidth}
        \vspace{0pt}
        \includegraphics[width=\linewidth]{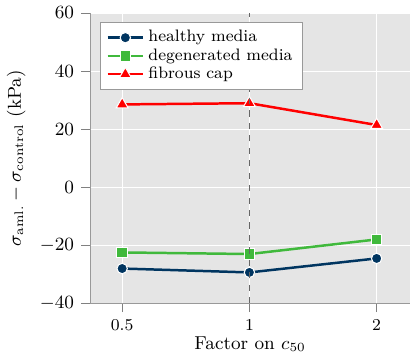}
        \par\vspace{-1.6\baselineskip}\noindent\llap{\footnotesize(a)\hspace{-0.6em}}
        \vspace{0.8\baselineskip}
    \end{minipage}
    \hspace{0.05\linewidth}
    \begin{minipage}[b]{0.40\linewidth}
        \vspace{0pt}
        \includegraphics[width=\linewidth]{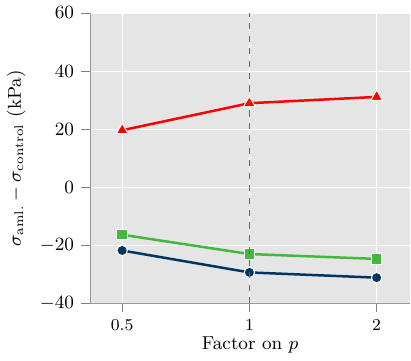}
        \par\vspace{-1.6\baselineskip}\noindent\llap{\footnotesize(b)\hspace{-0.6em}}
        \vspace{0.8\baselineskip}
    \end{minipage}
        \caption{\colsr{Effect of the drug-tissue interaction parameters on the drug-induced change in von Mises stress in the external iliac artery~($\sigma_\text{aml.} - \sigma_\text{control}$ at $t=1500 \, \text{s}$, boundary concentration $2 \, \mathrm{\mu M}$).
        (a) $c_{50}$ and (b) $p$ are halved and doubled, each with the other parameter kept at its fitted value; the dashed line marks the fitted values~(Table~\ref{tab:drug_interaction_parameters}).
        The lines connect the volume-weighted medians of each component.
        A less sensitive or flatter drug response~(doubled $c_{50}$ or halved $p$) reduces the stress changes, whereas a more sensitive or steeper response hardly changes them} \label{fig:drug_response}}
    \end{figure*}

    \colsr{Finally, the response of the SMCs to the drug is governed by the parameters $c_{50}$ and $p$ of the drug-tissue interaction~(Table~\ref{tab:drug_interaction_parameters}), which were fitted to data of rat aortic SMCs~\cite{WanJiaLiLai:2008:eos} and may differ in human atherosclerotic arteries.
    We therefore halved and doubled each of them at the boundary concentration of $2\,\mathrm{\mu M}$, one at a time with the other kept at its fitted value~(Fig.~\ref{fig:drug_response}).
    A more sensitive or steeper drug response~(halved $c_{50}$ or doubled $p$) changes the median stress increase in the fibrous cap only slightly, to $28.5$~kPa and $31.1$~kPa, respectively, since the response is already close to saturation.
    A less sensitive or flatter response~(doubled $c_{50}$ or halved $p$) lowers it to $21.4$~kPa and $19.6$~kPa, respectively, together with a smaller stress reduction in the healthy media~(median $-24.5$~kPa and $-21.8$~kPa instead of $-29.3$~kPa).
    Hence, the drug-induced shift of load from the medial layers to the fibrous cap persists for all considered drug-tissue interaction parameters.}
    \subsection{\colsr{Effect of axial stretch}}\label{sec:axial_stretch}
    \colsr{Since arteries may experience axial stretch,
    we repeated the simulations of the external iliac artery while additionally including axial stretches prior to the reorientation-growth-SMC protocol.
    These stretches were chosen to reflect a moderate axial strain of 5\% and a larger one of 10\% and were prescribed by axially fixing the bottom face and applying an axial displacement to the top face together with the luminal pressure during the load ramp, while all other boundary conditions, parameters and the loading protocol remained unchanged.
    \begin{figure*}[t]
    \centering
    \begin{minipage}[b]{0.47\linewidth}
        \vspace{0pt}
        \includegraphics[width=\linewidth]{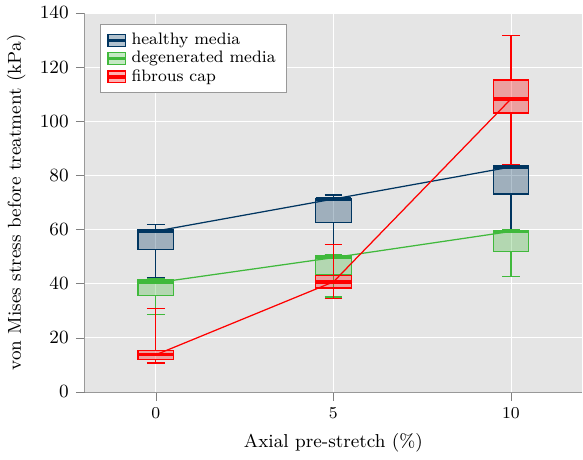}
        \par\vspace{-1.6\baselineskip}\noindent\llap{\footnotesize(a)\hspace{-0.6em}}
        \vspace{0.8\baselineskip}
    \end{minipage}%
    \hfill
    \begin{minipage}[b]{0.47\linewidth}
        \vspace{0pt}
        \includegraphics[width=\linewidth]{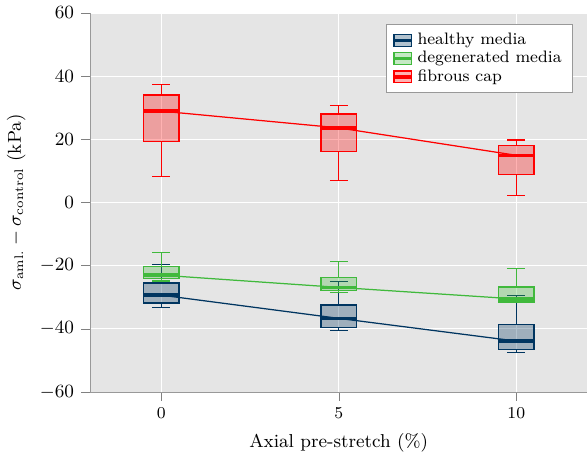}
        \par\vspace{-1.6\baselineskip}\noindent\llap{\footnotesize(b)\hspace{-0.6em}}
        \vspace{0.8\baselineskip}
    \end{minipage}
        \caption{\colsr{Effect of axial pre-stretch on the von Mises stress in the external iliac artery, for 0\%, 5\%, and 10\% axial pre-stretch.
        (a) shows the von Mises stress before drug treatment~($t=860 \, \text{s}$), and (b) the difference in von Mises stress between the simulations with and without drug effects~($\sigma_\text{aml.} - \sigma_\text{control}$) at $t=1500 \, \text{s}$.
        Each box spans the 25th to 75th percentile of the component's volume, with the median indicated by a thick horizontal line, whereas the whiskers represent the 5th and 95th percentiles, so $90\%$ of the component's volume lies between them.
        With increasing axial pre-stretch, the stresses before treatment increase in all components, most strongly in the fibrous cap~(change of the median from $14$ to $108$~kPa).
        The drug-induced stress reduction in the healthy media becomes larger~(median $-29$ to $-44$~kPa), whereas the stress increase in the fibrous cap becomes smaller~(median $29$ to $15$~kPa).
        However, as shown in Table~\ref{table:axial_stretch_components}, the circumferential stress increase only becomes smaller by approx. 20\%.} \label{fig:axial_stretch}}
    \end{figure*}
    \begin{table}[t]
        \caption{\colsr{Difference in circumferential ($\sigma_{\theta\theta}$) and axial ($\sigma_{zz}$) stress between the drug and control cases ($\sigma_\text{aml.} - \sigma_\text{control}$, median in kPa) for 0, 5 and 10\% axial pre-stretch, and the cross-sectional area $A$~(mm$^2$) of each component \label{table:axial_stretch_components}}}
        \centering
        \setlength{\tabcolsep}{4pt}
        \begin{tabular}{rcccc} \toprule
            & \begin{tabular}{@{}c@{}}Healthy\\ media\end{tabular} & \begin{tabular}{@{}c@{}}Degen.\\ media\end{tabular} & \begin{tabular}{@{}c@{}}Fibrous\\ cap\end{tabular} & Adventitia \\ \midrule
            $A$ & $1.88$ & $1.77$ & $0.51$ & $5.90$ \\ \midrule
            \multicolumn{5}{l}{Circumferential} \\
            0\% & $-3.2$ & $-7.5$ & $20.9$ & $1.6$ \\
            5\% & $-3.6$ & $-6.9$ & $17.9$ & $1.8$ \\
            10\% & $-4.4$ & $-6.3$ & $16.6$ & $2.1$ \\ \midrule
            \multicolumn{5}{l}{Axial} \\
            0\% & $-37.9$ & $-25.9$ & $34.0$ & $0.5$ \\
            5\% & $-43.6$ & $-29.4$ & $29.1$ & $0.7$ \\
            10\% & $-49.3$ & $-32.7$ & $21.6$ & $1.2$ \\ \bottomrule
        \end{tabular}
    \end{table}
    Fig.~\ref{fig:axial_stretch}(a) shows the median, 5th and 95th percentiles of the von Mises stresses before drug treatment for the three cases.
    The axial pre-stretch increases the von Mises stresses in all plaque components.
    Fig.~\ref{fig:axial_stretch}(b) depicts the difference in von Mises stresses between the drug and control cases in the medial layers and the fibrous cap.
    With pre-stretch, too, amlodipine reduces the stresses in the medial layers and increases those in the fibrous cap.
    The stress reduction in the medial layers increases with increasing axial pre-stretch, from a volume-weighted median value of $-29.3$~kPa at 0\% to $-36.7$~kPa at 5\% and $-43.9$~kPa at 10\% in the healthy media.
    A similar trend is also seen in the degenerated media, where the values change from $-23.0$~kPa at 0\% to $-26.9$~kPa at 5\% and $-30.5$~kPa at 10\% axial pre-stretch.
    However, the median fibrous cap stress increase also drops from $28.9$~kPa~(0\% case) to $23.6$~kPa and $14.9$~kPa for the 5\% and 10\% cases, respectively.
    This change is rather counter-intuitive, as we would expect that an increased unloading of the medial layers would shift more load onto the fibrous cap.
    Therefore, instead of just analyzing the von Mises stress, individual stress components are investigated and the difference in axial and circumferential components of the Cauchy stress tensor between the drug and control runs is shown in Table~\ref{table:axial_stretch_components}.
    As can be seen, the circumferential stress gain decreases with axial pre-stretch much less (only by approx. 20\%) than the von Mises stress gain.
    This is also in line with the changes in the other components.
    A majority of the increased medial unloading seen in Fig.~\ref{fig:axial_stretch}(b) can be attributed to axial stress.
    This reduced axial loading does not need to be shifted to the other layers, as it simply results in a reduced reaction force axially.
    However, the circumferential loading is primarily a result of the luminal pressure, which does not change.
    From Table~\ref{table:axial_stretch_components}, it is evident that the change in circumferential stress in the adventitia also increases with pre-stretch, which in combination with its relatively large cross-sectional area contributes to it carrying a higher proportion of the circumferential load.
    This is likely a result of adventitial collagen fibers being recruited at lower circumferential stretches due to axial pre-stretching.
    For the forthcoming simulations, we adopt the 0\% pre-stretch since it produces the largest drug-induced cap stress and can be considered as the worst-case scenario.}

\end{appendices}

\clearpage  %
\bibliographystyle{references_style}
\bibliography{references_merged}

\end{document}